\documentclass[11pt]{article}

\usepackage[letterpaper,left=1in,top=1in,right=1in,bottom=1.2in,nohead]{geometry}

\usepackage[pdftex, colorlinks=true, urlcolor=MyDarkBlue, citecolor=MyDarkRed, linkcolor=MyDarkGreen]{hyperref}
\usepackage[utf8]{inputenc}
\usepackage{amsmath,amsthm,epsf,amscd,amssymb,verbatim}
\usepackage[square,sort,comma,numbers]{natbib}
\usepackage{graphicx}
\usepackage{mathrsfs}
\usepackage{color}
\usepackage{titlesec}
\usepackage{enumerate}
\usepackage[shortlabels]{enumitem}
\usepackage{datetime}
\usepackage[normalem]{ulem}
\usepackage[font=footnotesize]{caption}
\usepackage{abstract}

\usepackage{chemarr}
\usepackage{extarrows}
\usepackage[all]{xy}
\usepackage{enumerate}
\usepackage{tabularx}

\usepackage[outline]{contour}

\usepackage{wrapfig}
\usepackage{fancyhdr}

\usepackage[low-sup]{subdepth}

\usepackage[framemethod=TikZ]{mdframed}

\usepackage[caption=false,labelformat=simple]{subfig}

\usepackage{times}
\usepackage{relsize}

\newcommand{\hide}[1]{}

\usepackage{color}
\definecolor{MyDarkRed}{rgb}{0.5,0,0.1}
\definecolor{MyDarkBlue}{rgb}{0.1,0.1,0.5}
\definecolor{MyDarkGreen}{rgb}{0.1,0.5,0.1}
\definecolor{MyRed}{rgb}{1.0,0,0}
\definecolor{MyBlue}{rgb}{0,0,1.0}
\definecolor{MyGreen}{rgb}{0,0.8,0}
\definecolor{lightgray}{rgb}{0.96,0.96,0.96}
\definecolor{darkgray}{rgb}{0.4,0.4,0.4}

\newtheorem{lemma}{Lemma}
\newtheorem{cor}{Corollary}
\newtheorem{prop}{Proposition}

\theoremstyle{definition}
\newtheorem{definition}{Definition}
\theoremstyle{remark}

\newtheorem{mynote}{Note}

\newcommand{\spn}{{~\mathrm{span}}}

\newcommand{\wipe}[1]{}

\newcommand{\conjsym}{{\dagger}}

\newcommand{\corrected}[1]{{#1}} 

\definecolor{MyRed}{rgb}{1.0,0,0}
\definecolor{MyDarkRed}{rgb}{0.6,0,0}
\definecolor{MyBlue}{rgb}{0,0,1.0}
\definecolor{MyDarkBlue}{rgb}{0,0,0.6}
\definecolor{MyPurple}{rgb}{0.5,0,0.5}

\ifx\hidespecial\undefined
	\newcommand{\todo}[1]{ {\color{MyRed}\textbf{TODO:} #1} }
	\newcommand{\todoassign}[2]{ {\color{red}\textbf{{#1}:} {#2}} }
	
\else
	\newcommand{\todo}[1]{}
	\newcommand{\todoassign}[2]{}
	
\fi

\ifx\hideproof\undefined
	\usepackage{framed}
	\definecolor{leftBarGray}{rgb}{0.8,0.8,0.8}
	\newenvironment{myleftbar}{%
	\MakeFramed {\advance\hsize-\width \FrameRestore}}%
	{\endMakeFramed}
	
\else
	\usepackage{environ}
	\NewEnviron{quoteproof}{}
\fi

\usepackage{textgreek}

\titleformat*{\section}{\Large\bfseries}

\titleformat*{\subsection}{\large\bfseries}

\titleformat*{\subsubsection}{\normalfont\bfseries}

\newcommand{\pubnote}[1]{}

\newlength{\sectiongaps}
\newlength{\sectiongapsafter}
\newlength{\subsubsectiongaps}
\allowdisplaybreaks

\newcommand{\subspacetext}{{sp}}
\newcommand{\couplingtext}{{\mathcal{CP}}}
\newcommand{\externaldegtext}{{\mathcal{ED}}}
\newcommand{\maxextdeg}{{\mathcal{MED}}}

\newcommand\Item[1][]{%
	\ifx\relax#1\relax  \item \else \item[#1] \fi
	\abovedisplayskip=0pt\abovedisplayshortskip=0pt~\vspace*{-\baselineskip}}

\usepackage{etoolbox}

\mdfsetup{%
	middlelinewidth=2pt,
	roundcorner=10pt}

\definecolor{MyLightGray8}{rgb}{0.8,0.8,0.8} 
\definecolor{MyLightGray9}{rgb}{0.9,0.9,0.9} 
\definecolor{MyLightGray95}{rgb}{0.95,0.95,0.95} 

\definecolor{proofcolor}{rgb}{0.96,0.96,0.96}

\BeforeBeginEnvironment{proof}{%
	\begin{mdframed}[topline=false,
				rightline=false,
				bottomline=false,
				leftline=false,
				innerleftmargin=1ex,
				innerrightmargin=1ex,
				innertopmargin=1ex,
				innerbottommargin=1ex,
				backgroundcolor=proofcolor,
				]\hyphenchar\font=\defaulthyphenchar 
}
\AfterEndEnvironment{proof}{%
\end{mdframed}\hyphenchar\font=-1
}

\definecolor{lemmacolor}{rgb}{1.0,0.95,0.9}

\BeforeBeginEnvironment{lemma}{%
	\begin{mdframed}[topline=false,
		rightline=false,
		bottomline=false,
		leftline=false,
		innerleftmargin=1ex,
		innerrightmargin=1ex,
		innertopmargin=1ex,
		innerbottommargin=1ex,
		backgroundcolor=lemmacolor,
		]\hyphenchar\font=\defaulthyphenchar 
	}
	\AfterEndEnvironment{lemma}{%
	\end{mdframed}\hyphenchar\font=-1
}

\definecolor{corcolor}{rgb}{0.95,0.95,0.9}

\BeforeBeginEnvironment{cor}{%
	\begin{mdframed}[topline=false,
		rightline=false,
		bottomline=false,
		leftline=false,
		innerleftmargin=1ex,
		innerrightmargin=1ex,
		innertopmargin=1ex,
		innerbottommargin=1ex,
		backgroundcolor=corcolor,
		]\hyphenchar\font=\defaulthyphenchar 
	}
	\AfterEndEnvironment{cor}{%
	\end{mdframed}\hyphenchar\font=-1
}

\definecolor{propcolor}{rgb}{1.0,0.93,0.94}

\BeforeBeginEnvironment{prop}{%
	\begin{mdframed}[topline=false,
		rightline=false,
		bottomline=false,
		leftline=false,
		innerleftmargin=1ex,
		innerrightmargin=1ex,
		innertopmargin=1ex,
		innerbottommargin=1ex,
		backgroundcolor=propcolor,
		]\hyphenchar\font=\defaulthyphenchar 
	}
	\AfterEndEnvironment{prop}{%
	\end{mdframed}\hyphenchar\font=-1
}

\definecolor{defcolor}{rgb}{0.92,0.92,0.95}

\BeforeBeginEnvironment{definition}{%
	\begin{mdframed}[topline=false,
		rightline=false,
		bottomline=false,
		leftline=false,
		innerleftmargin=1ex,
		innerrightmargin=1ex,
		innertopmargin=1ex,
		innerbottommargin=1ex,
		backgroundcolor=defcolor,
		]\hyphenchar\font=\defaulthyphenchar 
	}
	\AfterEndEnvironment{definition}{%
	\end{mdframed}\hyphenchar\font=-1
}

\title{On Some Bounds on the Perturbation of Invariant Subspaces of Normal Matrices with Application to a Graph Connection Problem}
\date{}
\author{Subhrajit Bhattacharya \thanks{Lehigh University, U.S.A.  ~$\bullet$~
		\href{https://www.lehigh.edu/\%7Esub216/}{https://www.lehigh.edu/$\sim$sub216/} ~$\bullet$~
		email: \href{mailto:sub216@lehigh.edu}{\texttt{sub216@lehigh.edu}}.}}

\begin{document}

\maketitle

\begin{abstract}
	\noindent
	We provide upper bounds on the perturbation of invariant subspaces of normal matrices 
	\corrected{measured using a metric on the space of vector subspaces of $\mathbb{C}^n$
	in terms of the spectrum of both the unperturbed \& perturbed matrices, as well as, spectrum of the unperturbed matrix only. The results presented give tighter bounds than the Davis-Khan $\sin\Theta$ theorem.}
%
%
	We apply the result to a graph perturbation problem.
\end{abstract}

\tableofcontents

\section{Introduction}

Classical results on perturbation of invariant subspaces of a matrix 
\corrected{
usually take one of the two forms:
(1) perturbation measured in terms of a natural metric in the space of vector subspaces (usually expressed as the sine of the angle between subspaces), with upper bound described in terms of the perturbation in the matrices as well as the spectra of both the unperturbed and perturbed matrices (for example, the Davis-Kahan $\sin\Theta$ Theorem~\cite{10.1137/0707001} -- see Section VIII.3 of~\cite{bhatia1996matrix} where generalization of this theorem is given for normal matrices.); or,
(2) perturbation measured in terms of bounds on norms of matrices that relate an invarient subspace with its perturbation in a more complex manner (which, in general, is not a natural metric in the space of vector subspaces) 
although the upper bound is based on the spectrum of the unperturbed matrix only (see, for example, \cite{10.2307/2028728,10.1137/0707001} or Chapter V of \cite{Stewart90}).}

\corrected{
In this paper we first derive an upper-bound reminiscent of the Davis-Kahan $\sin\Theta$ Theorem, but generalized for normal matrices and with modestly tighter bound (Proposition~\ref{prop:full-main}).
Then we use some geometric methods to derive a bound on perturbation measured in terms of a natural metric in the space of subspaces, but with upper-bounds in terms of spectrum of the unperturbed matrix only (Proposition~\ref{prop:subspace-purturbation-hat-tilde-free}) when the spectrum is well-clustered (a relation formally described as ``separation-preserving perturbation'').
In this later case our proposed result also allows easy identification of the perturbed invariant subspace (Lemma~\ref{lemma:sep-preserving-M-tilde}).
}

\begin{definition}[Notations] \label{def:notations}
\corrected{ 
Throughout the paper we assume $M, \widetilde{M} \in \mathbb{C}^{n\times n}$ to be normal matrices unless specified otherwise, and by ``eigenvectors'' we will refer to their right eigenvectors.
The eigenvalues (not necessarily distinct) and corresponding unit eigenvectors (for degenerate eigenspaces, any orthonormal basis thereof) of $M$ be $\lambda_j$ and $\mathbf{u}_j$ for $j=1,2,\cdots,n$.
Likewise, the eigenvalues and corresponding unit eigenvectors of $\widetilde{M}$ be $\widetilde{\lambda}_j$ and $\widetilde{\mathbf{u}}_j$ for $j=1,2,\cdots,n$.
We will usually consider the eigenvectors to be column vectors in $\mathbb{C}^{n\times 1}$.
Let $U = [\mathbf{u}_1,\mathbf{u}_2,\cdots,\mathbf{u}_n]$ and $\widetilde{U} = [\widetilde{\mathbf{u}}_1,\widetilde{\mathbf{u}}_2,\cdots,\widetilde{\mathbf{u}}_n]$ be the unitary matrices that diagonalize $M$ and $\widetilde{M}$ respectively.
For notational convenience, define $N = \{1,2,\cdots,n\}$.

As a convention, we choose primed lower-case Latin letters to index variables (eigenvalues or eigenvectors) with tilde on them.
Given a set $S \subseteq N$, we define the set $\mathbf{u}_S = \{\mathbf{u}_j \,|\, j\in S\}$. Likewise $\widetilde{\mathbf{u}}_S = \{\widetilde{\mathbf{u}}_{j'} \,|\, {j'}\in S\}$. 
Define the multi-sets $\lambda_S = \{\lambda_j \,|\, j\in S\}$ and $\widetilde{\lambda}_S = \{\widetilde{\lambda}_{j'} \,|\, {j'}\in S\}$ (by asserting that these are multi-sets, we allow multiplicity in the values, thus ensuring these sets have the same number of elements as $S$).
%
%
We also define the \emph{complement} of $S$ as $S^c = N-S$.
}
\end{definition}

\noindent
The outline of the paper \corrected{is as follows}:
\begin{enumerate}
	\item \corrected{In in Section~\ref{sec:d-subspace} we describe} a natural metric, \corrected{$d_\mathrm{\subspacetext}$,} on $Gr(q,\mathbb{C}^n)$ (the space of $q$-dimensional complex vector subspaces of $\mathbb{C}^n$) to measure perturbation of invariant subspaces of $n\times n$ normal matrices. \corrected{This metric is equivalent to the Frobenius norm of the $\sin\Theta$ matrix between subspaces of $\mathbb{C}^n$.}
	\item Some geometry lemmas are proven in Section~\ref{sec:geomrtry}, \corrected{which are used in Section~\ref{sec:well-clusterted-spectrum} for 
	deriving bounds on the perturbation of invariant subspaces in terms of the spectrum of the unperturbed matrix only (when the spectrum is well-clustered).}
	\item 
		\corrected{In Section~\ref{sec:partitioned-spectrum} we describe} upper-bound on the distance between 
		
	\noindent  \hspace{-0.6em}
	\begin{tabular}{p{0.65\columnwidth}p{0.35\columnwidth}} \vspace{-1.2em}
		\corrected{
		 invariant subspaces in terms of the spectrum of both the unperturbed and perturbed matrices.}
		\corrected{Some of these results give improvements on the Davis-Kahan $\sin\Theta$ Theorem for normal matrices (although the Davis-Kahan $\sin\Theta$ is usually stated for Hermitian matrices, there exists generalizations of the theorem for normal matrices -- see Section VIII.3 of~\cite{bhatia1996matrix}).}
			As an example, 
		\corrected{for any $J,\widetilde{J} \subseteq N$, with $|J| = |\widetilde{J}| = q$,
			\textbf{Proposition~\ref{prop:full-main}} states,}
		&
		\vspace{-4em}
		{\includegraphics[width=0.26\textwidth, clip=true, trim=0 0 0 0]{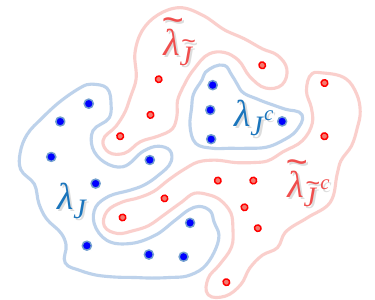}}
	\end{tabular}
%
	\corrected{
		\begin{align} 
			d_\mathrm{\subspacetext} \left( \spn(\mathbf{u}_J) , \spn(\widetilde{\mathbf{u}}_{\widetilde{J}}) \right) ~~\leq~~ \sqrt{ \frac{1}{q}\, \mathlarger{\mathlarger{\sum}}_{j\in J} ~\frac{ 
					\left\| (\widetilde{M} - M) \mathbf{u}_{{j}} \right\|^2_2 ~-~ \kappa_j \,\min_{j'\in \widetilde{J}}\left|\widetilde{\lambda}_{{j'}} - \lambda_j\right|^2 
				}{ 
					\min_{j'\in {\widetilde{J}}^c} \left|\widetilde{\lambda}_{{j'}} - \lambda_j\right|^2 - \kappa_j \,\min_{j'\in \widetilde{J}}\left|\widetilde{\lambda}_{{j'}} - \lambda_j\right|^2 
			} }\nonumber
		\end{align}
	with
%
		\begin{equation}\nonumber
		\kappa_j = \left\{ \begin{array}{ll}
		0, & \text{if ~$\| (\widetilde{M} - M) \mathbf{u}_{{j}} \|_2 \geq \min_{j'\in {\widetilde{J}}^c} |\widetilde{\lambda}_{{j'}} - \lambda_j|$} \\
		1, & \text{if ~$\| (\widetilde{M} - M) \mathbf{u}_{{j}} \|_2 < \min_{j'\in {\widetilde{J}}^c} |\widetilde{\lambda}_{{j'}} - \lambda_j|$}
		\end{array}\right.
		\end{equation}
		This is a tighter upper bound than the Davis-Kahan $\sin\Theta$ Theorem, and as a consequence
		leads to the rediscovery of a couple of slight variations on the Davis-Kahan $\sin\Theta$ Theorem 
		in Corollary~\ref{cor:kappa-zero}, where, as an example, one result states
		\[\displaystyle d_\mathrm{\subspacetext} \left( \spn(\mathbf{u}_J) , \spn(\widetilde{\mathbf{u}}_{\widetilde{J}}) \right) 
			~\leq~ 
			\frac{  \displaystyle \min \left( 1, \sqrt{\frac{n-q}{q}} \right) }{ \displaystyle \max \left( \mathrm{sep}\left( \lambda_J, \widetilde{\lambda}_{{\widetilde{J}}^c} \right), ~\mathrm{sep}\left( \lambda_{J^c}, \widetilde{\lambda}_{\widetilde{J}} \right) \right) } \left\| \widetilde{M} - M \right\|_2 \]
		where, $\mathrm{sep}(P,Q) = \min_{p\in P,\atop q\in Q}|p-q|$ simply measures the min-min distance between the sets
		(this is unlike the Davis-Kahan $\sin\Theta$ Theorem generalized for normal matrices, where it is necessary to find a `strip' or `annulus' of width $\delta$ separating $\lambda_J$ and $\widetilde{\lambda}_{{\widetilde{J}}^c}$ -- see Theorem VIII.3.1 of~\cite{bhatia1996matrix}).
	}
%
%
	\item 
	\corrected{The next set of main results of this paper appear in Section~\ref{sec:well-clusterted-spectrum}, which}
	 formalizes the notion of well-clustered spectrum in \textbf{Lemma~\ref{lemma:sep-preserving-M-tilde}}, followed by \textbf{Proposition~\ref{prop:subspace-purturbation-hat-tilde-free}} that provides the upper bound on the perturbation of an invariant subspace in terms of the spectrum of the unperturbed matrix only.
	These results rely on the geometry lemmas from Section~\ref{sec:geomrtry}.
	As an example, one of the results of Proposition~\ref{prop:subspace-purturbation-hat-tilde-free} \corrected{states that if $\|\widetilde{M} - M\|_2 < \frac{1}{2} \mathrm{sep}(\lambda_J, \lambda_{J^c})$, then}
	\begin{align*}
	d_\mathrm{\subspacetext} \left( \spn(\mathbf{u}_J) , \spn(\widetilde{\mathbf{u}}_{\widehat{J}}) \right)
	& ~~\leq~~
	\frac{1}{\sqrt{q}}\, \min \left(
	\sqrt{\mathlarger{\mathlarger{\sum}}_{j\in J} \left(\frac{ 
		\left\| (\widetilde{M} - M) \mathbf{u}_{{j}} \right\|_2 
	}{ 
		\min_{k\in {J^c}} |{\lambda}_{k} - \lambda_j| ~-~ \|\widetilde{M} - M\|_2  
	} \right)^{\!\!2}} ~,~ 
	\right. \nonumber \\ & \qquad\qquad \qquad\qquad  \left.
	\sqrt{\mathlarger{\mathlarger{\sum}}_{j\in {J^c}} \left( \frac{ 
		\left\| (\widetilde{M} - M) \mathbf{u}_{{j}} \right\|_2 
	}{ 
		\min_{k\in J} |{\lambda}_{k} - \lambda_j| ~-~ \|\widetilde{M} - M\|_2  
	} \right)^{\!\!2} } \right)
	\\
	& \quad ~~\leq~~
	\min\left( 1, \sqrt{\frac{n-q}{q}} \right) ~\frac{ 
		\left\| \widetilde{M} - M \right\|_2 
	}{ 
		\mathrm{sep}(\lambda_J, \lambda_{J^c})  ~-~ \|\widetilde{M} - M\|_2  
	} 
	\end{align*}
	where 
	\corrected{$\widehat{J} = \{j' \,|\, \min_{j\in N} |\widetilde{\lambda}_{j'} - \lambda_j| = \min_{j\in J} |\widetilde{\lambda}_{j'} - \lambda_j| \}$ is the set of indices corresponding to the eigenvalues of $\widetilde{M}$ that are closer to $\lambda_J$ than to $\lambda_{J^c}$.}
	\item Section~\ref{sec:graph} demonstrates an application to the perturbation of null-space of a matrix in context of a graph perturbation problem.
\end{enumerate}

\section{Preliminaries}

\subsection{A Metric on $Gr(q,\mathbb{C}^n)$} \label{sec:d-subspace}


\begin{definition}[Subspace Distance] \label{def:d-subspace}
	Suppose $X,Y \subseteq \mathbb{C}^n$ are $q$-dimensional vector sub-spaces of $\mathbb{C}^n$. \newline Let $\{\mathbf{x}_j\}_{j=1,2,\cdots,q}$ and $\{\mathbf{y}_j\}_{j=1,2,\cdots,q}$ be orthonormal basis on $X$ and $Y$.
	The \emph{\underline{s}ubs\underline{p}ace distance} between $X$ and $Y$ is defined as
	\begin{equation}
		d_\text{\subspacetext} (X, Y) ~=~ \frac{1}{\sqrt{2q}} \| \mathbf{X}\mathbf{X}^\conjsym - \mathbf{Y}\mathbf{Y}^\conjsym \|_F
	\end{equation}
	where, \begin{equation}\begin{array}{ll}
	\mathbf{X} = [\mathbf{x}_1, \mathbf{x}_2, \cdots, \mathbf{x}_q] ~~~~\text{and} &
	\mathbf{Y} = [\mathbf{y}_1, \mathbf{y}_2, \cdots, \mathbf{y}_q] 
	\end{array} \label{eq:XY-matrix-def} \end{equation}
	are the $n\times q$ matrices in which the columns represent the unit vectors $\{\mathbf{x}_j\}_{j=1,2,\cdots,q}$ and $\{\mathbf{y}_j\}_{j=1,2,\cdots,q}$.
\end{definition}

Note that the matrices $\mathbf{X}\mathbf{X}^\conjsym$ and $\mathbf{Y}\mathbf{Y}^\conjsym$ are the projection operators on $X$ and $Y$ respectively.
The space of difference of such projection operators is well-studied in literature (see~\cite{ANDRUCHOW20141634,golub1996matrix} for example), and the norms of such differences have been used as metric on $Gr(q,\mathbb{C}^n)$ (see~\cite{damle2020uniform} for example).
\corrected{In fact this metric is equivalent to the Frobenius norm of the $\sin\Theta$ matrix between subspaces of $\mathbb{C}^n$ that is used for measuring perturbation of invariant subspaces in context of the Davis-Kahan $\sin\Theta$ Theorem.}
We choose 
the Frobenius norm \corrected{for measuring the distance between the projection operators, and use a scaling factor of $\frac{1}{\sqrt{2q}}$ for convenience and some additional properties of the metric. The following lemmas outline some elementary and mostly standard properties of this metric.} 

Let $X^\perp$ and $Y^\perp$ are orthogonal complements of $X$ and $Y$ respectively in $\mathbb{C}^n$.
Let $\{\mathbf{x}_j\}_{j=q+1,q+2,\cdots,n}$ and $\{\mathbf{y}_k\}_{k=q+1,q+2,\cdots,n}$ be orthonormal basis for $X^\perp$ and $Y^\perp$ respectively.
Define 
 \begin{equation}\begin{array}{ll}
\mathbf{X}^\perp = [\mathbf{x}_{q+1}, \mathbf{x}_{q+2}, \cdots, \mathbf{x}_n] ~~~~\text{and} &
\mathbf{Y}^\perp = [\mathbf{y}_{q+1}, \mathbf{y}_{q+2}, \cdots, \mathbf{y}_n]
\end{array} \label{eq:XperpYperp-matrix-def} \end{equation}

\begin{lemma}[Equivalent Forms of $d_\text{\subspacetext}$] ~ \label{lemma:d-subspace-equivalent-forms}
\begin{enumerate}
	\item $\displaystyle d_\text{\normalfont{\subspacetext}} (X, Y) ~~=~~ \sqrt{1 - \frac{1}{q} \|\mathbf{X}^\conjsym \mathbf{Y}\|_F^2} ~~=~~ \sqrt{1 - \frac{1}{q} \sum_{j=1}^q \sum_{k=1}^q \left| \mathbf{x}_j^\conjsym \mathbf{y}_{k} \right|^2}$
	\item $\displaystyle d_\text{\normalfont{\subspacetext}} (X, Y) ~~=~~ \sqrt{\frac{1}{q}  \|{\mathbf{X}^\perp}^\conjsym \,\mathbf{Y}\|_F^2} ~~=~~ \sqrt{\frac{1}{q}  \sum_{j=q+1}^n \sum_{k=1}^q \left| \mathbf{x}_j^\conjsym \mathbf{y}_{k} \right|^2}$
\end{enumerate}
\end{lemma}

\begin{proof} ~
	\begin{enumerate}
		\item In the following we use the definition $\|\mathbf{A}\|_F^2 = \mathrm{tr}(\mathbf{A}^\conjsym \mathbf{A})$ and the property that $\mathrm{tr}(\mathbf{A}\mathbf{B}) = \mathrm{tr}(\mathbf{B}\mathbf{A})$.
		\begin{align*}
		& \left( d_\text{\normalfont{\subspacetext}} (X, Y) \right)^2 \\
		& ~~=~~ \frac{1}{2q} \| \mathbf{X}\mathbf{X}^\conjsym - \mathbf{Y}\mathbf{Y}^\conjsym \|_F^2  \\
		& ~~=~~ \frac{1}{2q} \mathrm{tr} \left( \left( \mathbf{X}\mathbf{X}^\conjsym - \mathbf{Y}\mathbf{Y}^\conjsym \right)^\conjsym\left( \mathbf{X}\mathbf{X}^\conjsym - \mathbf{Y}\mathbf{Y}^\conjsym \right) \right) \\
		& ~~=~~ \frac{1}{2q} \mathrm{tr} (\mathbf{X}\mathbf{X}^\conjsym \mathbf{X}\mathbf{X}^\conjsym) + \mathrm{tr} (\mathbf{Y}\mathbf{Y}^\conjsym\mathbf{Y}\mathbf{Y}^\conjsym)
		- \mathrm{tr} (\mathbf{X}\mathbf{X}^\conjsym \mathbf{Y}\mathbf{Y}^\conjsym) - \mathrm{tr} ( \mathbf{Y}\mathbf{Y}^\conjsym \mathbf{X}\mathbf{X}^\conjsym) \\
		& ~~=~~ \frac{1}{2q} \left( \mathrm{tr} (\mathbf{X}\mathbf{X}^\conjsym) + \mathrm{tr} (\mathbf{Y}\mathbf{Y}^\conjsym) 
		- 2 \mathrm{tr} ( \mathbf{Y}^\conjsym \mathbf{X}\mathbf{X}^\conjsym \mathbf{Y}) \right) \qquad\text{\small (since $\mathbf{X}^\conjsym \mathbf{X} = \mathbf{Y}^\conjsym \mathbf{Y} = I$.)} \\
		& ~~=~~ 1 -  \frac{1}{q} \|\mathbf{X}^\conjsym \mathbf{Y}\|_F^2 \qquad\text{\small (since $\mathrm{tr} \left(  \mathbf{X} \mathbf{X}^\conjsym \right) = \mathrm{tr} \left( \mathbf{X}^\conjsym \mathbf{X} \right) = \sum_{j=1}^q \mathbf{x}_j^\conjsym \mathbf{x}_j = q$, and likewise for $\mathbf{Y}$.)} \\
		& ~~=~~ 1 - \frac{1}{q} \sum_{j=1}^q \sum_{k=1}^q \left| \mathbf{x}_j^\conjsym \mathbf{y}_{k} \right|^2
		\end{align*}
		
		\item Note that $[\mathbf{X}, \mathbf{X}^\perp]$ is a $n \times n$ unitary matrix with columns being the vectors of the orthonormal basis $\{\mathbf{x}_i\}_{i=1,2,\cdots,n}$. Thus, $[\mathbf{X}, \mathbf{X}^\perp] [\mathbf{X}, \mathbf{X}^\perp]^\conjsym = \mathbf{X} \mathbf{X}^\conjsym + \mathbf{X}^\perp {\mathbf{X}^\perp}^\conjsym = I$.
		Thus,
		\begin{align*}
		\left( d_\mathrm{\subspacetext} \left( X, Y \right) \right)^2 & ~=~ 1 - \frac{1}{q} \|\mathbf{X}^\conjsym \mathbf{Y}\|_F^2 \nonumber \\
		& ~=~ 1 - \frac{1}{q} \mathrm{tr} \left( \mathbf{Y}^\conjsym \mathbf{X} \mathbf{X}^\conjsym \mathbf{Y} \right) \nonumber \\
		& ~=~ 1 - \frac{1}{q} \mathrm{tr} \left( \mathbf{Y}^\conjsym \left(I - \mathbf{X}^\perp {\mathbf{X}^\perp}^\conjsym \right) \mathbf{Y} \right) \nonumber \\
		& ~=~ 1 - \frac{1}{q} \mathrm{tr} \left( \mathbf{Y}^\conjsym \mathbf{Y} \right) + \frac{1}{q} \mathrm{tr} \left( \mathbf{Y}^\conjsym \mathbf{X}^\perp {\mathbf{X}^\perp}^\conjsym \mathbf{Y} \right) \nonumber \\ 
		& ~=~ 1 - \frac{1}{q} q + \frac{1}{q} \mathrm{tr} \left( \mathbf{Y}^\conjsym \mathbf{X}^\perp {\mathbf{X}^\perp}^\conjsym \mathbf{Y} \right)  \nonumber \\
		& ~=~ \frac{1}{q} \|{\mathbf{X}^\perp}^\conjsym \mathbf{Y} \|_F^2 ~=~ \frac{1}{q} \sum_{j=q+1}^n \sum_{k=1}^q \left| \mathbf{x}_j^\conjsym \mathbf{y}_{k} \right|^2 
		\end{align*}
	\end{enumerate}
\end{proof}


\begin{lemma}[Properties of $d_\text{\subspacetext}$] ~ \label{lemma:d-subspace-properties}
	\begin{enumerate}
		\item The value of $d_\mathrm{\subspacetext} \left( X, Y \right)$ is independent of the choice of basis on $X$ or $Y$ (or the basis on $X^\perp$ or $Y^\perp$, if using the equivalent form in Lemma~\ref{lemma:d-subspace-equivalent-forms}.2).
		\item $d_\mathrm{\subspacetext}$ is a metric on $Gr(q, \mathbb{C}^n)$ (the space of $q$-dimensional complex subspaces of $\mathbb{C}^n$). 
		\item $ \sqrt{q} ~ d_\mathrm{\subspacetext} \left( X, Y \right) ~=~ \sqrt{n-q} ~ d_\mathrm{\subspacetext} \left( X^\perp, Y^\perp \right)  $
		\item $d_\mathrm{\subspacetext} \left( X, Y \right) \leq 1$, with equality holding iff $X$ and $Y$ are orthogonal subspaces (which is possible only if $q \leq n/2$).
	\end{enumerate}
\end{lemma}

\begin{proof} ~
	\begin{enumerate}
		\item Suppose $\{\mathbf{x}'_j\}_{j=1,2,\cdots,q}$ and $\{\mathbf{y}'_j\}_{j=1,2,\cdots,q}$ be a different set of orthonormal bases on $X$ and $Y$ respectively.
		Define $\mathbf{X}' = [\mathbf{x}'_1, \mathbf{x}'_2, \cdots, \mathbf{x}'_q],
		\mathbf{Y}' = [\mathbf{y}'_1, \mathbf{y}'_2, \cdots, \mathbf{y}'_q]$.
		Thus there exists $q\times q$ unitary matrices $R_X, R_Y \in U(q)$ such that $\mathbf{X} = \mathbf{X}' R_X$ and $\mathbf{Y} = \mathbf{Y}' R_Y$. Then,
		\begin{align*}
		& \left( d_\text{\normalfont{\subspacetext}} (X, Y) \right)^2 \\
		& ~~=~~ \frac{1}{2q} \| \mathbf{X}\mathbf{X}^\conjsym - \mathbf{Y}\mathbf{Y}^\conjsym \|_F^2  \\
		& ~~=~~ \frac{1}{2q} \| \left(\mathbf{X}' R_X\right) \left(\mathbf{X}' R_X\right)^\conjsym - \left(\mathbf{Y}' R_Y\right) \left(\mathbf{Y}' R_Y\right)^\conjsym \|_F^2  \\
		& ~~=~~ \frac{1}{2q} \| \mathbf{X}' {\mathbf{X}'}^\conjsym - \mathbf{Y}' {\mathbf{Y}'}^\conjsym \|_F^2
		\end{align*}
		For the equivalent form in Lemma~\ref{lemma:d-subspace-equivalent-forms}.2 we can use the orthonormal basis $\{\mathbf{x}'_j\}_{j=q+1,q+2,\cdots,n}$ and $\{\mathbf{y}'_k\}_{k=q+1,q+2,\cdots,n}$ for $X^\perp$ and $Y^\perp$ respectively and analogously derive at the equivalent form using the primed basis.
		
		\item \emph{Non-negativity} and \emph{symmetry} properties are obvious from the definition of $d_\text{\subspacetext}$. 
		
		If $X$ and $Y$ are the same subspaces, we can choose the same basis for them (since the value of $d_\text{\normalfont{\subspacetext}} (X, Y)$ is independent of the choice of basis on $X$ and $Y$), doing so makes it obvious that $d_\text{\normalfont{\subspacetext}} (X, Y) = 0$.
		
		\emph{Triangle inequality} holds due to the fact that Frobenius norm of difference of matrices is a metric on $\mathbb{C}^{n\times n}$.

		\item Note that $X^\perp$ and $Y^\perp$ are $(n-q)$-dimensional subspaces of $\mathbb{C}^n$. Furthermore, $X$ is the orthogonal complement of $X^\perp$. As a consequence, due to Lemma~\ref{lemma:d-subspace-equivalent-forms}.2.,
		\begin{align*}
		 \displaystyle d_\text{\normalfont{\subspacetext}} (X^\perp, Y^\perp)
		& ~=~~ \sqrt{ \frac{1}{n-q} } \|{\mathbf{X}}^\conjsym \,\mathbf{Y}^\perp\|_F \\
		& ~=~~ \sqrt{ \frac{1}{n-q} } \|{\mathbf{Y}^\perp}^\conjsym \, \mathbf{X}\|_F \qquad\text{\small (since $\|\mathbf{A}\|_F = \|\mathbf{A}^\conjsym\|_F$.)} \\
		& ~=~~ \sqrt{\frac{1}{n-q}} \sqrt{q}~ d_\text{\normalfont{\subspacetext}} (Y, X)  ~~=~~ \sqrt{\frac{q}{n-q}}  d_\text{\normalfont{\subspacetext}} (X, Y)
		\end{align*}
		
		\item The last property is obvious from the result of Lemma~\ref{lemma:d-subspace-equivalent-forms}.1.
	\end{enumerate}
\end{proof}
%


\subsection{Some Results Involving Set Distances} \label{sec:geomrtry}

\corrected{In this section we provide some geometry results that will be used in Section~\ref{sec:well-clusterted-spectrum} for computing the upper bounds on the perturbation of invariant subspaces in terms of the spectrum of the unperturbed matrix only.}
For the purpose of this paper and for simplicity, we consider only closed subsets of metric spaces in the following lemmas, although all these results can potentially be generalized for subsets that are open or/and closed in the metric space.

\begin{definition} \label{def:geomrtry-distances}
	Given closed subsets, $A,B$, of a metric space, $(\Psi,d)$, we define
	\begin{enumerate}
		\item \emph{Separation between the sets:} \[ \mathrm{sep}(A,B) ~=~ \min_{a\in A,\atop b\in B} d(a,b) \]
		\item \emph{Hausdorff distance between the sets:} \[ d_H(A,B) ~=~ \max \left( \max_{a\in A} \min_{b\in B} d(a,b) ~,~ \max_{b\in B} \min_{a\in A} d(a,b) \right) \]
		\item \emph{Diameter of a set:} \[ \mathrm{diam}(A) ~=~ \max_{a\in A,\atop a'\in A} d(a,a') \]
	\end{enumerate}
\end{definition}


\begin{lemma} \label{lemma:sep-diam-triangle-inequality}
If $(\Psi,d)$ is a metric space, then for any closed subsets, $P, Q, R \subseteq \Psi$,
\begin{equation}
\mathrm{sep}(P,Q) ~~\leq~~ \mathrm{sep}(P,R) ~+~ \mathrm{sep}(R,Q) ~+~ \mathrm{diam}(R)
\end{equation}
\end{lemma}

\begin{proof}
Let $(p^{*},r_1) \in \arg\!\min_{p\in P, \atop r\in R} d(p,r)$ (that is, $p^{*}\in P, ~r_1 \in R$ are a pair of points such that $d(p^{*},r_1) = \min_{p\in P, \atop r\in R} d(p,r) = \mathrm{sep}(P,R)$).
Likewise, let $(q^{*},r_2) \in \arg\!\min_{q\in Q, \atop r\in R} d(q,r)$ (that is, $d(q^{*}, r_2) = \mathrm{sep}(R,Q)$).
Then,
\begin{align}
\mathrm{sep}(P,Q) &~\leq~ d(p^{*},q^{*}) \quad\text{\small (since $\mathrm{sep}(P,Q) = \min_{p\in P,\atop q\in Q} d(p,q)$)} \nonumber \\
             &~\leq~ d(p^{*},r_1) ~+~ d(r_1,q^{*}) \quad\text{\small (triangle inequality.)} \nonumber \\
             &~=~ \mathrm{sep}(P,R) ~+~ d(r_1,q^{*}) \nonumber \\
             &~\leq~ \mathrm{sep}(P,R) ~+~ d(r_1,r_2) ~+~ d(q^{*}, r_2) \quad\text{\small (triangle inequality.)} \nonumber \\
             &~=~ \mathrm{sep}(P,R) ~+~ \mathrm{sep}(R,Q) ~+~ d(r_1,r_2) \nonumber \\
             &~\leq~ \mathrm{sep}(P,R) ~+~ \mathrm{sep}(R,Q) ~+~  \mathrm{diam}(R)
\end{align}
\end{proof}


\begin{lemma} \label{lemma:distance-lemma}
If $(\Psi,d)$ is a connected path metric space, then for any closed subsets, $P, Q, \widetilde{Q} \subseteq \Psi$,
\begin{equation}
\mathrm{sep}(P,Q)
~~\leq~~ 
\mathrm{sep}(P,\widetilde{Q})
~+~ 
d_H(\widetilde{Q}, Q)
\end{equation}
\end{lemma}

\begin{proof}~

\noindent \hspace{-0.5em}
\begin{tabular}{p{0.55\columnwidth}p{0.35\columnwidth}} 
\vspace{0.5em}
	Let $(p_0,q^{*}) \in \arg\!\min_{p\in P, \atop q\in Q} d(p,q)$ (that is, $p_0\in P, ~q^{*} \in Q$ are a pair of points such that $d(p_0,q^{*}) = \min_{p\in P, \atop q\in Q} d(p,q)$).\newline
	Likewise, let $(p_1,\widetilde{q}^{*}) \in \arg\!\min_{p\in P, \atop q'\in \widetilde{Q}} d(p,q')$.\newline
	Furthermore, let $\widetilde{q}^\dagger \in \arg\!\min_{q' \in \widetilde{Q}} d(q^{*}, q')$ and $q^\dagger \in \arg\!\min_{q \in Q} d(q, \widetilde{q}^{*})$. 
	&
\vspace{-2em}
	~~~~{\includegraphics[width=0.4\columnwidth, trim=0 0 0 0, clip=true]{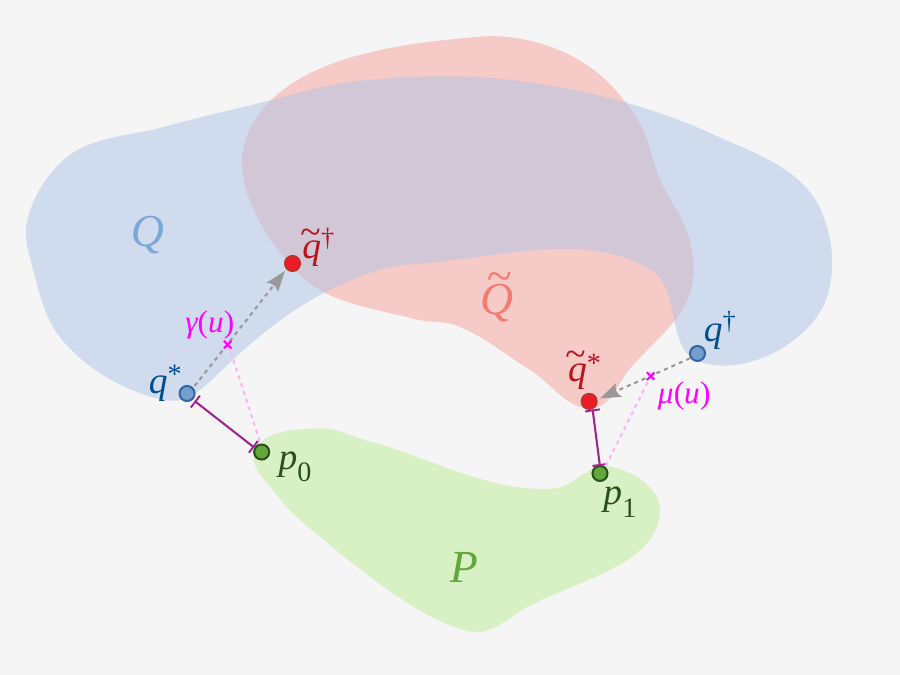}}
\end{tabular}


\vspace{-1em}
Consider a shortest path, $\gamma:[0,1]\rightarrow \Psi$, connecting $q^{*}$ and $\widetilde{q}^\dagger$, and parameterized by the normalized distance from $q^{*}$, so that $\gamma(0) = q^{*}$, $\gamma(1) = \widetilde{q}^\dagger$ and 
\begin{equation} \label{eq:d-qstar-gamma-lin}
d(q^{*},\gamma(u)) ~=~ u \,d(q^{*},\widetilde{q}^\dagger)
\end{equation}

Likewise, $\mu:[0,1]\rightarrow \Psi$ be the shortest path connecting $q^\dagger$ and $\widetilde{q}^{*}$,
nd parameterized by the normalized distance from $q^\dagger$,
so that $\mu(0) = q^\dagger$, $\mu(1) = \widetilde{q}^{*}$ and $d(q^\dagger,\mu(u)) = u \,d(q^\dagger,\widetilde{q}^{*})$.
Consequently, since $\mu(u)$ is a point on the shortest path connecting $q^\dagger$ and $\widetilde{q}^{*}$, we have 
\begin{equation} \label{eq:d-muu-tildeqstar-affine}
d(\mu(u),\widetilde{q}^{*}) ~=~ d(q^\dagger,\widetilde{q}^{*}) - d(q^\dagger,\mu(u)) ~=~ (1-u)\,d(q^\dagger,\widetilde{q}^{*})
\end{equation}

Define $f:[0,1]\rightarrow\mathbb{R}$ as $f(t) = d(p_0, \gamma(t))$, and $g:[0,1]\rightarrow\mathbb{R}$ as $g(t) = d(p_1, \mu(t))$.
It's easy to note that both $f$ and $g$ are continuous.

As a consequence, we have the following
\begin{align*}
f(0) = d(p_0, q^{*}) = \min_{p\in P, \atop q\in Q} d(p,q) ~\leq~ d(p_1, q^\dagger) = g(0) \\
g(1) = d(p_1, \widetilde{q}^{*}) = \min_{p\in P, \atop q'\in \widetilde{Q}} d(p,q') ~\leq~ d(p_0, \widetilde{q}^\dagger) = f(1)
\end{align*}
Thus, by intermediate value theorem, there exists a $u\in[0,1]$ such that $f(u)=g(u)$. That is,
\begin{equation} \label{eq:intermediate-value}
d(p_0, \gamma(u)) ~=~ d(p_1, \mu(u)), \qquad\text{for some $u\in [0,1]$.}
\end{equation}
Using this we have,
\begin{align*}
\min_{p\in P, \atop q\in Q} d(p,q) & ~=~ d(p_0,q^{*}) \\
& ~\leq~~ d(p_0,\gamma(u)) + d(q^{*},\gamma(u)) \qquad\text{\small (triangle inequality.)}\\
& ~=~~  d(p_1, \mu(u)) ~+~ d(q^{*},\gamma(u)) \qquad\text{\small (using \eqref{eq:intermediate-value}.)} \\
& ~\leq~~ d(p_1, \widetilde{q}^{*}) + d(\mu(u), \widetilde{q}^{*}) ~+~ d(q^{*},\gamma(u)) \qquad\text{\small (triangle inequality.)}\\
& ~=~ \min_{p\in P, \atop q'\in \widetilde{Q}} d(p,q') ~+~ d(\mu(u), \widetilde{q}^{*}) + d(q^{*},\gamma(u)) \\
& ~=~ \min_{p\in P, \atop q'\in \widetilde{Q}} d(p,q') ~+~ (1-u)\,d(q^\dagger,\widetilde{q}^{*}) + u \,d(q^{*},\widetilde{q}^\dagger) \qquad\text{\small (using \eqref{eq:d-qstar-gamma-lin} and \eqref{eq:d-muu-tildeqstar-affine}.)}\\
& ~\leq~ \min_{p\in P, \atop q'\in \widetilde{Q}} d(p,q') ~+~ \max \left( d(q^\dagger,\widetilde{q}^{*}) \,,\, d(q^{*},\widetilde{q}^\dagger) \right) \\ 
& ~=~~ \min_{p\in P, \atop q'\in \widetilde{Q}} d(p,q') ~+~ \max \left( \min_{q \in Q} d(q, \widetilde{q}^{*}) \,,\, \min_{q' \in \widetilde{Q}} d(q^{*}, q') \right) \quad\text{\small (definitions of $q^\dagger$ and $\widetilde{q}^\dagger$.)}\\
%
& ~\leq~~ \min_{p\in P, \atop q'\in \widetilde{Q}} d(p,q') ~+~ \max \left( \max_{q' \in \widetilde{Q}} \, \min_{q\in Q} d(q,q') ~,~ \max_{q\in Q} \, \min_{q' \in \widetilde{Q}} d(q,q') \right) \\
& ~=~ \mathrm{sep}(P,\widetilde{Q}) ~+~ d_H(\widetilde{Q}, Q)
\end{align*}

\end{proof}


\begin{lemma} \label{lemma:sep-preserving-purturbation}
	Suppose $P, Q, \widetilde{R}$ are closed subsets of a metric space, $(\Psi,d)$, such that

\noindent	
\begin{tabular}{p{0.55\columnwidth}p{0.35\columnwidth}}
		\vspace{-2em}
    \begin{equation} \displaystyle 
		\max_{r'\in \widetilde{R}} \min_{s\in P\cup Q} d(s,r') \,+\, d_H (P\cup Q, \widetilde{R})
	    ~<~ 
	    \mathrm{sep} (P, Q) 
	    \label{eq:lemma-sep-preserving-purturbation}
	\end{equation}
    
    \noindent Define, $\widetilde{P}, \widetilde{Q} \subseteq \widetilde{R}$, such that
    \begin{eqnarray} 
    \widetilde{P}  = \{r' \in \widetilde{R} \,|\, \min_{s\in P\cup Q} d(s,r') = \min_{p\in P} d(p,r') \},
    ~~~~\text{and,} \nonumber \\
    \widetilde{Q}  = \{r' \in \widetilde{R} \,|\, \min_{s\in P\cup Q} d(s,r') = \min_{q\in Q} d(q,r') \} \quad
    \label{eq:tilde-P-Q}
    \end{eqnarray}
    &
    \vspace{1em}
    {\includegraphics[width=0.4\columnwidth, trim=0 0 0 0, clip=true]{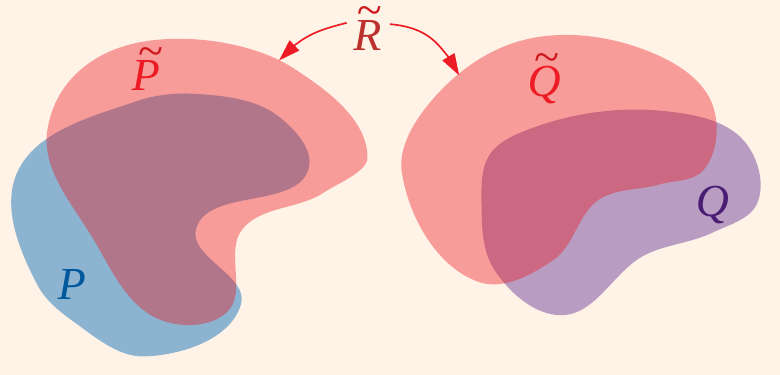}}
\end{tabular}

\vspace{-1em}
\noindent Then
    \begin{enumerate}
    	\item $\{\widetilde{P}, \widetilde{Q}\}$ constitutes a partition of $\widetilde{R}$,
    	
    	\item $\arg\!\min_{s \in P \cup Q} d(s,p') \subseteq P, ~\forall p' \in \widetilde{P}$,
    	 ~~and,~~ $\arg\!\min_{s \in P \cup Q} d(s,q') \subseteq Q, ~\forall q' \in \widetilde{Q}$. \newline 
    	 {\small (consequently, $\displaystyle \min_{s \in P \cup Q} d(s,p') = \min_{s \in P} d(s,p'), ~\forall p' \in \widetilde{P}$,
    	 ~and,~ $\displaystyle \min_{s \in P \cup Q} d(s,q') = \min_{s \in Q} d(s,q'), ~\forall q' \in \widetilde{Q}$.)}
    	
    	\item $\arg\!\min_{r' \in \widetilde{R}} d(p,r') \subseteq \widetilde{P}, ~\forall p \in P$,
    	~~and,~~ $\arg\!\min_{r' \in \widetilde{R}} d(q,r') \subseteq \widetilde{Q}, ~\forall q \in Q$. \newline 
    	{\small (consequently, $\displaystyle \min_{r' \in \widetilde{R}} d(p,r') = \min_{r' \in \widetilde{P}} d(p,r'), ~\forall p \in P$,
    	~and,~ $\displaystyle \min_{r' \in \widetilde{R}} d(q,r') = \min_{r' \in \widetilde{Q}} d(q,r'), ~\forall q \in Q$.)}
    	
    	\item $d_H(P,\widetilde{P}) \leq d_H(P\cup Q, \widetilde{R})$,~~ $d_H(Q,\widetilde{Q}) \leq d_H(P\cup Q,\widetilde{R})$, \newline \phantom{xxxx} and,~ $\max\left( d_H(P,\widetilde{P}) , d_H(Q,\widetilde{Q}) \right) = d_H(P\cup Q, \widetilde{R})$.
    	
    	\item If $(\Psi,d)$ is a connected path metric space, then ~$ \displaystyle
    	\mathrm{sep}(\widetilde{P}, \widetilde{Q}) ~\geq~  \mathrm{sep} (P, Q) - 2 
		~d_H(P\cup Q, \widetilde{R})
    	$
    \end{enumerate}
If the above holds, 
we say 
``$\widetilde{R}$ is a separation-preserving perturbation of $P$ and $Q$'', and call $\{\widetilde{P}, \widetilde{Q}\}$ to be the ``separation-preserving partition of $\widetilde{R}$''.
\end{lemma}

\begin{proof}~

\begin{enumerate}

\item 
We first prove that $\{\widetilde{P}, \widetilde{Q}\}$ constitutes of a partition of $\widetilde{R}$.
\begin{quote}
\emph{Proof for $\widetilde{P} \cup \widetilde{Q} = \widetilde{R}$:}
For a fixed $r'\in \widetilde{R}$, an element of $\,\arg\!\!\min_{s\in P \cup Q} d(s,r')$ is either in $P$ or in $Q$.
In the former case the point $r'$ will belong to $\widetilde{P}$, while in the later case it will belong to $\widetilde{Q}$ (with the possibility that it belongs to both) due to the definition \eqref{eq:tilde-P-Q}.
Thus there does not exist a point $r'\in \widetilde{R}$ that does not belong to either $\widetilde{P}$ or $ \widetilde{Q}$.

\emph{Proof for $\widetilde{P} \cap \widetilde{Q} = \emptyset$:}
We prove this by contradiction. If possible, let $\rho' \in \widetilde{P} \cap \widetilde{Q}$.
Since $\rho' \in \widetilde{P}$, due to definition \eqref{eq:tilde-P-Q}, there exists a $p_1 \in P$ such that $\min_{s\in P\cup Q} d(s,\rho') = d(p_1,\rho')$.
Likewise, there exists a $q_1 \in Q$ such that $\min_{s\in P\cup Q} d(s,\rho') = d(q_1,\rho')$.
Thus,
\begin{align*}
& \qquad \begin{array}{rcl}
	\displaystyle 2 \min_{s\in P\cup Q} d(s,\rho') & = &	d(p_1,\rho') + d(q_1,\rho') \\
	& \geq & d(p_1,q_1) \qquad\text{\small (tringle inequality.)} \\
	& \geq & \displaystyle \min_{p\in P,\atop q\in Q} d(p,q) \qquad\text{\small (since $p_1\in P, q_1\in Q$.)}
\end{array} \\
& \begin{array}{rl}
\Rightarrow & \displaystyle 2\, \max_{r'\in \widetilde{R}} \min_{s\in P\cup Q} d(s,r')  ~\geq~ \min_{p\in P,\atop q\in Q} d(p,q) \\
\Rightarrow & \displaystyle \max_{r'\in \widetilde{R}} \min_{s\in P\cup Q} d(s,r') \,+\, d_H (P\cup Q, \widetilde{R}) 
				~\geq~ \mathrm{sep}(P,Q)
\end{array}
\end{align*}
This contradicts the assumption \eqref{eq:lemma-sep-preserving-purturbation} of the Lemma. Hence there cannot exist a $\rho' \in \widetilde{P} \cap \widetilde{Q}$. Thus $\widetilde{P} \cap \widetilde{Q} = \emptyset$.
\end{quote}

\item
We next prove $\arg\!\min_{s \in P \cup Q} d(s,p') \subseteq P, ~\forall p' \in \widetilde{P}$.
We do this by contradiction.

If possible, suppose there exists a $p' \in \widetilde{P}$ such that ~$\arg\!\min_{s \in P \cup Q} d(s,p') \not\subseteq P$. Then there exists a $q\in Q$ such that $\min_{s \in P \cup Q} d(s,p') = d(q,p')$.
But $d(q,p') \geq \min_{s\in Q} d(s,p') \geq \min_{s \in P \cup Q} d(s,p')$. This implies $\min_{s \in P \cup Q} d(s,p') = \min_{s\in Q} d(s,p')$.
Due to definition of $\widetilde{Q}$ in \eqref{eq:tilde-P-Q} this implies $p'\in \widetilde{Q}$. However, we have already shown that $\widetilde{P} \cap \widetilde{Q} = \emptyset$. This leads to a contradiction. Thus $\arg\!\min_{s \in P \cup Q} d(s,p') \subseteq P, ~\forall p' \in \widetilde{P}$.

Likewise we can prove $\arg\!\min_{s \in P \cup Q} d(s,q') \subseteq Q, ~\forall q' \in \widetilde{Q}$.

\item
We next prove $\arg\!\min_{r' \in \widetilde{R}} d(p,r') \subseteq \widetilde{P}, ~\forall p \in P$.
We do this by contradiction.

If possible, suppose there exists a $p_3\in P$ such that $\arg\!\min_{r' \in \widetilde{R}} d(p_3,r') \not\subseteq \widetilde{P}$.
Then there exists a $\rho' \in \widetilde{Q}$ such that $\min_{r' \in \widetilde{R}} d(p_3,r') = d(p_3,\rho')$.

Again, due to the definition of $\widetilde{Q}$ in \eqref{eq:tilde-P-Q}, for any $\rho' \in \widetilde{Q}$ there exists a $q_3\in Q$ such that ~$d(q_3,\rho') = \min_{s\in P\cup Q} d(s,\rho')$.

Thus,
\begin{align*}
& \qquad \begin{array}{rcl}
\displaystyle \min_{r' \in \widetilde{R}} d(p_3,r') + \min_{s\in P\cup Q} d(s,\rho') & = & d(p_3,\rho') + d(q_3,\rho') \\
& \geq & d(p_3,q_3) \qquad\text{\small (tringle inequality.)} \\
& \geq & \displaystyle \min_{p\in P,\atop q\in Q} d(p,q) \qquad\text{\small (since $p_3\in P, q_3\in Q$.)}
\end{array}\\
& \begin{array}{rl}
\Rightarrow & \displaystyle \max_{s\in P\cup Q} \min_{r' \in \widetilde{R}} d(s,r') + \max_{r' \in \widetilde{R}} \min_{s\in P\cup Q} d(s,r') ~\geq~ \min_{p\in P,\atop q\in Q} d(p,q) \\
\Rightarrow &  \displaystyle d_H (P\cup Q, \widetilde{R}) \,+\, \max_{r'\in \widetilde{R}} \min_{s\in P\cup Q} d(s,r') 
			~\geq~ \mathrm{sep}(P,Q)
\end{array}
\end{align*}
This contradicts the assumption \eqref{eq:lemma-sep-preserving-purturbation} of the Lemma. Hence there cannot exist a $p_3\in P$ such that $\arg\!\min_{r' \in \widetilde{R}} d(p_3,r') \not\subseteq \widetilde{P}$. Thus $\arg\!\min_{r' \in \widetilde{R}} d(p,r') \subseteq \widetilde{P}, ~\forall p \in P$.

Likewise we can prove $\arg\!\min_{r' \in \widetilde{R}} d(q,r') \subseteq \widetilde{Q}, ~\forall q \in Q$.

\item
Since $\arg\!\min_{s \in P \cup Q} d(s,p') \subseteq P, ~\forall p' \in \widetilde{P}$, we have $\min_{s \in P \cup Q} d(s,p') = \min_{p \in P} d(p,p'), ~\forall p' \in \widetilde{P}$. Thus,
$\max_{p'\in \widetilde{P}} \min_{p \in P} d(p,p') = \max_{p'\in \widetilde{P}} \min_{s \in P \cup Q} d(s,p') $.

Likewise, since $\arg\!\min_{r' \in \widetilde{R}} d(p,r') \subseteq \widetilde{P}, ~\forall p \in P$, we have \newline $\max_{p \in P} \min_{p'\in \widetilde{P}} d(p,p') = \max_{p \in P} \min_{r'\in \widetilde{R}} d(p,r')$.

Thus,
\begin{align}
d_H(P,\widetilde{P}) & ~=~ \max \left(\, \max_{p \in P} \min_{p'\in \widetilde{P}} d(p,p') ~,~ \max_{p'\in \widetilde{P}} \min_{p \in P} d(p,p') \,\right) \nonumber \\
& ~=~  \max \left(\, \max_{p \in P} \min_{r'\in \widetilde{R}} d(p,r') ~,~ \max_{p'\in \widetilde{P}} \min_{s \in P \cup Q} d(s,p') \,\right) \label{eq:dH-P-Ptilde} \\
& ~\leq~  \max \left(\, \max_{s \in P \cup Q} \min_{r'\in \widetilde{R}} d(s,r') ~,~ \max_{r'\in \widetilde{R}} \min_{s \in P \cup Q} d(s,r') \,\right) \quad\text{\small (since $P \subseteq P\cup Q, ~\widetilde{P} \subseteq \widetilde{R}$.)} \nonumber \\
& ~=~ d_H(P\cup Q,\widetilde{R}) \nonumber
\end{align}

Similarly we can show,
\begin{align}
d_H(Q,\widetilde{Q}) 
& ~=~  \max \left(\, \max_{q \in Q} \min_{r'\in \widetilde{R}} d(q,r') ~,~ \max_{q'\in \widetilde{Q}} \min_{s \in P \cup Q} d(s,q') \,\right) \label{eq:dH-Q-Qtilde} \\
& ~\leq~ d_H(P\cup Q,\widetilde{R}) \nonumber
\end{align}
Again, from \eqref{eq:dH-P-Ptilde} and \eqref{eq:dH-Q-Qtilde},
\begin{align*}
\max \left( d_H(P,\widetilde{P}) , d_H(Q,\widetilde{Q})  \right)
& ~=~ \max\left( \max_{p \in P} \min_{r'\in \widetilde{R}} d(p,r') ~,~ \max_{q \in Q} \min_{r'\in \widetilde{R}} d(q,r') ~, \right. \\
& \qquad\quad\qquad \left. \max_{p'\in \widetilde{P}} \min_{s \in P \cup Q} d(s,p')  ~,~ \max_{q'\in \widetilde{Q}} \min_{s \in P \cup Q} d(s,q')  \right) \\
& ~=~ \max\left( \max_{p \in P\cup Q} \min_{r'\in \widetilde{R}} d(p,r') ~,~ \max_{p'\in \widetilde{P}\cup \widetilde{Q}} \min_{s \in P \cup Q} d(s,p') \right) \\
& ~=~ d_H(P\cup Q,\widetilde{R}) \qquad\text{\small (since $\widetilde{P}\cup \widetilde{Q} = \widetilde{R}$)}
\end{align*}

\item 
\begin{align*}
\mathrm{sep}(\widetilde{P}, \widetilde{Q}) & \geq~  \mathrm{sep} (\widetilde{P}, Q) - 
d_H(Q, \widetilde{Q})
\quad\text{\small (using Lemma~\ref{lemma:distance-lemma}.)} \\ 
& \geq~  \mathrm{sep} (P, Q) 
- d_H(P, \widetilde{P})
- d_H(Q, \widetilde{Q})
 \quad\text{\small (using Lemma~\ref{lemma:distance-lemma}.)} \\
& \geq~  \mathrm{sep} (P, Q) - 2 
~d_H(P\cup Q, \widetilde{R}) \\
& \qquad\qquad\text{\small (since $d_H(P,\widetilde{P}) \leq d_H(P\cup Q, \widetilde{R})$ ~and~ $d_H(Q,\widetilde{Q}) \leq d_H(P\cup Q,\widetilde{R})$.)} 
\end{align*}

\end{enumerate}

\end{proof}


\begin{cor} \label{corr:sep-preserving-purturbation}
	If $P, Q, \widetilde{R}$ are closed subsets of a metric space, $(\Psi,d)$, such that
	$d_H (P\cup Q, \widetilde{R})
	~<~ \frac{1}{2} \, 
	\mathrm{sep} (P, Q)$, then
	$\widetilde{R}$ is a separation-preserving perturbation of $P$ and $Q$.
	
	As a consequence, the separation-preserving partition, $\{\widetilde{P}, \widetilde{Q}\}$, of $\widetilde{R}$ as defined in \eqref{eq:tilde-P-Q}
	satisfies properties `1' to `4' in Lemma~\ref{lemma:sep-preserving-purturbation}, as well as property `5' (if $(\Psi,d)$ is a connected path metric space) with an additional inequality: \[ \mathrm{sep}(\widetilde{P}, \widetilde{Q})  ~\geq~  \mathrm{sep} (P, Q) - 2 ~d_H(P\cup Q, \widetilde{R}) ~>~ 0\]
\end{cor}

\begin{proof}
	The result follows directly from Lemma~\ref{lemma:sep-preserving-purturbation} by observing that
	\[ \max_{r'\in \widetilde{R}} \min_{s\in P\cup Q} d(s,r') \,+\, d_H (P\cup Q, \widetilde{R})
	~~\leq~~ 2 ~d_H (P\cup Q, \widetilde{R})
	~~<~~ 
	\mathrm{sep} (P, Q) \]
\end{proof}

\section{\corrected{Results on Perturbation Upper Bounds}}

\corrected{Throughout this section we use the notations and conventions described in Definition~\ref{def:notations}.}

\subsection{Elementary Results on Spectrum Perturbation}

In this section we provide some elementary results relating the norm of the matrix perturbation and perturbation of the eigenvalues and eigenvectors.

\begin{lemma} \label{lemma:D-matrix-equality}
	Define $D \in \mathbb{C}^{n\times n}$ such that $D_{jj'} = (\widetilde{\lambda}_{j'} - \lambda_j) ~\mathbf{u}_j^\conjsym \widetilde{\mathbf{u}}_{j'} $.
	Then
	\begin{equation}
	D ~=~ U^\conjsym (\widetilde{M} - M) \widetilde{U}
	\end{equation}
	Equivalently,
	\begin{equation} \label{eq:MD-element-wise-equality}
	(\widetilde{\lambda}_{j'} - \lambda_j) ~\mathbf{u}_j^\conjsym \widetilde{\mathbf{u}}_{j'} ~=~ \mathbf{u}_j^\conjsym (\widetilde{M} - M) \widetilde{\mathbf{u}}_{j'}, \qquad \forall j,j' \in N
	\end{equation}
	The later relation in fact holds even when $\widetilde{M}$ is not normal but $\widetilde{\mathbf{u}}_{j'}$ is simply a right eigenvector of $\widetilde{M}$ with corresponding eigenvalue $\widetilde{\lambda}_{j'}$.
\end{lemma}

\begin{proof}
	First we note that since $M$ is normal with $\mathbf{u}_j$ a right eigenvector and corresponding eigenvalue $\lambda_j$, $\mathbf{u}_j^\conjsym$ is a left eigenvector of $M$ with the same eigenvalue. Thus,
	\begin{align*}
	\mathbf{u}_j^\conjsym (\widetilde{M} - M) \widetilde{\mathbf{u}}_{j'} & =~ \mathbf{u}_j^\conjsym ~\widetilde{M}  \widetilde{\mathbf{u}}_{j'} - \mathbf{u}_j^\conjsym M~ \widetilde{\mathbf{u}}_{j'} 
	~=~ \mathbf{u}_j^\conjsym \widetilde{\lambda}_{j'}  \widetilde{\mathbf{u}}_{j'} - \lambda_j \mathbf{u}_j^\conjsym \widetilde{\mathbf{u}}_{j'} 
	~=~ (\widetilde{\lambda}_{j'} - \lambda_j) ~\mathbf{u}_j^\conjsym \widetilde{\mathbf{u}}_{j'}
	\end{align*}
	This proves \eqref{eq:MD-element-wise-equality}.

%
	We note that if both $M$ and $\widetilde{M}$ are normal, the L.H.S. of \eqref{eq:MD-element-wise-equality} is the $(j,j')$-th element of $U^\conjsym (\widetilde{M} - M) \widetilde{U}$ and the R.H.S. is $D_{jj'}$. 
\end{proof}

\begin{cor} \label{cor:2-norm-sum-equality}
	\begin{eqnarray}
	\left\| \widetilde{M} - M \right\|^2_2 ~~\geq~~ \left\| (\widetilde{M} - M) \widetilde{\mathbf{u}}_{j'} \right\|^2_2 & = & \sum_{j=1}^n \left|\widetilde{\lambda}_{j'} - \lambda_j\right|^2 ~\left| \mathbf{u}_j^\conjsym \widetilde{\mathbf{u}}_{j'}\right|^2, \quad\forall j' \in N \nonumber \\
	\left\| \widetilde{M} - M \right\|^2_2 ~~\geq~~  \left\| (\widetilde{M} - M) \mathbf{u}_{j} \right\|^2_2 & = & \sum_{j'=1}^n \left|\widetilde{\lambda}_{j'} - \lambda_j\right|^2 ~\left| \mathbf{u}_j^\conjsym \widetilde{\mathbf{u}}_{j'}\right|^2, \quad\forall j \in N
	\end{eqnarray}
	The first relation holds even when $\widetilde{M}$ is not normal, while the second relation holds even when $M$ is not normal.
\end{cor}

\begin{proof}
	The inequalities follows from the definition of induced $2$-norm for matrices.
	
	\vspace{0.5em} \noindent
	When $M$ is normal, $\{\mathbf{u}_{j}\}_{j\in N}$ forms an orthonormal basis in $\mathbb{C}^n$.
	Noting that \eqref{eq:MD-element-wise-equality} is a scalar equation, multiplying on both sides with $\mathbf{u}_{j}$ and summing over $j$, we get
{\small\begin{equation*}
		\sum_{j=1}^n \left( (\widetilde{\lambda}_{j'} - \lambda_j) ~\mathbf{u}_j^\conjsym \widetilde{\mathbf{u}}_{j'} \right) \mathbf{u}_{j} 
		~=~ \sum_{j=1}^n \mathbf{u}_{j}  \left( \mathbf{u}_j^\conjsym (\widetilde{M} - M) \widetilde{\mathbf{u}}_{j'}  \right) 
		~=~ \left( \sum_{j=1}^n  \mathbf{u}_{j}  \mathbf{u}_j^\conjsym \right) (\widetilde{M} - M) \widetilde{\mathbf{u}}_{j'} 
		~=~ I ~(\widetilde{M} - M) \widetilde{\mathbf{u}}_{j'}
\end{equation*}}
	 Taking the $2$-norm on both sides of the above gives the first equality.
	
	Switching the roles of tilde and non-tilde terms in Lemma~\ref{lemma:D-matrix-equality} and the above gives the second relation.
\end{proof}

\begin{cor}~ \label{cor:M-diff-lambda-diff}
	\begin{enumerate}
		\item $\left\| \widetilde{M} - M \right\|_2 \geq \left\| (\widetilde{M} - M) \widetilde{\mathbf{u}}_{j'} \right\|_2 \geq \min_{j\in N} \left|\widetilde{\lambda}_{j'} - \lambda_j\right|, ~~\forall j'\in N$, 
		~~and~~ \newline
		\phantom{.} \qquad\qquad 
		 $\left\| \widetilde{M} - M \right\|_2 \geq \left\| (\widetilde{M} - M) {\mathbf{u}}_{j} \right\|_2 \geq \min_{j'\in N} \left|\widetilde{\lambda}_{j'} - \lambda_j\right|, ~~\forall j\in N$.
		 \newline 
		The first relation holds even when $\widetilde{M}$ is not normal, while the second relation holds even when $M$ is not normal.
		
		\item \noindent The following results are a consequence of the Bauer-Fike Theorem for normal matrices~\cite{BauF60b}:
		\begin{align*}
		\left\| \widetilde{M} - M \right\|_2 
		&~\geq~~ \max_{j\in N} \left\| (\widetilde{M} - M) \widetilde{\mathbf{u}}_{j'} \right\|_2
		~~\geq~~ \max_{j'\in N} \,\min_{j\in N} \,\left|\widetilde{\lambda}_{j'} - \lambda_j\right| \\
		\left\| \widetilde{M} - M \right\|_2 
		&~\geq~~ \max_{j\in N} \left\| (\widetilde{M} - M) {\mathbf{u}}_{j} \right\|_2
		~~\geq~~ \max_{j\in N} \,\min_{j'\in N} \,\left|\widetilde{\lambda}_{j'} - \lambda_j\right|
		\end{align*}
		Once again, the first relation holds even when $\widetilde{M}$ is not normal, while the second relation holds even when $M$ is not normal.
	\end{enumerate}
\end{cor}

\begin{proof}
	From the result of Corollary~\ref{cor:2-norm-sum-equality}, 
	when $M$ is normal (and $\widetilde{M}$ is not necessarily normal), for all $j' \in N$,
	\begin{align*}
	\left\| \widetilde{M} - M \right\|^2_2 & \geq~~ \left\| (\widetilde{M} - M) \widetilde{\mathbf{u}}_{j'} \right\|^2_2 \\
	& =~ \sum_{j=1}^n \left|\widetilde{\lambda}_{j'} - \lambda_j\right|^2 ~\left| \mathbf{u}_j^\conjsym \widetilde{\mathbf{u}}_{j'}\right|^2 \\
	& \geq~ \min_{j\in N} \left|\widetilde{\lambda}_{j'} - \lambda_j\right|^2 \sum_{j=1}^n \left| \mathbf{u}_j^\conjsym  \widetilde{\mathbf{u}}_{j'}\right|^2 \\
	&=~ \min_{j\in N} \left|\widetilde{\lambda}_{j'} - \lambda_j\right|^2 ~\|\widetilde{\mathbf{u}}_{j'}\|^2  ~~\text{\small (since $\{\mathbf{u}_j\}_{j\in N}$ forms an orthonormal basis.)} \\
	&=~ \min_{j\in N} \left|\widetilde{\lambda}_{j'} - \lambda_j\right|^2
	\end{align*}
%
%
	Since this is true for any $j'\in N$, it follows that 
	$\left\| \widetilde{M} - M \right\|_2 \geq \max_{j'\in N} \,\min_{j\in N} \,\left|\widetilde{\lambda}_{j'} - \lambda_j\right|$.
	
	A similar set of results can be derived with the tilde and non-tilde terms exchanged. 
\end{proof}





\subsection{Distance Between Invariant Subspaces of Normal Matrices with Partitioned Spectra} \label{sec:partitioned-spectrum}


Suppose $J, \widetilde{J} \subseteq N$ such that $|J| = |\widetilde{J}| = q$.
%
We are interested in understanding how much the invariant space $\spn(\mathbf{u}_J)$ of $M$ differs from the invariant space $\spn(\widetilde{\mathbf{u}}_{\widetilde{J}})$ of $\widetilde{M}$.
\corrected{The results in this section are variations and modest improvements on the Davis-Kahan $\sin\Theta$ Theorem~\cite{10.1137/0707001} (see Section VIII.3 of~\cite{bhatia1996matrix}, for example).}
In Proposition~\ref{prop:full-main} and the two corollaries that follow, we present results of the form
\[
d_\mathrm{\subspacetext} \left( \spn(\mathbf{u}_J) , \spn(\widetilde{\mathbf{u}}_{\widetilde{J}}) \right) ~~\leq~~
\mathscr{F}(\widetilde{M} \!-\! M,\, \mathbf{u}_N,\, \lambda_N,\, \widetilde{\lambda}_N;\, J,\, \widetilde{J})
\]
where $\mathscr{F}$ is a function specific to the exact statement of the proposition or corollary.


%

For a given invariant subspace $\spn(\mathbf{u}_J)$ of $M$, we can consider all the possible $q$-dimensional invariant subspaces of $\widetilde{J}$ and choose the one that is closest to $\spn(\mathbf{u}_J)$ as its perturbation. \corrected{As a consequence,
for any of these results we can write}
\[
\min_{\widetilde{J}\in S_{q,n}} d_\mathrm{\subspacetext} \left( \spn(\mathbf{u}_J) , \spn(\widetilde{\mathbf{u}}_{\widetilde{J}}) \right) ~~\leq~~
\min_{\widetilde{J}\in S_{q,n}} \mathscr{F}(\widetilde{M} \!-\! M,\, \mathbf{u}_N,\, \lambda_N,\, \widetilde{\lambda}_N;\, J,\, \widetilde{J})
\]
where $S_{q,n}$ is the set of all $q$-element subsets of $N = \{1,2,\cdots,n\}$.
This gives a combinatorial means of finding the $q$-dimensional invariant subspace of $\widetilde{M}$ that is closest to $\spn(\mathbf{u}_J)$.




\begin{definition}
For $a,b,c \in \mathbb{R}$ with $a \leq \min(b,c)$ we define~
$	[a,\min(b,c_{-})] = \left\{ 
					\begin{array}{ll}
						{[a,b]} & \text{if $c>b$} \\
						{[a,c)} & \text{if $c\leq b$}
					\end{array} 
					\right.  $.
\end{definition}

\begin{prop} \label{prop:full-main}
	For any $J, \widetilde{J} \subseteq N$ with $|J| = |\widetilde{J}|=q$, 
\begin{align} \label{eq:prop-full-main}
	d_\mathrm{\subspacetext} \left( \spn(\mathbf{u}_J) , \spn(\widetilde{\mathbf{u}}_{\widetilde{J}}) \right) ~~\leq~~ \sqrt{ \frac{1}{q}\, \mathlarger{\mathlarger{\sum}}_{j\in J} ~\frac{ 
		\left\| (\widetilde{M} - M) \mathbf{u}_{{j}} \right\|^2_2 ~-~ \kappa_j \,\min_{j'\in \widetilde{J}}\left|\widetilde{\lambda}_{{j'}} - \lambda_j\right|^2 
	}{ 
		\min_{j'\in {\widetilde{J}}^c} \left|\widetilde{\lambda}_{{j'}} - \lambda_j\right|^2 - \kappa_j \,\min_{j'\in \widetilde{J}}\left|\widetilde{\lambda}_{{j'}} - \lambda_j\right|^2 
	} }
\end{align}
for any $\kappa_j \in \Big[ 0,\min\left(1, \left( {\small \frac{ \min_{j'\in {\widetilde{J}}^c} \left|\widetilde{\lambda}_{{j'}} - \lambda_j\right|^2}{ \min_{j'\in \widetilde{J}}\left|\widetilde{\lambda}_{{j'}} - \lambda_j\right|^2}} \right)_{\!\!-} \, \right) \Big], ~j\in J$.

The tightest bound in \eqref{eq:prop-full-main} is obtained by choosing
\begin{equation}
\kappa_j = \left\{ \begin{array}{ll}
0, & \text{if ~$\| (\widetilde{M} - M) \mathbf{u}_{{j}} \|_2 \geq \min_{j'\in {\widetilde{J}}^c} |\widetilde{\lambda}_{{j'}} - \lambda_j|$} \\
1, & \text{if ~$\| (\widetilde{M} - M) \mathbf{u}_{{j}} \|_2 < \min_{j'\in {\widetilde{J}}^c} |\widetilde{\lambda}_{{j'}} - \lambda_j|$}
\end{array}\right.
\end{equation}
\end{prop}

\begin{proof}
From Corollary~\ref{cor:2-norm-sum-equality}, for all ${j} \in N$
{\small \begin{align}
		&		\left\| (\widetilde{M} - M) \mathbf{u}_{{j}} \right\|^2_2 
		~~=~~ \sum_{j'\in \widetilde{J}} \left|\widetilde{\lambda}_{{j'}} - \lambda_j\right|^2 ~\left| \mathbf{u}_j^\conjsym \widetilde{\mathbf{u}}_{{j'}}\right|^2 ~+~ \sum_{j'\in {\widetilde{J}}^c} \left|\widetilde{\lambda}_{{j'}} - \lambda_j\right|^2 ~\left| \mathbf{u}_j^\conjsym \widetilde{\mathbf{u}}_{{j'}}\right|^2 \nonumber \\
		& \qquad\qquad\qquad\qquad \geq~~ \min_{j'\in \widetilde{J}}\left|\widetilde{\lambda}_{{j'}} - \lambda_j\right|^2 ~\sum_{j'\in \widetilde{J}} \left| \mathbf{u}_j^\conjsym \widetilde{\mathbf{u}}_{{j'}}\right|^2 ~+~ \min_{j'\in {\widetilde{J}}^c} \left|\widetilde{\lambda}_{{j'}} - \lambda_j\right|^2 ~\sum_{j'\in {\widetilde{J}}^c} ~\left| \mathbf{u}_j^\conjsym \widetilde{\mathbf{u}}_{{j'}}\right|^2 \label{eq:first-ineq} \\
		& \qquad\qquad\qquad\qquad =~~ \min_{j'\in \widetilde{J}}\left|\widetilde{\lambda}_{{j'}} - \lambda_j\right|^2 ~\left(1 - \sum_{j'\in {\widetilde{J}}^c} \left| \mathbf{u}_j^\conjsym \widetilde{\mathbf{u}}_{{j'}}\right|^2 \right) ~+~ \min_{j'\in {\widetilde{J}}^c} \left|\widetilde{\lambda}_{{j'}} - \lambda_j\right|^2 ~\sum_{j'\in {\widetilde{J}}^c} ~\left| \mathbf{u}_j^\conjsym \widetilde{\mathbf{u}}_{{j'}}\right|^2 \nonumber \\
		& \qquad\qquad\qquad\qquad \geq~~ \kappa_j \,\min_{j'\in \widetilde{J}}\left|\widetilde{\lambda}_{{j'}} - \lambda_j\right|^2 ~\left(1 - \sum_{j'\in {\widetilde{J}}^c} \left| \mathbf{u}_j^\conjsym \widetilde{\mathbf{u}}_{{j'}}\right|^2 \right) ~+~ \min_{j'\in {\widetilde{J}}^c} \left|\widetilde{\lambda}_{{j'}} - \lambda_j\right|^2 ~\sum_{j'\in {\widetilde{J}}^c} ~\left| \mathbf{u}_j^\conjsym \widetilde{\mathbf{u}}_{{j'}}\right|^2 \nonumber \\
		& \qquad\qquad\qquad\qquad\qquad\qquad\qquad\qquad\qquad\qquad\qquad\qquad\qquad\qquad\qquad \text{for any $\kappa_j \in [0,1]$.} \nonumber \\
		%
		\Rightarrow~~ &
		\left( \min_{j'\in {\widetilde{J}}^c} \left|\widetilde{\lambda}_{{j'}} - \lambda_j\right|^2 - \kappa_j \,\min_{j'\in \widetilde{J}}\left|\widetilde{\lambda}_{{j'}} - \lambda_j\right|^2  \right) \sum_{j'\in {\widetilde{J}}^c} ~\left| \mathbf{u}_j^\conjsym \widetilde{\mathbf{u}}_{{j'}}\right|^2
		~~\leq~~ 
		\left\| (\widetilde{M} - M) \mathbf{u}_{{j}} \right\|^2_2 ~-~ \kappa_j \,\min_{j'\in \widetilde{J}}\left|\widetilde{\lambda}_{{j'}} - \lambda_j\right|^2 \label{eq:prop-full-main-before-division} \\
		\Rightarrow~~ & \sum_{j'\in {\widetilde{J}}^c} ~\left| \mathbf{u}_j^\conjsym \widetilde{\mathbf{u}}_{{j'}}\right|^2
		~~\leq~~ \frac{ 
			\left\| (\widetilde{M} - M) \mathbf{u}_{{j}} \right\|^2_2 ~-~ \kappa_j \,\min_{j'\in \widetilde{J}}\left|\widetilde{\lambda}_{{j'}} - \lambda_j\right|^2 
		}{ 
			\min_{j'\in {\widetilde{J}}^c} \left|\widetilde{\lambda}_{{j'}} - \lambda_j\right|^2 - \kappa_j \,\min_{j'\in \widetilde{J}}\left|\widetilde{\lambda}_{{j'}} - \lambda_j\right|^2 
		} \nonumber \\
		& \qquad\qquad\qquad\qquad\qquad\qquad\qquad\qquad \text{for any $\kappa_j \in \Big[ 0,\min\left(1, \left( {\small \frac{ \min_{j'\in {\widetilde{J}}^c} \left|\widetilde{\lambda}_{{j'}} - \lambda_j\right|^2}{ \min_{j'\in \widetilde{J}}\left|\widetilde{\lambda}_{{j'}} - \lambda_j\right|^2}} \right)_{\!\!-} \, \right) \Big]$.} 
		\label{eq:prop-main-for-a-j}
\end{align}}
In the last step, we ensured that $\min_{j'\in {\widetilde{J}}^c} \left|\widetilde{\lambda}_{{j'}} - \lambda_j\right|^2 - \kappa_j \,\min_{j'\in \widetilde{J}}\left|\widetilde{\lambda}_{{j'}} - \lambda_j\right|^2 $ is positive by restricting the domain of $\kappa_j$ appropriately.

Thus from \eqref{eq:prop-main-for-a-j},
\begin{align}
	\left( d_\mathrm{\subspacetext} \left( \spn(\mathbf{u}_J) , \spn(\widetilde{\mathbf{u}}_{\widetilde{J}}) \right) \right)^2 & ~=~
 \frac{1}{q}\, \sum_{j\in J \atop j'\in {\widetilde{J}}^c} ~\left| \mathbf{u}_j^\conjsym \widetilde{\mathbf{u}}_{{j'}}\right|^2 
 \qquad \text{\small (due to Lemma~\ref{lemma:d-subspace-equivalent-forms}.2.)} \nonumber \\
& \qquad \leq~~ \frac{1}{q}\, \mathlarger{\mathlarger{\sum}}_{j\in J} ~\frac{ 
	\left\| (\widetilde{M} - M) \mathbf{u}_{{j}} \right\|^2_2 ~-~ \kappa_j \,\min_{j'\in \widetilde{J}}\left|\widetilde{\lambda}_{{j'}} - \lambda_j\right|^2 
}{ 
	\min_{j'\in {\widetilde{J}}^c} \left|\widetilde{\lambda}_{{j'}} - \lambda_j\right|^2 - \kappa_j \,\min_{j'\in \widetilde{J}}\left|\widetilde{\lambda}_{{j'}} - \lambda_j\right|^2 
}
\end{align}
for any $\kappa_j \in \Big[ 0,\min\left(1, \left( {\small \frac{ \min_{j'\in {\widetilde{J}}^c} \left|\widetilde{\lambda}_{{j'}} - \lambda_j\right|^2}{ \min_{j'\in \widetilde{J}}\left|\widetilde{\lambda}_{{j'}} - \lambda_j\right|^2}} \right)_{\!\!-} \, \right) \Big], ~j\in J$.

\vspace{0.5em}
Additionally, we note that
\begin{align*} &\| (\widetilde{M} - M) \mathbf{u}_{{j}} \|_2 < \min_{j'\in {\widetilde{J}}^c} |\widetilde{\lambda}_{{j'}} - \lambda_j|
~~\Rightarrow~~
\min_{j'\in {\widetilde{J}}^c} |\widetilde{\lambda}_{{j'}} - \lambda_j| > \min_{j'\in \widetilde{J}} |\widetilde{\lambda}_{{j'}} - \lambda_j| \\
&\qquad\qquad\qquad\qquad\qquad\qquad \text{\small (since, due to Corollary~\ref{cor:M-diff-lambda-diff}, $\| (\widetilde{M} - M) \mathbf{u}_{{j}} \|_2 \geq \min_{j'\in {N}} |\widetilde{\lambda}_{{j'}} - \lambda_j|$.)}
\end{align*}
Thus, when $\| (\widetilde{M} - M) \mathbf{u}_{{j}} \|_2 < \min_{j'\in {\widetilde{J}}^c} |\widetilde{\lambda}_{{j'}} - \lambda_j|$ the valid domain of $\kappa_j$ is $[0,1]$.
%
The statement about the tightest bound then follows from the fact that the function $f(\kappa) = \frac{a-\kappa c}{b-\kappa c}, \kappa \in [0,d]$ (with $d < \frac{b}{c}$) is minimized with $\kappa=0$ when $a\geq b$, and with $\kappa=d$ when $a < b$.

\end{proof}

The key achievement in the above proposition is to provide an upper bound on the distance (in terms of $d_\mathrm{\subspacetext}$) between the invariant subspaces $\spn(\mathbf{u}_J)$ and $\spn(\widetilde{\mathbf{u}}_{\widetilde{J}})$ in terms of the distance between the matrices $M$ and $\widetilde{M}$ and their eigenvalues.
For a given/fixed matrix perturbation, $(\widetilde{M} - M)$, and an appropriately chosen $\widetilde{J}$, the inequality \eqref{eq:prop-full-main} 
can be interpreted as a relation between the perturbation in the eigenvalues, $\{\lambda_j\,|\, j\in J\}$, 
and the perturbation in the invariant space $\spn(\mathbf{u}_J)$.
This relationship, 
in general, can be expected to be an inverse one -- with higher perturbation in the eigenvalues we will have a lower (upper-bound on the) perturbation in the invariant space, and vice versa.

It is easy to note that equality in \eqref{eq:prop-full-main} holds when 
\begin{itemize}
\item[i.] $\frac{ \min_{j'\in {\widetilde{J}}^c} \left|\widetilde{\lambda}_{{j'}} - \lambda_j\right|^2}{ \min_{j'\in \widetilde{J}}\left|\widetilde{\lambda}_{{j'}} - \lambda_j\right|^2} > 1, \,\forall j\in J$, allowing us to choose $\kappa_j=1, \,\forall j\in J$, and,
\item[ii.] $\min_{j'\in \widetilde{J}}\left|\widetilde{\lambda}_{{j'}} - \lambda_{j_1}\right| = \min_{j'\in \widetilde{J}}\left|\widetilde{\lambda}_{{j'}} - \lambda_{j_2}\right|, \,\forall j_1,j_2\in J$, \newline \phantom{.}\qquad $\min_{j'\in {\widetilde{J}}^c}\left|\widetilde{\lambda}_{{j'}} - \lambda_{j_1}\right| = \min_{j'\in {\widetilde{J}}^c}\left|\widetilde{\lambda}_{{j'}} - \lambda_{j_2}\right|, \,\forall j_1,j_2\in {J^c}$ \newline (these conditions hold, for example, when $\widetilde{\lambda}_{\widetilde{J}}$ and $\widetilde{\lambda}_{{\widetilde{J}}^c}$ are small translations of ${\lambda}_{{J}}$ and ${\lambda}_{{{J^c}}}$ respectively in $\mathbb{C}$.)
\end{itemize}

In Proposition~\ref{prop:full-main}, without loss of generality, we can interchange the roles of $J$ and ${J^c}$ (likewise $\widetilde{J}$ and ${\widetilde{J}}^c$). Observing that $\spn(\mathbf{u}_{J^c})$ and $\spn(\widetilde{\mathbf{u}}_{{\widetilde{J}}^c})$ are $(n-q)$ dimensional sub-spaces of $\mathbb{C}^n$ which are orthogonal complements of $\spn(\mathbf{u}_J)$ and $\spn(\widetilde{\mathbf{u}}_{\widetilde{J}})$ respectively, we then obtain
\begin{align} \label{eq:prop-full-main-J{J^c}-interchanged}
d_\mathrm{\subspacetext} \left( \spn(\mathbf{u}_J) , \spn(\widetilde{\mathbf{u}}_{\widetilde{J}}) \right) 
& =~~ \sqrt{\frac{n-q}{q}} ~ d_\mathrm{\subspacetext} \left( \spn(\mathbf{u}_{J^c}) , \spn(\widetilde{\mathbf{u}}_{{\widetilde{J}}^c}) \right) ~~~~\text{\small (due to Lemma~\ref{lemma:d-subspace-properties}.3)} \nonumber \\
& ~~\leq~~ \sqrt{ \frac{1}{q}\, \mathlarger{\mathlarger{\sum}}_{j\in {J^c}} ~\frac{ 
	\left\| (\widetilde{M} - M) \mathbf{u}_{{j}} \right\|^2_2 ~-~ \kappa_j \,\min_{j'\in {\widetilde{J}}^c}\left|\widetilde{\lambda}_{{j'}} - \lambda_j\right|^2 
}{ 
	\min_{j'\in \widetilde{J}} \left|\widetilde{\lambda}_{{j'}} - \lambda_j\right|^2 - \kappa_j \,\min_{j'\in {\widetilde{J}}^c}\left|\widetilde{\lambda}_{{j'}} - \lambda_j\right|^2 
}}
\end{align}
for any $\kappa_j \in \Big[ 0,\min\left(1, \left( {\small \frac{ \min_{j'\in \widetilde{J}} \left|\widetilde{\lambda}_{{j'}} - \lambda_j\right|^2}{ \min_{j'\in {\widetilde{J}}^c}\left|\widetilde{\lambda}_{{j'}} - \lambda_j\right|^2}} \right)_{\!\!-} \, \right) \Big], ~j\in {J^c}$.

\begin{cor}~ \label{cor:main-simplified}
For any ~$\kappa_J \in \Big[ 0,\min\left(1, \left( {\small \frac{ \min_{j\in J \atop j'\in {\widetilde{J}}^c} \left|\widetilde{\lambda}_{{j'}} - \lambda_j\right|^2}{ \max_{j\in J} \min_{j'\in \widetilde{J}}\left|\widetilde{\lambda}_{{j'}} - \lambda_j\right|^2}} \right)_{\!\!-} \, \right) \Big]$~
and

\qquad\qquad\qquad\qquad\qquad $\kappa_{J^c} \in \Big[ 0,\min\left(1, \left( {\small \frac{ \min_{j\in {J^c} \atop j'\in \widetilde{J}} \left|\widetilde{\lambda}_{{j'}} - \lambda_j\right|^2}{ \max_{j\in {J^c}} \min_{j'\in {\widetilde{J}}^c}\left|\widetilde{\lambda}_{{j'}} - \lambda_j\right|^2}} \right)_{\!\!-} \, \right) \Big]$,
\begin{itemize}
\Item[1.]
{\small \begin{align} \label{eq:corr-prop-main-simplified}
	d_\mathrm{\subspacetext} \left( \spn(\mathbf{u}_J) , \spn(\widetilde{\mathbf{u}}_{\widetilde{J}}) \right)  
	& \leq~~
	\frac{1}{q} \min \left( \sqrt{ \frac{ 
		\sum_{j\in J} \left\| (\widetilde{M} - M) \mathbf{u}_{{j}} \right\|^2_2 ~-~ \kappa_J \sum_{j\in J} \,\min_{j'\in \widetilde{J}}\left|\widetilde{\lambda}_{{j'}} - \lambda_j\right|^2 
	}{ 
		\min_{j\in J \atop j'\in {\widetilde{J}}^c} \left|\widetilde{\lambda}_{{j'}} - \lambda_j\right|^2 - \kappa_J \max_{j\in J} \,\min_{j'\in \widetilde{J}}\left|\widetilde{\lambda}_{{j'}} - \lambda_j\right|^2
	} } ~, \right. \nonumber \\
& \qquad\qquad\quad \left. \sqrt{ \frac{ 
	\sum_{j\in {J^c}} \left\| (\widetilde{M} - M) \mathbf{u}_{{j}} \right\|^2_2 ~-~ \kappa_{J^c} \sum_{j\in {J^c}} \,\min_{j'\in {\widetilde{J}}^c}\left|\widetilde{\lambda}_{{j'}} - \lambda_j\right|^2 
}{ 
	\min_{j\in {J^c} \atop j'\in \widetilde{J}} \left|\widetilde{\lambda}_{{j'}} - \lambda_j\right|^2 - \kappa_{J^c} \max_{j\in {J^c}} \,\min_{j'\in {\widetilde{J}}^c}\left|\widetilde{\lambda}_{{j'}} - \lambda_j\right|^2
} } \right)
	\end{align}}

\Item[2.]
{\small \begin{align} \label{eq:prop-frob}
	& d_\mathrm{\subspacetext} \left( \spn(\mathbf{u}_J) , \spn(\widetilde{\mathbf{u}}_{\widetilde{J}}) \right)
	\nonumber \\
	& ~~\leq \sqrt{ \frac{ \displaystyle
	\frac{1}{q} \left( \left\| \widetilde{M} - M \right\|^2_F 
		~-~ \left( \kappa_J \sum_{j\in J} \min_{{j'} \in \widetilde{J}}\left|\widetilde{\lambda}_{{j'}} - \lambda_j\right|^2 + \kappa_{J^c} \sum_{j\in {J^c}} \min_{{j'} \in {\widetilde{J}}^c}\left|\widetilde{\lambda}_{{j'}} - \lambda_j\right|^2 \right) \right)
	}{  \displaystyle 
		\left(  \min_{{j'} \in {\widetilde{J}}^c \atop j\in J} \left|\widetilde{\lambda}_{{j'}} \!-\! \lambda_j\right|^2 \!\!+\! \min_{{j'} \in \widetilde{J} \atop j\in {J^c}} \left|\widetilde{\lambda}_{{j'}} \!-\! \lambda_j\right|^2 \right)
		\!-\! \left( \kappa_J \, \max_{j\in J} \,\min_{{j'} \in \widetilde{J}}\left|\widetilde{\lambda}_{{j'}} \!-\! \lambda_j\right|^2 \!\!+\, \kappa_{J^c} \, \max_{j\in {J^c}} \,\min_{{j'} \in {\widetilde{J}}^c}\left|\widetilde{\lambda}_{{j'}} \!-\! \lambda_j\right|^2 \right) }  }
\end{align}}
\end{itemize} 
\end{cor}

\begin{proof}
With $\kappa_j \in \Big[ 0,\min\left(1, \left( {\small \frac{ \min_{j'\in {\widetilde{J}}^c} \left|\widetilde{\lambda}_{{j'}} - \lambda_j\right|^2}{ \min_{j'\in \widetilde{J}}\left|\widetilde{\lambda}_{{j'}} - \lambda_j\right|^2}} \right)_{\!\!-} \, \right) \Big], ~j\in J$,
\begin{align} 
	& q~ \left( d_\mathrm{\subspacetext} \left( \spn(\mathbf{u}_J) , \spn(\widetilde{\mathbf{u}}_{\widetilde{J}}) \right) \right)^2 \nonumber \\
	& \leq~~  \mathlarger{\mathlarger{\sum}}_{j\in J} ~\frac{ 
		\left\| (\widetilde{M} - M) \mathbf{u}_{{j}} \right\|^2_2 ~-~ \kappa_j \,\min_{j'\in \widetilde{J}}\left|\widetilde{\lambda}_{{j'}} - \lambda_j\right|^2 
	}{ 
		\min_{j'\in {\widetilde{J}}^c} \left|\widetilde{\lambda}_{{j'}} - \lambda_j\right|^2 - \kappa_j \,\min_{j'\in \widetilde{J}}\left|\widetilde{\lambda}_{{j'}} - \lambda_j\right|^2 
	}
	~~~~\text{\small (due to Proposition~\ref{prop:full-main})} \nonumber \\
	& \leq~~ \frac{ 
		\sum_{j\in J} \left\| (\widetilde{M} - M) \mathbf{u}_{{j}} \right\|^2_2 ~-~ \sum_{j\in J} \kappa_j \,\min_{j'\in \widetilde{J}}\left|\widetilde{\lambda}_{{j'}} - \lambda_j\right|^2 
	}{ 
		\min_{j\in J} \left( \min_{j'\in {\widetilde{J}}^c} \left|\widetilde{\lambda}_{{j'}} - \lambda_j\right|^2 - \kappa_j \,\min_{j'\in \widetilde{J}}\left|\widetilde{\lambda}_{{j'}} - \lambda_j\right|^2 \right)
	} 
	~~~~\text{\small (since $\sum_{k\in S} \frac{c_k}{d_k} \leq \frac{\sum_{k\in S} c_k}{\min_{k\in S} d_k}$)} \nonumber \\
	& \leq~~ \frac{ 
		\sum_{j\in J} \left\| (\widetilde{M} - M) \mathbf{u}_{{j}} \right\|^2_2 ~-~ \sum_{j\in J} \kappa_j \,\min_{j'\in \widetilde{J}}\left|\widetilde{\lambda}_{{j'}} - \lambda_j\right|^2 
	}{ 
		 \min_{j\in J \atop j'\in {\widetilde{J}}^c} \left|\widetilde{\lambda}_{{j'}} - \lambda_j\right|^2 - \max_{j\in J} \kappa_j \,\min_{j'\in \widetilde{J}}\left|\widetilde{\lambda}_{{j'}} - \lambda_j\right|^2
	} 
	~~~~\text{\small ($\min_{k\in S} (c_k - d_k) \geq \min_{k\in S} c_k - \max_{k\in S} d_k$.)}
\end{align}
We next choose $\kappa_j = \kappa_k, \,\forall j,k \in J$ and denote this value by
\[ \textstyle
\kappa_J ~\in~ \mathlarger{\mathlarger{\bigcap}}_{j\in J} \Big[ 0,\min\left(1, \left( {\small \frac{ \min_{j'\in {\widetilde{J}}^c} \left|\widetilde{\lambda}_{{j'}} - \lambda_j\right|^2}{ \min_{j'\in \widetilde{J}}\left|\widetilde{\lambda}_{{j'}} - \lambda_j\right|^2}} \right)_{\!\!-} \, \right) \Big]
~\supseteq~
\Big[ 0,\min\left(1, \left( {\small \frac{ \min_{j\in J \atop j'\in {\widetilde{J}}^c} \left|\widetilde{\lambda}_{{j'}} - \lambda_j\right|^2}{ \max_{j\in J} \min_{j'\in \widetilde{J}}\left|\widetilde{\lambda}_{{j'}} - \lambda_j\right|^2}} \right)_{\!\!-} \, \right) \Big]
\]
Thus,
\begin{align}  \label{eq:d-sp-simplified-J}
q~ \left( d_\mathrm{\subspacetext} \left( \spn(\mathbf{u}_J) , \spn(\widetilde{\mathbf{u}}_{\widetilde{J}}) \right) \right)^2 
~~\leq~~
\frac{ 
	\sum_{j\in J} \left\| (\widetilde{M} - M) \mathbf{u}_{{j}} \right\|^2_2 ~-~ \kappa_J \sum_{j\in J} \,\min_{j'\in \widetilde{J}}\left|\widetilde{\lambda}_{{j'}} - \lambda_j\right|^2 
}{ 
	\min_{j\in J \atop j'\in {\widetilde{J}}^c} \left|\widetilde{\lambda}_{{j'}} - \lambda_j\right|^2 - \kappa_J \max_{j\in J} \,\min_{j'\in \widetilde{J}}\left|\widetilde{\lambda}_{{j'}} - \lambda_j\right|^2
} 
\end{align}
for any $\kappa_J \in \Big[ 0,\min\left(1, \left( {\small \frac{ \min_{j\in J \atop j'\in {\widetilde{J}}^c} \left|\widetilde{\lambda}_{{j'}} - \lambda_j\right|^2}{ \max_{j\in J} \min_{j'\in \widetilde{J}}\left|\widetilde{\lambda}_{{j'}} - \lambda_j\right|^2}} \right)_{\!\!-} \, \right) \Big]$.

By interchanging the roles of $J$ and ${J^c}$ (accordingly, $\widetilde{J}$ and ${\widetilde{J}}^c$), and noting that $\spn(\mathbf{u}_{J^c})$ and $\spn(\widetilde{\mathbf{u}}_{{\widetilde{J}}^c})$ are $(n-q)$ dimensional sub-spaces of $\mathbb{C}^n$, we get
\begin{align} \label{eq:d-sp-{J^c}}
(n-q)\, \left( d_\mathrm{\subspacetext} \left( \spn(\mathbf{u}_{J^c}) , \spn(\widetilde{\mathbf{u}}_{{\widetilde{J}}^c}) \right) \right)^2 
~~\leq~~
\frac{ 
	\sum_{j\in {J^c}} \left\| (\widetilde{M} - M) \mathbf{u}_{{j}} \right\|^2_2 ~-~ \kappa_{J^c} \sum_{j\in {J^c}} \,\min_{j'\in {\widetilde{J}}^c}\left|\widetilde{\lambda}_{{j'}} - \lambda_j\right|^2 
}{ 
	\min_{j\in {J^c} \atop j'\in \widetilde{J}} \left|\widetilde{\lambda}_{{j'}} - \lambda_j\right|^2 - \kappa_{J^c} \max_{j\in {J^c}} \,\min_{j'\in {\widetilde{J}}^c}\left|\widetilde{\lambda}_{{j'}} - \lambda_j\right|^2
} 
\end{align}
for any $\kappa_{J^c} \in \Big[ 0,\min\left(1, \left( {\small \frac{ \min_{j\in {J^c} \atop j'\in \widetilde{J}} \left|\widetilde{\lambda}_{{j'}} - \lambda_j\right|^2}{ \max_{j\in {J^c}} \min_{j'\in {\widetilde{J}}^c}\left|\widetilde{\lambda}_{{j'}} - \lambda_j\right|^2}} \right)_{\!\!-} \, \right) \Big]$.

But, since $\spn(\mathbf{u}_J)$ and $\spn(\mathbf{u}_{J^c})$ are orthogonal complements (likewise, $\spn(\mathbf{u}_{\widetilde{J}})$ and $\spn(\mathbf{u}_{{\widetilde{J}}^c})$ are orthogonal complements), using Lemma~\ref{lemma:d-subspace-properties} we can write \eqref{eq:d-sp-{J^c}} as
\begin{align} \label{eq:d-sp-simplified-J2}
q\, \left( d_\mathrm{\subspacetext} \left( \spn(\mathbf{u}_J) , \spn(\widetilde{\mathbf{u}}_{\widetilde{J}}) \right) \right)^2  
~~\leq~~
\frac{ 
	\sum_{j\in {J^c}} \left\| (\widetilde{M} - M) \mathbf{u}_{{j}} \right\|^2_2 ~-~ \kappa_{J^c} \sum_{j\in {J^c}} \,\min_{j'\in {\widetilde{J}}^c}\left|\widetilde{\lambda}_{{j'}} - \lambda_j\right|^2 
}{ 
	\min_{j\in {J^c} \atop j'\in \widetilde{J}} \left|\widetilde{\lambda}_{{j'}} - \lambda_j\right|^2 - \kappa_{J^c} \max_{j\in {J^c}} \,\min_{j'\in {\widetilde{J}}^c}\left|\widetilde{\lambda}_{{j'}} - \lambda_j\right|^2
} 
\end{align}
for any $\kappa_{J^c} \in \Big[ 0,\min\left(1, \left( {\small \frac{ \min_{j\in {J^c} \atop j'\in \widetilde{J}} \left|\widetilde{\lambda}_{{j'}} - \lambda_j\right|^2}{ \max_{j\in {J^c}} \min_{j'\in {\widetilde{J}}^c}\left|\widetilde{\lambda}_{{j'}} - \lambda_j\right|^2}} \right)_{\!\!-} \, \right) \Big]$.

Combining \eqref{eq:d-sp-simplified-J} and \eqref{eq:d-sp-simplified-J2} proves part `1.'

Again, adding \eqref{eq:d-sp-simplified-J} and \eqref{eq:d-sp-simplified-J2},
\begin{align*}
& q~ \left( 
\min_{j\in J \atop j'\in {\widetilde{J}}^c} \left|\widetilde{\lambda}_{{j'}} - \lambda_j\right|^2 
+ \min_{j\in {J^c} \atop j'\in \widetilde{J}} \left|\widetilde{\lambda}_{{j'}} - \lambda_j\right|^2 
- \kappa_J \max_{j\in J} \,\min_{j'\in \widetilde{J}}\left|\widetilde{\lambda}_{{j'}} - \lambda_j\right|^2
- \kappa_{J^c} \max_{j\in {J^c}} \,\min_{j'\in {\widetilde{J}}^c}\left|\widetilde{\lambda}_{{j'}} - \lambda_j\right|^2
\right) \\ 
& \qquad \times~ \left( d_\mathrm{\subspacetext} \left( \spn(\mathbf{u}_J) , \spn(\widetilde{\mathbf{u}}_{\widetilde{J}}) \right) \right)^2  \\
& 
\qquad\qquad\qquad
\leq~ 
\sum_{j\in J} \left\| (\widetilde{M} - M) \mathbf{u}_{{j}} \right\|^2_2 ~+~ \sum_{j\in {J^c}} \left\| (\widetilde{M} - M) \mathbf{u}_{{j}} \right\|^2_2 \\
& \qquad\qquad\qquad\qquad\qquad\qquad ~-~ \kappa_J \sum_{j\in J} \,\min_{j'\in \widetilde{J}}\left|\widetilde{\lambda}_{{j'}} - \lambda_j\right|^2 
 ~-~ \kappa_{J^c} \sum_{j\in {J^c}} \,\min_{j'\in {\widetilde{J}}^c}\left|\widetilde{\lambda}_{{j'}} - \lambda_j\right|^2 
\end{align*}
%
The part `2.' of the result then follows 
by observing that
\[
\sum_{{j} \in {J}} \left\| (\widetilde{M} - M) \mathbf{u}_{j} \right\|^2_2 ~+~ \sum_{{j} \in {{J^c}}} \left\| (\widetilde{M} - M) \mathbf{u}_{j} \right\|^2_2 ~=~ \left\| (\widetilde{M} - M) U \right\|^2_F ~=~ \left\| \widetilde{M} - M \right\|^2_F
\]
\end{proof}

\begin{cor}[Davis-Kahan~\cite{10.1137/0707001} -- see Section VIII.3 of~\cite{bhatia1996matrix}] \label{cor:kappa-zero}~
	
\begin{itemize}
\item[1.]	~~$
	\displaystyle d_\mathrm{\subspacetext} \left( \spn(\mathbf{u}_J) , \spn(\widetilde{\mathbf{u}}_{\widetilde{J}}) \right) 
~\leq~ 
\frac{ \displaystyle \min \left( 1, \sqrt{ \frac{n-q}{q} } \right) }{ \displaystyle \max \left( \mathrm{sep}\left( \lambda_J, \widetilde{\lambda}_{{\widetilde{J}}^c} \right), ~\mathrm{sep}\left( \lambda_{J^c}, \widetilde{\lambda}_{\widetilde{J}} \right) \right) } \left\| \widetilde{M} - M \right\|_2
	$
\item[2.]
~~$
d_\mathrm{\subspacetext} \left( \spn(\mathbf{u}_J) , \spn(\widetilde{\mathbf{u}}_{\widetilde{J}}) \right) 
~\leq~ 
\frac{\displaystyle \frac{1}{\sqrt{q}}~ \left\| \widetilde{M} - M \right\|_F }{\displaystyle \sqrt{ \mathrm{sep}\left( \lambda_J, \widetilde{\lambda}_{{\widetilde{J}}^c} \right)^{\, 2} ~+~ \mathrm{sep}\left( \lambda_{J^c}, \widetilde{\lambda}_{\widetilde{J}} \right)^{\, 2} } }
\newline
\phantom{.}\qquad\qquad\qquad\qquad\qquad\qquad\qquad ~\leq~
\displaystyle 
\sqrt{ \frac{\displaystyle n/q }{\displaystyle \mathrm{sep}\left( \lambda_J, \widetilde{\lambda}_{{\widetilde{J}}^c} \right)^{\, 2} ~+~ \mathrm{sep}\left( \lambda_{J^c}, \widetilde{\lambda}_{\widetilde{J}} \right)^{\, 2} } } ~ \left\| \widetilde{M} - M \right\|_2
$
\end{itemize}
\end{cor}

\begin{proof}
	In \eqref{eq:d-sp-simplified-J}, setting $\kappa_J = 0$, we get
	\[ \left( d_\mathrm{\subspacetext} \left( \spn(\mathbf{u}_J) , \spn(\widetilde{\mathbf{u}}_{\widetilde{J}}) \right) \right)^2 ~\leq~
%
%
%
	\frac{\frac{1}{q} \sum_{j\in J} \left\| (\widetilde{M} - M) \mathbf{u}_{{j}} \right\|^2_2}{\mathrm{sep}\left( \lambda_J, \widetilde{\lambda}_{{\widetilde{J}}^c} \right)^{~2}} ~\leq~
	 \frac{ \left\| (\widetilde{M} - M) \right\|^2_2}{\mathrm{sep}\left( \lambda_J, \widetilde{\lambda}_{{\widetilde{J}}^c} \right)^{~2}} \]
	Interchanging the roles of the tilde and non-tilde terms in this result we analogously obtain
	\[ \left( d_\mathrm{\subspacetext} \left( \spn(\widetilde{\mathbf{u}}_{\widetilde{J}}), \spn(\mathbf{u}_J) \right) \right)^2 \leq \frac{ \left\| (\widetilde{M} - M) \right\|^2_2}{\mathrm{sep}\left( \widetilde{\lambda}_{\widetilde{J}}, \lambda_{J^c} \right)^{~2}}\]
	The above two together gives
	\begin{align} \label{eq:kappa-free-eqJ}
	\left( d_\mathrm{\subspacetext} \left( \spn(\widetilde{\mathbf{u}}_{\widetilde{J}}), \spn(\mathbf{u}_J) \right) \right)^2  ~~\leq~~ 
	\frac{\left\| (\widetilde{M} - M) \right\|^2_2}{ \max\left( \mathrm{sep}\left( \lambda_J, \widetilde{\lambda}_{{\widetilde{J}}^c} \right), ~\mathrm{sep}\left( \lambda_{J^c}, \widetilde{\lambda}_{\widetilde{J}} \right) \right)^{~2}}
	\end{align}
	
	In the above inequality, interchanging the roles of $J$ and ${J^c}$ (accordingly, $\widetilde{J}$ and ${\widetilde{J}}^c$),
	and observing that by Lemma~\ref{lemma:d-subspace-properties} 
	$d_\mathrm{\subspacetext} \left( \spn(\widetilde{\mathbf{u}}_{\widetilde{J}}), \spn(\mathbf{u}_J) \right) = \sqrt{\frac{n-q}{q}} d_\mathrm{\subspacetext} \left( \spn(\widetilde{\mathbf{u}}_{{\widetilde{J}}^c}), \spn(\mathbf{u}_{J^c}) \right)$, we obtain
	\begin{align}\label{eq:kappa-free-eq{J^c}}
	\left( d_\mathrm{\subspacetext} \left( \spn(\widetilde{\mathbf{u}}_{\widetilde{J}}), \spn(\mathbf{u}_J) \right) \right)^2 ~~\leq~~ \frac{n-q}{q} \,
	\frac{\left\| (\widetilde{M} - M) \right\|^2_2}{ \max\left( \mathrm{sep}\left( \lambda_J, \widetilde{\lambda}_{{\widetilde{J}}^c} \right), ~\mathrm{sep}\left( \lambda_{J^c}, \widetilde{\lambda}_{\widetilde{J}} \right) \right)^{~2}}
	\end{align}
	\eqref{eq:kappa-free-eqJ} and \eqref{eq:kappa-free-eq{J^c}} together concludes the proof of part `1.'
	
	
	The second result follows directly from part `2.' of Corollary~\ref{cor:main-simplified} by setting $\kappa_J = \kappa_{J^c} = 0$ and using the fact that for $Q\in\mathbb{C}^{n \times n}$, $ \|Q\|_F \leq \sqrt{n} \|Q\|_2$.
\end{proof}



\subsection{Bound on Perturbation of Invariant Subspace of a Normal Matrix with Well-clustered Spectrum} \label{sec:well-clusterted-spectrum}


In this section we specialize the earlier results for the situation when $\lambda_J$ and $\lambda_{J^c}$ are well-clustered (\emph{i.e.}, the separation between them is large) compared to the perturbation $(\widetilde{M} - M)$.
In the following Lemma we outline the conditions under which the the perturbed eigenvalues, $\widetilde{\lambda}_N$ will also remain well-clustered.

\begin{lemma} \label{lemma:sep-preserving-M-tilde}
For any $J \subseteq N$, 
define ${J^c}\!=\!N\!-\!J$.
If $\|\widetilde{M} - M\|_2 < \frac{1}{2} \mathrm{sep}(\lambda_J, \lambda_{J^c})$ then,

\begin{itemize}
	\item[1.] $\widetilde{\lambda}_N$ is a separation-preserving perturbation of $\lambda_J$ and $\lambda_{J^c}$.
More explicitly, defining
\begin{equation} \label{eq:def-J-{J^c}-hat}
\widehat{J} = \{j' \,|\, \min_{j\in N} |\widetilde{\lambda}_{j'} - \lambda_j| = \min_{j\in J} |\widetilde{\lambda}_{j'} - \lambda_j| \} ~~~~\text{and}~~~~ 
{\widehat{J}^c} = \{j' \,|\, \min_{j\in N} |\widetilde{\lambda}_{j'} - \lambda_j| = \min_{j\in {J^c}} |\widetilde{\lambda}_{j'} - \lambda_j| \}
\end{equation}
makes
$\{ \widetilde{\lambda}_{\widehat{J}}, \widetilde{\lambda}_{{\widehat{J}^c}} \}$ a separation-preserving partition of $\widetilde{\lambda}_N$, with \[ \mathrm{sep}(\widetilde{\lambda}_{\widehat{J}}, \widetilde{\lambda}_{{\widehat{J}^c}}) ~>~ \mathrm{sep}(\lambda_J, \lambda_{J^c}) - 2 \|\widetilde{M} - M\|_2 
\]

\item[2.] $|\widetilde{\lambda}_{\widehat{J}}| = |\lambda_J|$ ~(equivalently, $|\widetilde{\lambda}_{{\widehat{J}^c}}| = |\lambda_{J^c}|$), where $|\cdot|$ denotes the number of elements in the multi-sets (recall that $\lambda_J$ and $\widetilde{\lambda}_{\widehat{J}}$ are multi-sets, allowing them to contain multiple copies of non-distinct eigenvalues, if any, of $M$ and $\widetilde{M}$ respectively). 

\end{itemize}
\end{lemma}

\begin{proof} ~

\begin{itemize}
	\item[1.] We first observe that 
	\begin{equation} \label{eq:M-Mtilde-dH}
	\|\widetilde{M} - M\|_2 ~\geq~ \max \left( \max_{j\in N} \,\min_{j'\in N} \,\left|\widetilde{\lambda}_{j'} - \lambda_j\right| ~,~ \max_{j'\in N} \,\min_{j\in N} \,\left|\widetilde{\lambda}_{j'} - \lambda_j\right| \right) = d_H(\lambda_{N}, \widetilde{\lambda}_N)
	\end{equation}
	As a consequence, $d_H(\lambda_{N}, \widetilde{\lambda}_N) \leq \|\widetilde{M} - M\|_2 < \frac{1}{2} \mathrm{sep}(\lambda_J, \lambda_{J^c})$.
	Then the proof of the first part follows directly from Corollary~\ref{corr:sep-preserving-purturbation} by setting $P = \lambda_J$, $Q = \lambda_{J^c}$ and $\widetilde{R} = \widetilde{\lambda}_N$.
	
	\item[2.] 
	We prove the second part by contradiction.
	
	\noindent \hspace{-0.5em}
	\begin{tabular}{p{0.55\columnwidth}p{0.35\columnwidth}}
	If possible, let $|\widetilde{\lambda}_{\widehat{J}}| \neq |\lambda_J|$. Without loss of generality we will assume $|\widetilde{\lambda}_{\widehat{J}}| < |\lambda_J|$ (if the $|\widetilde{\lambda}_{\widehat{J}}| > |\lambda_J|$, we can show the contradiction for $|\widetilde{\lambda}_{{\widehat{J}^c}}| < |\lambda_{J^c}|$ instead).
	
	Define a path $\overline{M}:[0,1]\rightarrow \mathbb{R}^{n\times n}$ connecting $M$ and $\widetilde{M}$ as \[ \overline{M}(t) ~=~ t \widetilde{M} ~+~ (1-t)M \]
	&
	\vspace{-1em}
	{\includegraphics[width=0.35\columnwidth, trim=0 0 0 0, clip=true]{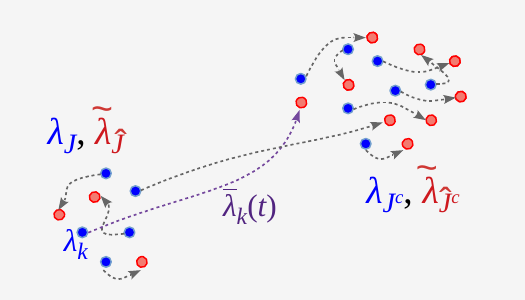}}
	\end{tabular}
	
	\vspace{-2em} \noindent
	Although $\overline{M}(t)$ is not necessarily normal for all $t$, its characteristic equation is a degree-$n$ polynomial equation in its eigenvalue with coefficient of the highest degree term equal to $1$ and other coefficients being polynomials in $t$.
	Since the roots of such a polynomial are continuous functions of the coefficients, the eigenvalues of $\widetilde{M}(t)$ are continuous functions of $t$. 
	Thus, we define $\overline{\lambda}_j:[0,1]\rightarrow \mathbb{C}$ to be the paths of the eigenvalues such that $\overline{\lambda}_j(0) = \lambda_j$ for all $j\in \{1,2,\cdots,n\}$.
	$\overline{\lambda}_j(1)$ are the eigenvalues of $\overline{M}(1) = \widetilde{M}$, so that $\overline{\lambda}_j(1) = \widetilde{\lambda}_{\sigma(j)}$ for some permutation $\sigma:\{1,2,\cdots,n\}\rightarrow \{1,2,\cdots,n\}$.
	
	Since $|\widetilde{\lambda}_{\widehat{J}}| < |\lambda_J|$, there exists at least one $k\in J$ 
	(with $\lambda_k = \overline{\lambda}_k(0) \in \lambda_J$)
	such that $\overline{\lambda}_k(1) \notin \widetilde{\lambda}_{\widehat{J}}$ (equivalently, $\overline{\lambda}_k(1) \in \widetilde{\lambda}_{{\widehat{J}^c}}$).
	
%
	Define $g(t) = \min_{j\in J} |\overline{\lambda}_k(t) - \lambda_j|$ and $h(t) = \min_{j\in {J^c}} |\overline{\lambda}_k(t) - \lambda_j|$. 
	Thus,
	\begin{align*}
	g(0) ~=~ \min_{j\in J} |\overline{\lambda}_k(0) - \lambda_j| &~=~ \min_{j\in J} |\lambda_k - \lambda_j| \\
	&~=~ 0 \quad\text{\small (since $\lambda_k \in \lambda_J$)} \\
	&~\leq h(0)
	\end{align*}
	
	
	Again,
	\begin{align*}
	h(1) &~=~ \min_{j\in {J^c}} |\overline{\lambda}_k(1) - \lambda_j|  \\
		&~\leq~ \min_{j\in J} |\overline{\lambda}_k(1) - \lambda_j| \qquad\text{\small (since $\overline{\lambda}_k(1) \in \widetilde{\lambda}_{{\widehat{J}^c}}$, from definition of ${\widehat{J}^c}$,} \\ & \qquad\qquad\qquad\qquad\qquad\qquad\text{\small $\min_{j\in {J^c}} |\overline{\lambda}_k(1) - \lambda_j| = \min_{j\in N} |\overline{\lambda}_k(1) - \lambda_j| $)} \\
		&~=~ g(1)
	\end{align*}
	Thus, by intermediate value theorem, there exists a ~$t'\in[0,1]$ such that $g(t') = h(t')$.
	That is, $\min_{j\in J} |\overline{\lambda}_k(t') - \lambda_j| = \min_{j\in {J^c}} |\overline{\lambda}_k(t') - \lambda_j|$. Equivalently,
	\begin{align} \label{eq:sep-equal}
	\mathrm{sep}(\lambda_J,\{\overline{\lambda}_k(t')\}) = \mathrm{sep}(\lambda_{J^c},\{\overline{\lambda}_k(t')\}),~~\text{for some}~ t'\in[0,1]
	\end{align}
	
	Now, 
	\begin{align}
	\|\overline{M}(t') - M\|_2 & ~~\geq~~ \min_{j\in N} |\overline{\lambda}_k(t') - \lambda_j| \qquad\text{\small (Corollary~\ref{cor:M-diff-lambda-diff}.1)} \nonumber \\
							& ~~=~~ \min\left( \min_{j\in J} |\overline{\lambda}_k(t') - \lambda_j| , \min_{j\in {J^c}} |\overline{\lambda}_k(t') - \lambda_j|  \right)  \nonumber \\
							& ~~=~~ \frac{1}{2}\left( \min_{j\in J} |\overline{\lambda}_k(t') - \lambda_j| + \min_{j\in {J^c}} |\overline{\lambda}_k(t') - \lambda_j|\right)  \nonumber \\ & \qquad\qquad\qquad\qquad\qquad\text{\small (since from \eqref{eq:sep-equal}, $\min_{j\in J} |\overline{\lambda}_k(t') - \lambda_j| = \min_{j\in {J^c}} |\overline{\lambda}_k(t') - \lambda_j|$)} \nonumber \\
							& ~~=~~ \frac{1}{2}\left( \mathrm{sep}(\lambda_J, \{\overline{\lambda}_k(t')\}) ~+~ \mathrm{sep}(\lambda_{J^c}, \{\overline{\lambda}_k(t')\}) ~+~ \mathrm{diam}( \{\overline{\lambda}_k(t')\}) \right) \nonumber \\ & \qquad\qquad\qquad\qquad\qquad\qquad\qquad\text{\small (since diameter of a point is zero.)} \nonumber \\
							& ~~\geq~~ \frac{1}{2} \mathrm{sep}(\lambda_J,\lambda_{J^c}) \qquad\text{\small (using Lemma~\ref{lemma:sep-diam-triangle-inequality})}
	\end{align}
	
	However, $\|\overline{M}(t') - M\|_2 ~=~ t' \|\widetilde{M} - M\|_2 ~<~ t' \frac{1}{2} \mathrm{sep}(\lambda_J,\lambda_{J^c}) ~\leq~ \frac{1}{2} \mathrm{sep}(\lambda_J,\lambda_{J^c})$.
	We thus end up with a contradiction.
%
\end{itemize}
\end{proof}



In the following propositions, we express the upper bounds on $d_\mathrm{\subspacetext} \left( \spn(\mathbf{u}_J) , \spn(\widetilde{\mathbf{u}}_{\widehat{J}}) \right)$ in terms of $(\widetilde{M}-M)$ and non-tilde terms only.

\begin{prop} \label{prop:subspace-purturbation-hat-tilde-free}
For any $J \subseteq N$ such that $|J|=q$, 
define ${J^c}\!=\!N\!-\!J$.
	If $\|\widetilde{M} - M\|_2 < \frac{1}{2} \mathrm{sep}(\lambda_J, \lambda_{J^c})$,
\begin{itemize}
\Item[1.]	
\begin{align} 
	d_\mathrm{\subspacetext} \left( \spn(\mathbf{u}_J) , \spn(\widetilde{\mathbf{u}}_{\widehat{J}}) \right)
	& ~~\leq~~
	\frac{1}{\sqrt{q}}\, \min \left(
	\sqrt{ \mathlarger{\mathlarger{\sum}}_{j\in J} \left(\frac{ 
		\left\| (\widetilde{M} - M) \mathbf{u}_{{j}} \right\|_2 
	}{ 
		 \min_{k\in {J^c}} |{\lambda}_{k} - \lambda_j| ~-~ \|\widetilde{M} - M\|_2  
	} \right)^{\!\!2} } ~,~ 
		\right. \nonumber \\ & \qquad\qquad \qquad\qquad  \left.
	\sqrt{ \mathlarger{\mathlarger{\sum}}_{j\in {J^c}} \left( \frac{ 
		\left\| (\widetilde{M} - M) \mathbf{u}_{{j}} \right\|_2 
	}{ 
		 \min_{k\in J} |{\lambda}_{k} - \lambda_j| ~-~ \|\widetilde{M} - M\|_2  
	} \right)^{\!\!2} } \right)
\\
	& \quad ~~\leq~~
	\min\left( 1, \sqrt{\frac{n-q}{q}} \right) ~ \frac{ 
		\left\| \widetilde{M} - M \right\|_2 
	}{ 
		 \mathrm{sep}(\lambda_J, \lambda_{J^c})  ~-~ \|\widetilde{M} - M\|_2  
	} 
	\end{align}
\Item[2.] ~

\begin{align}
d_\mathrm{\subspacetext} \left( \spn(\mathbf{u}_J) , \spn(\widetilde{\mathbf{u}}_{\widehat{J}}) \right)
& ~~\leq~~
\frac{\frac{1}{\sqrt{2q}}  \left\| \widetilde{M} - M \right\|_F}{\mathrm{sep}(\lambda_J, \lambda_{J^c})  ~-~ \|\widetilde{M} - M\|_2}
\end{align}
\end{itemize}
	where $\widehat{J}$ and ${\widehat{J}^c}$ are as defined in \eqref{eq:def-J-{J^c}-hat}.
\end{prop}

\begin{proof}
%
For any $j\in J$,
\begin{align}
\min_{j'\in {\widehat{J}^c}} |\widetilde{\lambda}_{j'} - \lambda_j| & ~~=~~ \mathrm{sep}(\{\lambda_j\}, \widetilde{\lambda}_{{\widehat{J}^c}}) \nonumber \\
&~~\geq~~
	\mathrm{sep}(\{\lambda_j\}, \lambda_{J^c}) ~-~ d_H(\lambda_{J^c}, \widetilde{\lambda}_{{\widehat{J}^c}}) \quad\text{\small (due to Lemma~\ref{lemma:distance-lemma})} \nonumber \\
& ~~\geq~~ \mathrm{sep}(\{\lambda_j\}, \lambda_{J^c}) ~-~ d_H(\lambda_N, \widetilde{\lambda}_N) \quad\text{\small (due to Lemma~\ref{lemma:sep-preserving-purturbation}.4., $d_H(\lambda_{J^c}, \widetilde{\lambda}_{{\widehat{J}^c}}) \leq d_H(\lambda_N, \widetilde{\lambda}_N)$.)} \nonumber \\
& ~~\geq~~ \mathrm{sep}(\{\lambda_j\}, \lambda_{J^c}) ~-~ \|\widetilde{M} - M\|_2 \quad\text{\small (using \eqref{eq:M-Mtilde-dH}.)} \nonumber \\
& ~~=~~ \min_{k\in {J^c}} |{\lambda}_{k} - \lambda_j| ~-~ \|\widetilde{M} - M\|_2 \label{eq:lambda-diff-inequality-hat-Mdiff}
\end{align}

Thus, 
in Proposition~\ref{prop:full-main} choosing $\kappa_j=0, \,\forall j\in J$, we get
\begin{align} 
\left( d_\mathrm{\subspacetext} \left( \spn(\mathbf{u}_J) , \spn(\widetilde{\mathbf{u}}_{\widehat{J}}) \right) \right)^2
& ~~\leq~~ \frac{1}{q}\, \mathlarger{\mathlarger{\sum}}_{j\in J} ~\frac{ 
	\left\| (\widetilde{M} - M) \mathbf{u}_{{j}} \right\|^2_2 
}{ 
	\min_{j'\in {\widehat{J}^c}} \left|\widetilde{\lambda}_{{j'}} - \lambda_j\right|^2 
} \nonumber \\
& ~~\leq~~ \frac{1}{q}\, \mathlarger{\mathlarger{\sum}}_{j\in J} ~\frac{ 
	\left\| (\widetilde{M} - M) \mathbf{u}_{{j}} \right\|^2_2 
}{ 
	\left( \min_{k\in {J^c}} |{\lambda}_{k} - \lambda_j| ~-~ \|\widetilde{M} - M\|_2  \right)^2
} 
\\
& ~~\leq~~  \frac{ 
	\frac{1}{q} \sum_{j\in J} \left\| (\widetilde{M} - M) \mathbf{u}_{{j}} \right\|^2_2 
}{ 
	\min_{j\in J} \left( \min_{k\in {J^c}} |{\lambda}_{k} - \lambda_j| ~-~ \|\widetilde{M} - M\|_2  \right)^2
} \nonumber \\ & \qquad\qquad\qquad\qquad\qquad\qquad\qquad\qquad\text{\small (since $\sum_{k\in S} \frac{c_k}{d_k} \leq \frac{\sum_{k\in S} c_k}{\min_{k\in S} d_k}$)} \nonumber \\
& ~~\leq~~  \frac{ 
	\left\| \widetilde{M} - M \right\|^2_2 
}{ 
	\left( \mathrm{sep}(\lambda_J, \lambda_{J^c})  ~-~ \|\widetilde{M} - M\|_2  \right)^2
} \label{eq:uJ-Jhat-inequality-2} \\ 
& \qquad\qquad\qquad\qquad\qquad \text{\small (since $\left\| \widetilde{M} - M \right\|_2 \geq \left\| (\widetilde{M} - M) \mathbf{u}_{{j}} \right\|_2$ ~and} \nonumber \\ & \qquad\qquad\qquad\qquad\qquad\qquad\qquad \text{$\min_{k\in J}(c_k - \alpha)^2 = (\min_{k\in J} c_k - \alpha)^2$)} \nonumber
\end{align}
%
In the above, switching the roles of $J$ and ${J^c}$ (likewise, $\widehat{J}$ and ${\widehat{J}^c}$), and noting that $\spn(\mathbf{u}_{J^c})$ and $\spn(\widetilde{\mathbf{u}}_{{\widehat{J}^c}})$ are $(n-q)$-dimensional subspaces of $\mathbb{C}^n$, we get
\begin{align}
\left( d_\mathrm{\subspacetext} \left( \spn(\mathbf{u}_{J^c}) , \spn(\widetilde{\mathbf{u}}_{{\widehat{J}^c}}) \right) \right)^2
& ~~\leq~~
\frac{1}{n-q}\, \mathlarger{\mathlarger{\sum}}_{j\in {J^c}} ~\frac{ 
	\left\| (\widetilde{M} - M) \mathbf{u}_{{j}} \right\|^2_2 
}{ 
	\left( \min_{k\in J} |{\lambda}_{k} - \lambda_j| ~-~ \|\widetilde{M} - M\|_2  \right)^2
} \nonumber \\
& \quad ~~\leq~~
 \frac{ 
	\left\| \widetilde{M} - M \right\|^2_2 
}{ 
	\left( \mathrm{sep}(\lambda_J, \lambda_{J^c})  ~-~ \|\widetilde{M} - M\|_2  \right)^2
}\nonumber
\end{align}
But since $\spn(\mathbf{u}_{J^c})$ and $\spn(\widetilde{\mathbf{u}}_{{\widehat{J}^c}})$ are orthogonal complements of $\spn(\mathbf{u}_J)$ and $\spn(\widetilde{\mathbf{u}}_{\widehat{J}})$ respectively, from Lemma~\ref{lemma:d-subspace-properties} we have $(n-q)\, ( d_\mathrm{\subspacetext} \left( \spn(\mathbf{u}_{J^c}) , \spn(\widetilde{\mathbf{u}}_{{\widehat{J}^c}}) \right) )^2  = q\, ( d_\mathrm{\subspacetext} \left( \spn(\mathbf{u}_J) , \spn(\widetilde{\mathbf{u}}_{\widehat{J}}) \right) )^2$.
This gives us from the above,
\begin{align}
\left( d_\mathrm{\subspacetext} \left( \spn(\mathbf{u}_J) , \spn(\widetilde{\mathbf{u}}_{\widehat{J}}) \right) \right)^2
& ~~\leq~~
\frac{1}{q}\, \mathlarger{\mathlarger{\sum}}_{j\in {J^c}} ~\frac{ 
	\left\| (\widetilde{M} - M) \mathbf{u}_{{j}} \right\|^2_2 
}{ 
	\left( \min_{k\in J} |{\lambda}_{k} - \lambda_j| ~-~ \|\widetilde{M} - M\|_2  \right)^2
} 
\\
& \quad ~~\leq~~
\frac{n-q}{q} \frac{ 
	\left\| \widetilde{M} - M \right\|^2_2 
}{ 
	\left( \mathrm{sep}(\lambda_J, \lambda_{J^c})  ~-~ \|\widetilde{M} - M\|_2  \right)^2
} \label{eq:uJ-Jhat-inequality-J{J^c}flipped-2}
\end{align}
Combining 
 \eqref{eq:uJ-Jhat-inequality-2} and \eqref{eq:uJ-Jhat-inequality-J{J^c}flipped-2} gives the first result of the proposition. 
 
 The second result can be obtained directly using Corollary~\ref{cor:kappa-zero}.2. and observing that due to 
 \eqref{eq:lambda-diff-inequality-hat-Mdiff}, $\mathrm{sep}(\lambda_J, \widetilde{\lambda}_{{\widehat{J}^c}}) \geq \mathrm{sep}(\lambda_J, \lambda_{J^c})  - \|\widetilde{M} - M\|_2$ (and analogously $\mathrm{sep}(\lambda_{J^c}, \widetilde{\lambda}_{\widehat{J}}) \geq \mathrm{sep}(\lambda_J, \lambda_{J^c})  - \|\widetilde{M} - M\|_2$).
\end{proof}

Assuming $q\leq n/2$,
it is worth noting that defining $\epsilon = \frac{1}{2} \mathrm{sep}(\lambda_J, \lambda_{J^c}) - \|\widetilde{M} - M\|_2$, the second inequality of the first result in the above proposition becomes $d_\mathrm{\subspacetext} \left( \spn(\mathbf{u}_J) , \spn(\widetilde{\mathbf{u}}_{\widehat{J}}) \right) \leq \frac{\|\widetilde{M} - M\|_2}{\|\widetilde{M} - M\|_2 ~+~ 2 \epsilon} $.
Thus, with $\epsilon \rightarrow 0$, this inequality becomes $d_\mathrm{\subspacetext} \left( \spn(\mathbf{u}_J) , \spn(\widetilde{\mathbf{u}}_{\widehat{J}}) \right) < 1$, rendering the result uninformative /redundant.
Thus, higher the separation between $\lambda_J$ and $\lambda_{J^c}$ (relative to $\|\widetilde{M} - M\|_2$), tighter will be the upper bound in the result of the proposition.

An interpretation of the result in the above proposition is that a perturbation, $\widetilde{M} - M$, of the matrix $M$ will result in perturbation in the invariant subspace $\spn(\mathbf{u}_J)$ such that the distance between the subspace and its purturbed counterpart is bounded by the upper bounds mentioned in the proposition.
One key feature of the proposition, however, is the the upper bound in the inequality does not depend on $\widehat{J}$. As a consequence, for any other size-$q$ subset, $\widetilde{J}$, of $N$ such that $\spn(\mathbf{u}_{\widetilde{J}})$ is closer to $\spn(\mathbf{u}_J)$ than $\spn(\mathbf{u}_{\widehat{J}})$ still satisfies the same upper bound. That is,
if $\|\widetilde{M} - M\|_2 < \frac{1}{2} \mathrm{sep}(\lambda_J, \lambda_{J^c})$, then
%
%
\begin{align}
\min_{\widetilde{J} \in S_{q,n}} d_\mathrm{\subspacetext} \left( \spn(\mathbf{u}_J) , \spn(\widetilde{\mathbf{u}}_{\widetilde{J}}) \right)
& ~~\leq~~
\min\left( 1, \sqrt{\frac{n-q}{q}} \right) ~ \frac{ 
	\left\| \widetilde{M} - M \right\|_2 
}{ 
	\mathrm{sep}(\lambda_J, \lambda_{J^c})  ~-~ \|\widetilde{M} - M\|_2  
} \nonumber \\
\min_{\widetilde{J} \in S_{q,n}} d_\mathrm{\subspacetext} \left( \spn(\mathbf{u}_J) , \spn(\widetilde{\mathbf{u}}_{\widetilde{J}}) \right)
& ~~\leq~~
 \frac{\frac{1}{\sqrt{2q}} \left\| \widetilde{M} - M \right\|_F}{\mathrm{sep}(\lambda_J, \lambda_{J^c})  ~-~ \|\widetilde{M} - M\|_2} 
\end{align}
where $S_{q,n}$ is the set of all $q$-element subsets of $N = \{1,2,\cdots,n\}$.

\section{Application to Null-space Perturbation in Context of a Graph Connection Problem} \label{sec:graph}

We consider a simple application of the above results 
in context of a graph theory problem.
Some definitions and basic properties of a weighted, undirected, simple graphs are listed below~\cite{godsil2001algebraic}:
\begin{enumerate}
	\item A graph, $G$, constitutes of a set of $n$ vertices, $\mathcal{V}(G) = \{v_1, v_2, \cdots, v_n\}$ and an edge set $\mathcal{E}(G) \subseteq \mathcal{V}(G) \times_{\text{sym}} \mathcal{V}(G)$ (where `$\times_{\text{sym}}$' represent the symmetric Cartesian product so that for the undirected graph the order of the vertices in a edge is irrelevant, making $(v_k, v_l) = (v_l, v_k)$).
	Each edge, $(v_k, v_l)\in \mathcal{E}(G)$, is assigned a positive real weight, $A_{kl} (= A_{lk})$. 
	Non-existent edges are implicitly assumed to have zero edge weight so that $A_{kl} = 0, \, \forall (v_k, v_l)\notin \mathcal{E}(G)$. The matrix $A\in \mathbb{R}^{n\times n}$ is called the \emph{weighted adjacency matrix} of the graph $G$, and is a symmetric matrix with zero diagonal for an undirected, simple graph.
	\item The \emph{weighted degree matrix}, $D$, is a $n\times n$ diagonal matrix in which the $k^\text{th}$ diagonal element is the sum of the elements in the $k^\text{th}$ row (equivalently, $k^\text{th}$ column) of $A$. Thus $D_{kk}$ is the sum of the weights of the edges emanating from $v_k$ (also called the \emph{degree} of the vertex).
	\item The \emph{weighted Laplacian matrix} of the graph is defined as $L = D - A$. An eigenvector of $L$ is a $n$-dimensional real vector and can be interpreted as a distribution over the vertices (with the $k^\text{th}$ element of the vector being the value associated to $v_k\in \mathcal{E}(G)$).
	\item The eigenvalues of $L$ are non-negative. The null-space of $L$ for a graph with $q$ disjoint components is $q$-dimensional, with the null-space spanned by vectors corresponding to distributions that are uniform over the vertices of each of those components. Without loss of generality we index the eigenvalues in increasing order of their magnitudes so that $0 = \lambda_1 = \lambda_2 = \cdots = \lambda_q \leq \lambda_{q+1} \leq \lambda_{q+2} \leq \cdots \leq \lambda_n$. The corresponding unit eigenvectors be $\mathbf{u}_1, \mathbf{u}_2, \cdots, \mathbf{u}_n$.
	Note that since a graph has at least one connected component, $\lambda_1 = 0$ for any graph.
	Furthermore, without loss of generality, we choose 
	$\mathbf{u}_j$ to be a distribution that is uniformly positive over the vertices if $G_j$, and zero over the rest of the vertices in the graph.
	\item Define $J = \{1,2,\cdots,q\}$, so that $\spn(\mathbf{u}_J)$ is the null-space of $L$.
\end{enumerate}

If $G$ has $q$ disjoint components, we define $G_j, j=1,2,\cdots,q$ to be the subgraph constituting of the vertices and edges in the $j^\text{th}$ component only.
Thus, $\mathcal{V}(G) = \cup_{j=1}^q \mathcal{V}(G_j)$ and $\mathcal{E}(G) = \cup_{j=1}^q \mathcal{E}(G_j)$ (more compactly, we write $G = \cup_{j=1}^q G_j$).
We also define the collection of these subgraphs as \[ \mathsf{G}=\{G_1,G_2,\cdots, G_q\} \]

We are interested in understanding perturbation of the invariant subspace, $\spn(\mathbf{u}_J)$ (the null-space), of $L$ as new edges are established between the different disjoint components (henceforth also referred to as ``\emph{clusters}'') of the graph.
Let the graph constructed by establishing the inter-cluster edges be $\widetilde{G}$ with $\widetilde{A}$, $\widetilde{D}$ and $\widetilde{L}$ its adjacency, degree and Laplacian matrices respectively.
Note that since $\widetilde{G}$ is constructed by just adding edges between the subgraphs $\{G_j\}_{j=1,2,\cdots,q}$ of $G$, each of these subgraphs are induced subgraphs of $\widetilde{G}$.

\subsection{Computation of $\|(\widetilde{L} - L) \mathbf{u}_j\|_2$} 

For any induced subgraph, $H \in \widetilde{G}$, we consider the edges that connect vertices in $H$ to vertices not in $H$ (\emph{inter-cluster edges}). 
These are edges of the form $(v_k,v_l)$ such that $v_k\in \mathcal{V}(H), v_l\notin \mathcal{V}(H)$.
We define a few quantities involving the weights on such edges.

\begin{definition}~ \label{def:ext-deg-graph}

\begin{itemize}
	\item[1.] \emph{External Degree of a Vertex Relative to a Subgraph:} Given an subgraph, $H\subseteq \widetilde{G}$, and a vertex $v_k\in\mathcal{V}(H)$, the \emph{external degree of $v_k$ relative to $H$ in $\widetilde{G}$} is defined as the sum of the weights on edges connecting $v_k$ to vertices outside $H$:
	{\small \begin{align}
		\externaldegtext_{H,\widetilde{G}}(v_k) ~~=~~
				\sum_{\{l\,|\, v_l\notin\mathcal{V}(H)\}} \widetilde{A}_{kl}
		\end{align}}
	\item[2.] \emph{Coupling of a subgraph in a Graph:} Given a induced subgraph, $H\subseteq \widetilde{G}$, we define the \emph{coupling of $H$ in $\widetilde{G}$} as
	{\small \begin{align}
		\couplingtext_{\widetilde{G}}(H) 
			&~~=~~ { \frac{1}{|\mathcal{V}(H)|} \left( 
				\sum_{\{k\,|\,v_k\in\mathcal{V}(H)\}} \left( \externaldegtext_{H,\widetilde{G}}(v_k) \right)^2 ~+~ 
				\sum_{\{l\,|\,v_l\notin\mathcal{V}(H)\}} \left( \externaldegtext_{(\widetilde{G}-H),\widetilde{G}}(v_l) \right)^2
				\right) } \nonumber \\
			&~~=~~ { \frac{1}{|\mathcal{V}(H)|} \left( 
			\sum_{\{k\,|\,v_k\in\mathcal{V}(H)\}} \left( \sum_{\{l\,|\,v_l\notin\mathcal{V}(H)\}} \widetilde{A}_{kl} \right)^2 ~+~ 
			\sum_{\{l\,|\,v_l\notin\mathcal{V}(H)\}} \left( \sum_{\{k\,|\,v_k\in\mathcal{V}(H)\}} \widetilde{A}_{kl} \right)^2
			\right) }
		\end{align}}
	where $(\widetilde{G}-H)$ is the induced subgraph of $\widetilde{G}$ constituting of all the vertices not in $H$. That is, $\mathcal{V}(\widetilde{G}-H) = \{v\in\mathcal{V}(\widetilde{G})\,|\,v\notin\mathcal{V}(H)\}$ and $\mathcal{E}(\widetilde{G}-H) = \{(v,w)\in\mathcal{E}(\widetilde{G})\,|\,v,w\notin \mathcal{V}(H)\}$. 
	\item[3.] \emph{Maximum External Degree of Vertices in a Subgraph:} Given an subgraph, $H\subseteq \widetilde{G}$, the \emph{maximum external degree of vertices in $H$ in $\widetilde{G}$} is defined as the maximum value of the external degrees of vertices in $H$ relative to $H$ in $\widetilde{G}$:
	{\small \begin{align}
		\maxextdeg_{\widetilde{G}}(H) ~~=~~ \max_{v\in\mathcal{V}(H)} \externaldegtext_{H,\widetilde{G}}(v) ~~=~~
		\max_{\{k\,|\,v_k\in\mathcal{V}(H)\}}  \sum_{\{l\,|\, v_l\notin\mathcal{V}(H)\}} \widetilde{A}_{kl}
		\end{align}}
\end{itemize}
	Note that the computation of the above quantities require the knowledge of only the weights on edges connecting vertices in $H$ to vertices outside $H$ in $\widetilde{G}$.
\end{definition}

\begin{wrapfigure}{r}{0.55\columnwidth}
	\hspace{2em} {\includegraphics[width=0.4\columnwidth, trim=15 20 20 15, clip=true]{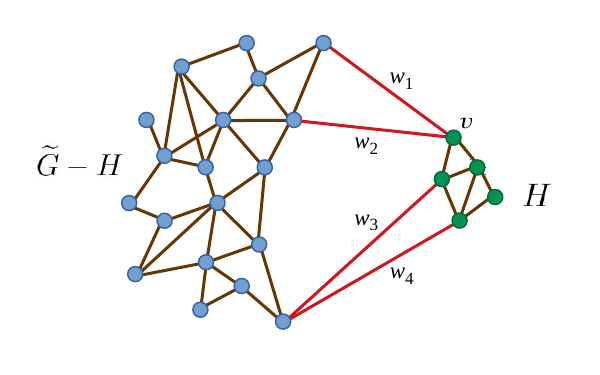}} \hspace{2em} 
	\caption{An example graph, $\widetilde{G}$, and induced subgraph, $H$. Weight values on the inter-cluster edges are written symbolically. In this example, 
	~$\externaldegtext_{H,\widetilde{G}}(v) = w_1 + w_2$, 
	~$\couplingtext_{\widetilde{G}}(H) = \frac{1}{5}\left( \left( (w_1 + w_2)^2 + w_3^2 + w_4^2 \right) + \left( w_1^2 + w_2^2 + (w_3 + w_4)^2 \right) \right)$, and,
	~$\maxextdeg_{\widetilde{G}}(H) = \max\left(w_1 \!+\! w_2,\, w_3,\, w_4 \right)$
} \label{fig:graph-ext-deg} 
\end{wrapfigure}

In the definition of $\couplingtext_{\widetilde{G}}$, referring to $H$ as a \emph{cluster} and considering the rest of the graph another cluster, the quantity within the innermost brackets is the sum of the weights on inter-cluster edges connected to a vertex, which is squared and summed over all the vertices that have at least one inter-cluster edge connected to it.
This quantity is then divided by the number of vertices in $H$.
Thus a large subgraph which is weakly connected to the rest of the graph will have a lower coupling value.

The following lemma provides bounds on $\couplingtext_{\widetilde{G}}(H)$ in terms of a simpler summation over the inter-cluster edge weights (or square thereof).

\begin{lemma} \label{lemma:coupling-inequalities}
\begin{align}
\frac{2}{|\mathcal{V}(H)|}  \sum_{\{k,l \,|\, v_k\in \mathcal{V}(H), \atop\qquad v_l\notin \mathcal{V}(H)\}} \!\!\!\! \widetilde{A}_{kl}^2 
\qquad \leq \qquad \couplingtext_{\widetilde{G}}(H) \qquad \leq \qquad
\frac{2}{|\mathcal{V}(H)|} \left( \sum_{\{k,l \,|\, v_k\in \mathcal{V}(H), \atop\qquad v_l\notin \mathcal{V}(H)\}} \!\!\!\! \widetilde{A}_{kl} \right)^2
\end{align}
\end{lemma}
\begin{proof}
	The proof follows directly using the fact that for a set of positive numbers, $\alpha_h,\, h\in S$, $\sum_{h\in S} \alpha_h^2 \leq (\sum_{h \in S} \alpha_h)^2$.
\end{proof}

\noindent
\textbf{Notations and Assumptions for the Rest of the Paper:}
In the rest of the paper we assume that $G$ is a graph with $q$ disjoint components, $\mathsf{G} = \{G_1,G_2,\cdots,G_q\}$, and $\widetilde{G}$ be the graph obtained by establishing edges between the components (so that each $G_j$ is an induced subgraph of both $G$ and $\widetilde{G}$). The Laplacian matrices of the two graphs be $L$ and $\widetilde{L}$ respectively.
Since $G$ has $q$ connected components, its null-space is $q$ dimensional (with corresponding eigenvalues $\lambda_1=\lambda_2=\cdots=\lambda_q$), for which we choose a basis $\{\mathbf{u}_j\}_{j=1,2\cdots,q}$ such that the distribution corresponding to $\mathbf{u}_j$ is uniform and positive on the vertices in $G_j$, and zero everywhere else.

\vspace{0.5em}
A weaker version of the following lemma appears in the author's prior work~\cite{zhang2020inter,zhan2021barrier-journal} and expresses the quantity $\|(\widetilde{L} - L) \mathbf{u}_j\|_2$ in terms of the weights on edges connecting vertices in $G_j$ to vertices outside $G_j$ in $\widetilde{G}$.

\begin{lemma} \label{lemma:diffnorm-graph-inter-cluster-weight}
For all $j\in\{1,2,\cdots,q\}$,
{\small	\begin{align}
	\|(\widetilde{L} - L) \mathbf{u}_j\|_2^2 
	&~~=~~ 
	\couplingtext_{\widetilde{G}}(G_j)
	\end{align}}
\end{lemma}

\begin{proof}
Suppose $v_k\in \mathcal{V}(G_j) \subseteq \mathcal{V}(G)$.
Since $\widetilde{D}_{kk}$ and $D_{kk}$ are the degrees of the vertex 
in the graphs $\widetilde{G}$ and $G$ respectively,
they are equal iff all the neighbors of $v_k$ are in $G_j$. Otherwise $\widetilde{D}_{kk} - D_{kk}$ is the net outgoing degree of the vertex $v_k$ from the subgraph $G_j$. That is, if $v_k\in \mathcal{V}(G_j)$,
\begin{align} \label{eq:D-diff}
\widetilde{D}_{kk} - D_{kk} & ~~=~~ \sum_{\{l\,|\, v_l\notin \mathcal{V}(G_j)\}} \widetilde{A}_{kl}
\end{align}

An edge $(v_k,v_l)$ exists in both $\widetilde{G}$ and $G$ (and have the same weight, \emph{i.e.}, $\widetilde{A}_{kl} = A_{kl}$)
iff $v_k$ and $v_l$ belong to the same subgraph $G_j$. Otherwise $A_{kl}=0$ (the edge is non-existent in $G$). Thus,
\begin{equation} \label{eq:A-diff}
\widetilde{A}_{kl} - A_{kl} ~=~ \left\{ \begin{array}{l} \widetilde{A}_{kl}, ~\text{if $v_k \in \mathcal{V}(G_j), v_l \notin \mathcal{V}(G_j)$} \\ 0, ~~\text{otherwise.} \end{array} \right.
\end{equation}

Next we consider the vector $\mathbf{u}_j$ (for $j=1,\cdots,q$), which by definition is non-zero and uniform only on vertices in the subgraph $G_j$.
Let $u_{lj}$ be the $l$-th element of the unit vector $\mathbf{u}_j$. Since $|\mathcal{V}(G_j)|$ of the elements of the vector are non-zero and uniform, we have,
\begin{equation} \label{eq:def-ulj}
u_{lj} = \left\{ \begin{array}{l} \frac{1}{\sqrt{|\mathcal{V}(G_j)|}}, ~~\text{if $v_l \in \mathcal{V}(G_j)$} \\ 0, ~~\text{otherwise} \end{array} \right.
\end{equation}

\noindent
Thus the $k$-th element of the vector $(\widetilde{L}-L) \mathbf{u}_j$,
{\small \begin{align} 
	 [(\widetilde{L}-{L}) {\mathbf{u}}_j]_k ~~ 
	= & \displaystyle \sum_{l} (\widetilde{D}_{kl} - \widetilde{A}_{kl} - {D}_{kl} + {A}_{kl}) {u}_{lj} \nonumber \\
	= & \displaystyle (\widetilde{D}_{kk} - {D}_{kk}) {u}_{kj} - \sum_{l} (\widetilde{A}_{kl} - {A}_{kl}) {u}_{lj} 
	\qquad\text{\small(since $\widetilde{D}$ and ${D}$ are diagonal marices)} \nonumber \\
	= & \displaystyle \left( \frac{1}{\sqrt{|\mathcal{V}(G_j)|}} \left\{ \begin{array}{l}  (\widetilde{D}_{kk} - {D}_{kk}), ~\text{if $v_k \in \mathcal{V}(G_j)$} \nonumber \\ 0, ~\text{otherwise} \end{array} \right.\right) \nonumber \\
	& \displaystyle \quad\qquad
	\qquad - \left( \frac{1}{\sqrt{|\mathcal{V}(G_j)|}} \sum_{\{l\,|\,v_l \in \mathcal{V}(G_j)\}} (\widetilde{A}_{kl} - {A}_{kl}) \right) 
	\qquad\text{\small(using \eqref{eq:def-ulj})} \nonumber \\
	= & \displaystyle \frac{1}{\sqrt{|\mathcal{V}(G_j)|}} \left( \left\{ \begin{array}{l}  \displaystyle \sum_{\{l\,|\, v_l\notin \mathcal{V}(G_j)\}} \!\! \widetilde{A}_{kl}, ~~~\text{if $v_k \in \mathcal{V}(G_j)$} \nonumber \\ 0, \qquad\qquad\text{otherwise} \end{array} \right.\right. \nonumber \\
	& \displaystyle \qquad\qquad - \left. \sum_{\{l\,|\,v_l \in \mathcal{V}(G_j)\}} \left\{ \begin{array}{l} \widetilde{A}_{kl},~~\text{if $v_k\notin \mathcal{V}(G_j)$} \nonumber \\ 0, ~\text{otherwise.} \end{array} \right. \right) 
	\quad\qquad\text{\small(using \eqref{eq:D-diff} and \eqref{eq:A-diff})} \nonumber \\
	= & \displaystyle \frac{1}{\sqrt{|\mathcal{V}(G_j)|}} 
	\left\{ \begin{array}{l}  
	\displaystyle \sum_{\{l\,|\, v_l\notin \mathcal{V}(G_j)\}} \!\!\!\! \widetilde{A}_{kl}, ~~~\text{if $v_k \in \mathcal{V}(G_j)$} \\
	\displaystyle -\sum_{\{l\,|\,v_l \in \mathcal{V}(G_j)\}} \widetilde{A}_{kl}, ~~~~~\text{if $v_k\notin \mathcal{V}(G_j)$} 
	\end{array} \right.
	%
	\end{align}}

\noindent
Thus, 
{\small \begin{align} 
	\|(\widetilde{L}-{L}) {\mathbf{u}}_j\|_2^2  
	& ~~=~~ \displaystyle \frac{1}{|\mathcal{V}(G_j)|}~~
	\left(
	\displaystyle \sum_{\{k\,|\,v_k \in \mathcal{V}(G_j)\}} \left( \sum_{\{l\,|\, v_l\notin \mathcal{V}(G_j)\}} \!\!\!\! \widetilde{A}_{kl} \right)^2 
~+~
	\displaystyle \sum_{\{k\,|\,v_k\notin \mathcal{V}(G_j)\}} \left( \sum_{\{l\,|\,v_l \in \mathcal{V}(G_j)\}} \widetilde{A}_{kl} \right)^2
	\right) 
	\end{align}}


\end{proof}

\begin{lemma} \label{lemma:laplacian-diff-2-norm-upperbound}
\begin{align}
\|\widetilde{L} - L\|_2 ~~\leq~~ 2 \, \max_{j\in\{1,\cdots,q\}} 
\maxextdeg_{\widetilde{G}}(G_j) 
\end{align}
\end{lemma}

\begin{proof}
Suppose $v_k \in \mathcal{V}(G_{\mathsf{j}(k)})$ (where $\mathsf{j}: \{1,2,\cdots,|\mathcal{V}(G)|\} \rightarrow \{1,2,\cdots,q\}$ maps the index of a vertex to the index of the subgraph in $\mathsf{G}$ that the vertex belongs to).
The sum of the elements of the $k^\text{th}$ row of $(\widetilde{A} - A)$ is
\begin{align}
\sum_{l} (\widetilde{A}_{kl} - {A}_{kl}) ~~=~ & \sum_{l} \left\{ \begin{array}{l} \widetilde{A}_{kl}, ~\text{if $v_l \notin \mathcal{V}(G_{\mathsf{j}(k)})$} \\ 0, ~~\text{otherwise.} \end{array} \right. \qquad\text{\small(Using \eqref{eq:A-diff}.)} \nonumber \\
=~ & \displaystyle \sum_{\{l\,|\, v_l\notin \mathcal{V}(G_{\mathsf{j}(k)})\}} \widetilde{A}_{kl} \nonumber \\
=~ & \externaldegtext_{G_{\mathsf{j}(k)},\widetilde{G}} (v_k) \qquad\text{\small(Definition~\ref{def:ext-deg-graph})} \label{eq:a-diff-row-sum}
\end{align}
Since, $(\widetilde{A} - A)$ is symmetric matrix, its $2$-norm is equal to its spectral radius, $\rho(\widetilde{A} - A)$. Furthermore, since all elements of  $(\widetilde{A} - A)$ are non-negative, using Perron-Frobenius theorem~\cite{doi:10.1137/1.9781611971262}, we get
\begin{align}
\| \widetilde{A} - A \|_2 ~~=~~ \rho(\widetilde{A} - A) ~~\leq & ~ \max_{k\in N} \,\externaldegtext_{G_{\mathsf{j}(k)},\widetilde{G}} (v_k) \nonumber \\
& =~ \max_{j\in\{1,\cdots,q\}} \max_{\{k\,|\atop v_k\in\mathcal{V}(G_j)\}} \,\externaldegtext_{G_{\mathsf{j}(k)},\widetilde{G}} (v_k) \nonumber \\ & \qquad\qquad\qquad{\small\begin{array}{l}\text{(since maximizing over all vertices in $\widetilde{G}$ is same as maximizing} \\ \text{ over the subgraphs, $G_j$, and for each subgraph maximizing over} \\ ~\text{ the vertices in the subgraph.)}\end{array}} \nonumber \\
& =~ \max_{j\in\{1,\cdots,q\}} \maxextdeg_{\widetilde{G}}(G_j)  \label{eq:a-diff-2-norm}
\end{align}

Again, since $(\widetilde{D} - D)$ is a diagonal matrix with positive diagonal elements (due to \eqref{eq:D-diff}), its $2$-norm is the maximum out of its diagonal elements. That is,
\begin{align}
\|\widetilde{D} - D\|_2 ~~=~~ \max_{k\in N} \, (\widetilde{D}_{kk} - {D}_{kk}) 
~~=~ & \displaystyle \max_{k\in N} \,  \sum_{\{l\,|\, v_l\notin \mathcal{V}(G_{\mathsf{j}(k)})\}} \widetilde{A}_{kl} \qquad\text{\small(Using \eqref{eq:D-diff}.)} \nonumber \\
& =~ \max_{k\in N} \externaldegtext_{G_{\mathsf{j}(k)},\widetilde{G}} (v_k) 
~~=~~ \max_{j\in\{1,\cdots,q\}} \maxextdeg_{\widetilde{G}}(G_j) \nonumber \\ & \qquad\qquad\qquad\text{\small(following similar steps as in \eqref{eq:a-diff-row-sum} and \eqref{eq:a-diff-2-norm}.)}
\end{align}

Thus,
\begin{align*}
\|\widetilde{L} - L\|_2 ~~=~~ \|(\widetilde{D} - D) - (\widetilde{A} - A)\|_2 ~~\leq~~ \|\widetilde{D} - D\|_2 + \|\widetilde{A} - A\|_2 ~~\leq~~ 2 \max_{j\in\{1,\cdots,q\}} \maxextdeg_{\widetilde{G}}(G_j)
\end{align*}

\end{proof}

In the following discussions,
without loss of generality, we assume that the eigenvalues of $\widetilde{L}$ are indexed in increasing order of magnitude, $0 = \widetilde{\lambda}_1 \leq \widetilde{\lambda}_2 \leq \cdots \leq \widetilde{\lambda}_n$.
The corresponding eigenvectors be $\widetilde{\mathbf{u}}_1,\widetilde{\mathbf{u}}_2,\cdots,\widetilde{\mathbf{u}}_n$.

%
%

\subsection{Bounds on Null-space Purturbation with Known Spectrum of $\widetilde{L}$}

%
The following proposition 
gives a bound on the perturbation of the null-space of $L$ upon introducing edges between the subgraphs in $\mathsf{G} = \{G_1,G_2,\cdots,G_q\}$ by considering the sub-space distance between the null-space of $L$ and a specific invariant sub-space of $\widetilde{L}$.

\begin{prop} 
\label{prop:general-coupling-result}
%
Choose $\widetilde{J} = \{1,2,\cdots,q\}$. Then
	\begin{align}
	d_\mathrm{\subspacetext} \left( \spn(\mathbf{u}_J) , \spn(\widetilde{\mathbf{u}}_{\widetilde{J}}) \right)
	& ~~~~\leq~~~~
		\displaystyle \frac{1}{\widetilde{\lambda}_{q+1}} ~\sqrt{ \frac{1}{q} \sum_{j=1}^q \couplingtext_{\widetilde{G}}(G_j) } 
\label{eq:graph-subspace-inequality-kappa-zero}
	\end{align}
\end{prop}

\begin{proof}
We first note that due to Lemma~\ref{lemma:diffnorm-graph-inter-cluster-weight}
$\sqrt{\couplingtext_{\widetilde{G}}(G_j)} =  \|(\widetilde{L} - L) \mathbf{u}_{{j}}\|_2, 
\,\forall j\in \{1,2,\cdots,q\}$.
The proof then follows from Proposition~\ref{prop:full-main} by setting $\kappa_j = 0,\,\forall j=1,2,\cdots q$ and noting that $\displaystyle \min_{j\in \{1,2,\cdots,q\} \atop j'\in \{q+1,q+2,\cdots,n\}} |\widetilde{\lambda}_{j'} - \lambda_j| = \widetilde{\lambda}_{q+1}$.
\end{proof}

The results of Proposition~\ref{prop:general-coupling-result} can be re-interpreted by considering $G$ to be the graph obtained by \emph{cutting} $\widetilde{G}$ into $q$-subgraphs. 
We call the set of subgraphs hence constructed upon performing the cut, $\mathsf{G} = \{G_1,G_2,\cdots, G_q\}$, a $q$-cut of $\widetilde{G}$.
%
Given a graph $\widetilde{G}$,
we consider all possible $q$-cuts of $\widetilde{G}$.
A $q$-cut, $\mathsf{G} = \{G_1,G_2,\cdots, G_q\}$, results in a graph, $G = \cup_{j=1}^q G_j$, with $q$ disjoint components.
%
The following corollary is then a direct consequence of the proposition. 

\begin{cor}
Given a graph $\widetilde{G}$ (with Laplacian $\widetilde{L}$ with eigenvalues $0 = \widetilde{\lambda}_1 \leq \cdots \leq \widetilde{\lambda}_n$ and corresponding eigenvectors $\widetilde{\mathbf{u}}_1,\cdots,\widetilde{\mathbf{u}}_n$), let $\mathscr{G}$ be the set of all $q$-cuts of $\widetilde{G}$.
We consider a $q$-cut such that the sum of the couplings of the resultant $q$ subgraphs in $\widetilde{G}$ is minimum. That is,
\begin{equation} \label{eqn:best-q-cut}
\mathsf{G}^{*} ~\in~ {\arg\!\min}_{\mathsf{G} \in \mathscr{G}} \sum_{G' \in \mathsf{G}} \couplingtext_{\widetilde{G}}(G')
\end{equation}
Let the corresponding graph, $\displaystyle G^{*} = \bigcup_{G'\in\mathsf{G}^{*}} G'$, have eigenvalues $0 = \lambda_1^{*} = \lambda_2^{*} = \cdots = \lambda_q^{*} \leq \lambda_{q+1}^{*} \leq \cdots \leq \lambda_{n}^{*}$ and corresponding eigenvectors $\mathbf{u}^{*}_1, \mathbf{u}^{*}_2, \cdots, \mathbf{u}^{*}_n$. Then,
	\begin{align}
	d_\mathrm{\subspacetext} \left( \spn(\{\widetilde{\mathbf{u}}_1,\cdots,\widetilde{\mathbf{u}}_q\}) , \spn(\{\mathbf{u}^{*}_1,\cdots,\mathbf{u}^{*}_q\}) \right)
	& ~~~~\leq~~~~
	\displaystyle \frac{1}{\widetilde{\lambda}_{q+1}} ~\sqrt{ \frac{1}{q} \sum_{G' \in \mathsf{G}^{*}} \couplingtext_{\widetilde{G}}(G') } 
	\end{align}
\end{cor}
The interpretation of the above corollary is that the ``best'' $q$-cut of a graph $\widetilde{G}$ (minimizing total inter-cluster coupling, as defined by \eqref{eqn:best-q-cut}) results in a graph such that the distance between the nullspace of the cut graph's Laplacian and 
the space spanned by the first $q$ eigenvectors 
of the Laplacian of $\widetilde{G}$
is bounded above by a quantity proportional to the total inter-cluster coupling (which was minimized in the first place).

\subsection{Bounds on Null-space Purturbation with Known Spectrum of $L$}

\begin{prop} \label{prop:null-space-purturbation-known-L-spectrum}
If ~$\max_{j\in\{1,\cdots,q\}} \maxextdeg_{\widetilde{G}}(G_j) ~<~ \frac{\lambda_{q+1}}{4}$, then
{\begin{align}
d_\mathrm{\subspacetext} \left( \spn(\mathbf{u}_J) , \spn(\widetilde{\mathbf{u}}_{\widehat{J}}) \right)
& ~\leq~
\frac{ 
	\sqrt{\displaystyle \frac{1}{q}~ \sum_{j=1}^q 
	\couplingtext_{\widetilde{G}}(G_j) }
}{\displaystyle 
	 \lambda_{q+1} ~-~ 
	2 \, \max_{k\in\{1,\cdots,q\}} \maxextdeg_{\widetilde{G}}(G_k)  
} 
\label{eq:null-space-purturbation-known-L-spectrum-1} 
~\leq~
\frac{ 
	2 \, \displaystyle\max_{j\in\{1,\cdots,q\}} \maxextdeg_{\widetilde{G}}(G_j)  
}{ 
	\lambda_{q+1} ~-~ 2 \, \displaystyle \max_{j\in\{1,\cdots,q\}} \maxextdeg_{\widetilde{G}}(G_j)   
} 
\end{align}}
where $\widehat{J} = \{1,2,\cdots,q\} = \{j' ~~|~~ \min_{j\in N} |\widetilde{\lambda}_{j'} - \lambda_j| = \widetilde{\lambda}_{j'} \} $.
\end{prop}
\begin{proof}
Recall, that the eigenvalues of the Laplacian, $L$, of $G$, are $(0=)\lambda_1=\lambda_2=\cdots=\lambda_q \leq \lambda_{q+1} \leq \cdots \leq \lambda_n$.
Let $J = \{1,2,\cdots,q\}$ so that ${J^c} = \{q+1,q+2,\cdots,n\}$ and $\mathrm{sep}(\lambda_J,\lambda_{J^c}) = \lambda_{q+1}$.
	
Using Lemma~\ref{lemma:laplacian-diff-2-norm-upperbound},
\begin{align} \label{eq:L-diff-sep-inequality}
\|\widetilde{L} - L\|_2 ~~\leq~~ 2 \, \max_{j\in\{1,\cdots,q\}} \maxextdeg_{\widetilde{G}}(G_j) 
 ~~<~~ \frac{\lambda_{q+1}}{2} ~~=~~ \frac{\mathrm{sep}(\lambda_J,\lambda_{J^c})}{2}
\end{align}
Thus the condition for Lemma~\ref{lemma:sep-preserving-M-tilde} and Proposition~\ref{prop:subspace-purturbation-hat-tilde-free} hold, and $\widetilde{L}$ is a separation preserving perturbation of $L$.
Hence, by Lemma~\ref{lemma:sep-preserving-M-tilde} there exists a separation preserving partition, $\{\widetilde{\lambda}_{\widehat{J}}, \widetilde{\lambda}_{\widehat{{J^c}}}\}$ of $\widetilde{\lambda}_{N}$ such that

\[
\widehat{J} = \{j' ~~|~~ \min_{j\in N} |\widetilde{\lambda}_{j'} - \lambda_j| = \min_{j\in J} |\widetilde{\lambda}_{j'} - \lambda_j| = \widetilde{\lambda}_{j'} \} \qquad\text{\small (since $\lambda_j = 0, \,\forall \,j\in J$)}
\]
Thus, for any $j'\in \widehat{J}$,
\begin{align}
 \widetilde{\lambda}_{j'} ~~=~~ \min_{j'\in N} |\widetilde{\lambda}_{j'} - \lambda_j| & ~\leq~~ \|\widetilde{L} - L\|_2 \quad\text{\small (due to Corollary~\ref{cor:M-diff-lambda-diff})} \nonumber \\
& ~\leq~~ \frac{\lambda_{q+1}}{2} \quad\text{\small (from \eqref{eq:L-diff-sep-inequality})} 
\end{align}
This implies that the elements of $\widetilde{\lambda}_{\widehat{J}}$ are closer to $0 (=\lambda_1=\lambda_2=\cdots=\lambda_q)$ than they are to $\lambda_{q+1}$.
Since $\widehat{J}$ has $q$-elements (due to Lemma~\ref{lemma:sep-preserving-M-tilde}.2) 
and is a unique set (by definition), we have $\widetilde{\lambda}_{\widehat{J}} = \{\widetilde{\lambda}_1, \widetilde{\lambda}_2, \cdots, \widetilde{\lambda}_q\}$ to be the set constituting of the lowest $q$ eigenvalues of $\widetilde{L}$.
Thus, $\widehat{J} = \{1,2,\cdots,q\}$.

Since we showed that $\|\widetilde{L} - L\|_2 ~\leq~ \frac{1}{2} \mathrm{sep}(\lambda_J,\lambda_{J^c})$, as direct consequence of Proposition~\ref{prop:subspace-purturbation-hat-tilde-free} we have the following

\vspace{-2em}
\begin{align}
\left( d_\mathrm{\subspacetext} \left( \spn(\mathbf{u}_J) , \spn(\widetilde{\mathbf{u}}_{\widehat{J}}) \right) \right)^2
& ~~\leq~~
\frac{1}{q}\,
\mathlarger{\mathlarger{\sum}}_{j\in J} \left(\frac{ 
	\left\| (\widetilde{L} - L) \mathbf{u}_{{j}} \right\|_2 
}{ 
	\min_{k\in {J^c}} |{\lambda}_{k} - \lambda_j| ~-~ \|\widetilde{L} - L\|_2  
} \right)^{\!\!2} \nonumber \\
& \quad ~~\leq~~
\frac{ 
\displaystyle \frac{1}{q}~ \sum_{j\in J} 
\couplingtext_{\widetilde{G}}(G_j) 
}{\displaystyle 
	\left( \lambda_{q+1} ~-~ 
	2 \, \max_{k\in\{1,\cdots,q\}} \maxextdeg_{\widetilde{G}}(G_k)   \right)^{\!\!2} 
} 
\nonumber \\ & \qquad\qquad\qquad\qquad\text{\small (using Lemma~\ref{lemma:diffnorm-graph-inter-cluster-weight} and Lemma~\ref{lemma:laplacian-diff-2-norm-upperbound} and the} \nonumber \\ & \quad\qquad\qquad\qquad\qquad\text{\small fact that $\min_{k\in {J^c}} |{\lambda}_{k} - \lambda_j| = \lambda_{q+1},\,\forall j\in \{1,2,\cdots,n\}$.)}\nonumber \\
& \quad ~~\leq~~
\left(\frac{ 
	2 \, \displaystyle\max_{j\in\{1,\cdots,q\}} \maxextdeg_{\widetilde{G}}(G_j)  
}{ 
	\lambda_{q+1} ~-~ 2 \, \displaystyle \max_{j\in\{1,\cdots,q\}} \maxextdeg_{\widetilde{G}}(G_j)   
} \right)^{\!\!2} 
\nonumber \\ & \qquad\qquad\qquad\text{\small (using Lemma~\ref{lemma:diffnorm-graph-inter-cluster-weight} and \ref{lemma:laplacian-diff-2-norm-upperbound}, 
	$\sum_{j\in J} \couplingtext_{\widetilde{G}}(G_j) = \sum_{j\in J} \left\| (\widetilde{L} - L) \mathbf{u}_{{j}} \right\|_2^2 $} \nonumber \\ & \qquad\qquad\qquad\qquad\qquad\qquad \text{\small $\leq q \|\widetilde{L} - L\|_2^2 \leq q (2 \, \displaystyle\max_{j\in\{1,\cdots,q\}} \maxextdeg_{\widetilde{G}}(G_j) )^2$.)} \nonumber
 \end{align}

\vspace{-1em}
\end{proof}

\vspace{1em}
\subsection{Example}

\begin{wrapfigure}{r}{0.58\columnwidth} \vspace{-3.5em}
\centering
	\subfloat[Graph $G$ with $12$ disjoint components.]{\includegraphics[width=0.26\textwidth, clip=true, trim=95 210 75 225]
		{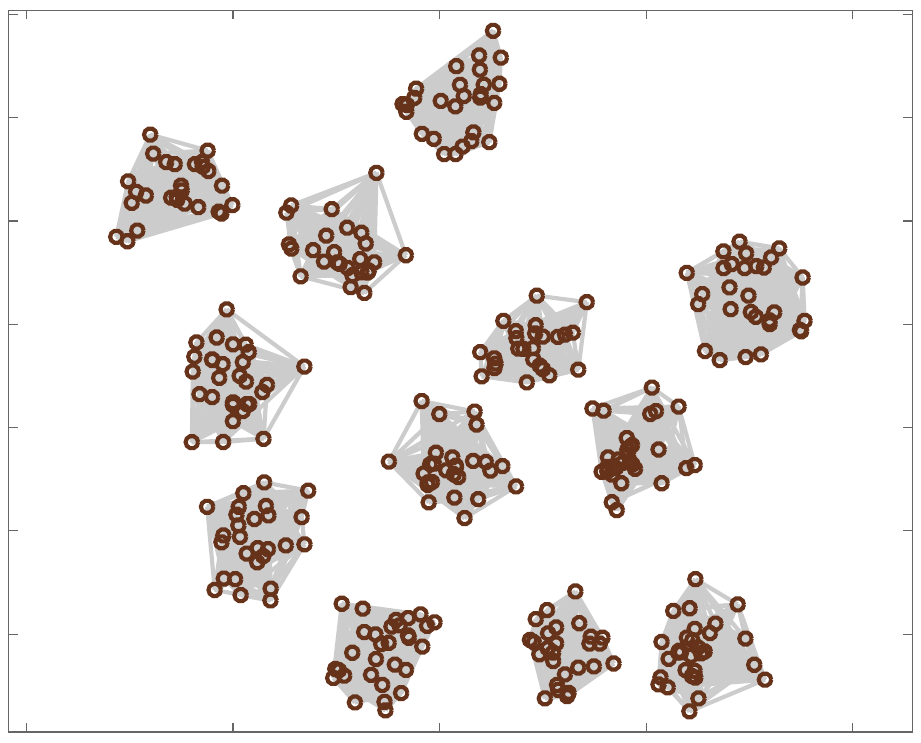} \label{fig:disjoint-graph-graph}} \hspace{0.2em}
	\subfloat[The spectrum of $L$ with first $12$ eigenvalues equal to zero.]{\includegraphics[width=0.28\textwidth, clip=true, trim=75 200 70 225]
		{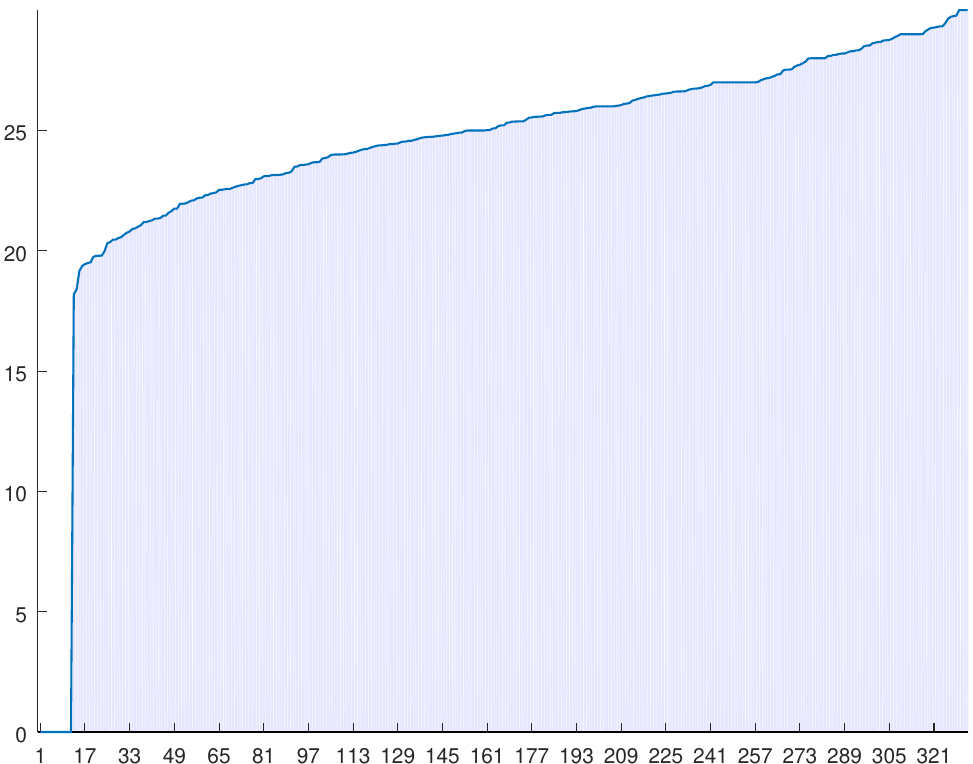} \label{}}
	 \vspace{-0.5em} \caption{Graph $G$ (immersed in $\mathbb{R}^2$ for visualization) and its spectrum. Each individual cluster in the graph is $G_j,\, j=1,2,\cdots,12$.} \vspace{-1em} \label{fig:disjoint-graph}
\end{wrapfigure}
As an illustration, we consider the graph, $G$, shown in Figure~\ref{fig:disjoint-graph} with $12$ disjoint components, thus $q=12$.
The graph is generated with $n = 333$ vertices clustered into $12$ components in a randomized manner, with only intra-cluster edges. The weight on every edge is chosen to be $1$. Figure~\ref{fig:disjoint-graph-graph} shows an immersion of the graph in $\mathbb{R}^2$ just for he purpose of visualization (the exact coordinates of the vertices has no significance).

\begin{figure}[h] \vspace{-1em}
\centering
\begin{tabular}{cc}
\hspace{-1.5em}
	\begin{tabular}{c}
	\subfloat[Graph $\widetilde{G}$ created by adding ranomized inter-cluster edges in $G$.]{\includegraphics[width=0.26\textwidth, clip=true, trim=95 210 75 225]
		{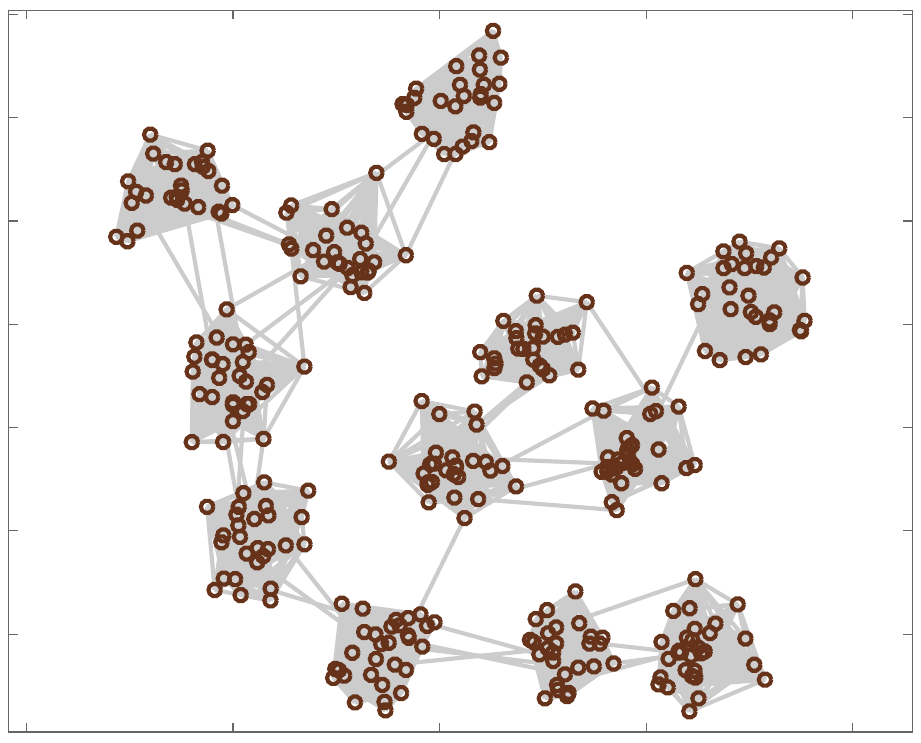} \label{fig:connected-graph-graph}} \\ 
	\subfloat[The spectrum of $\widetilde{L}$. Note the low values of the first $12$ eigenvalues.]{\includegraphics[width=0.28\textwidth, clip=true, trim=75 200 70 225]
		{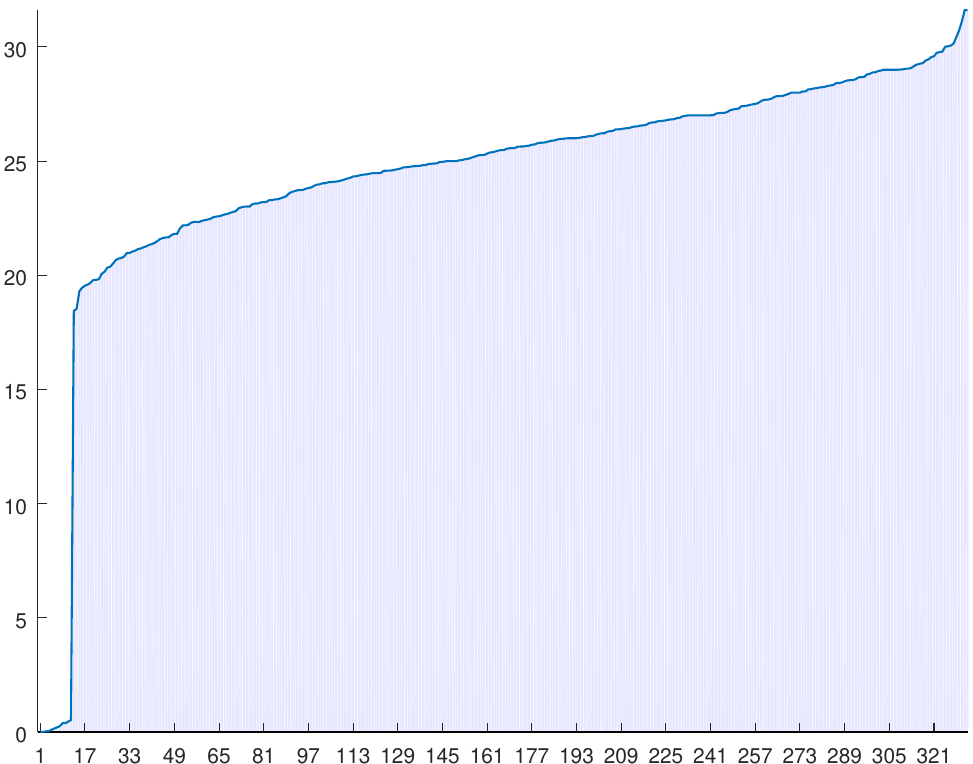} \label{}}
	\end{tabular}
\hspace{-1.5em} &
	\begin{tabular}{cccc}
	\subfloat[$\widetilde{\mathbf{u}}_1$.]{\includegraphics[width=0.165\textwidth, clip=true, trim=95 210 75 225]
		{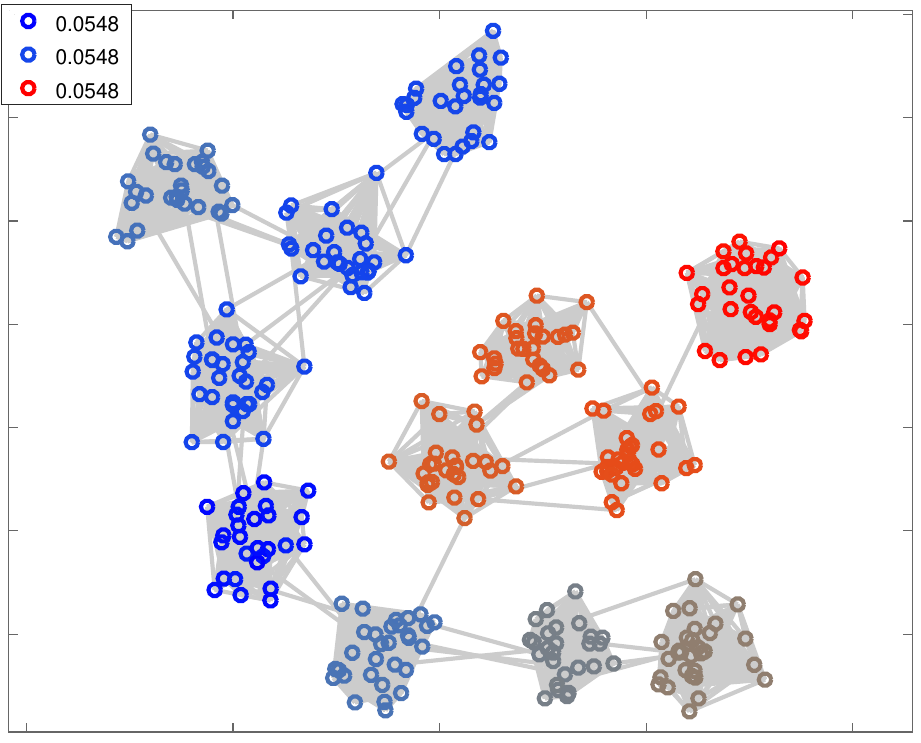} \label{}} \hspace{-1em} &
	\subfloat[$\widetilde{\mathbf{u}}_2$.]{\includegraphics[width=0.165\textwidth, clip=true, trim=95 220 75 215]
		{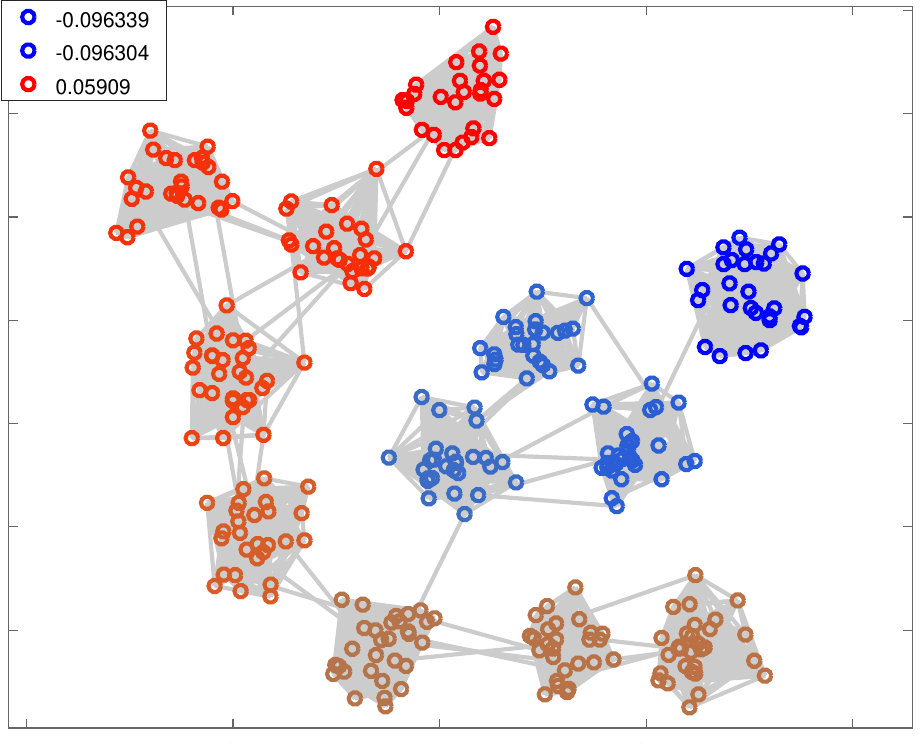} \label{}} \hspace{-1em} &
	\subfloat[$\widetilde{\mathbf{u}}_3$.]{\includegraphics[width=0.165\textwidth, clip=true, trim=95 220 75 215]
		{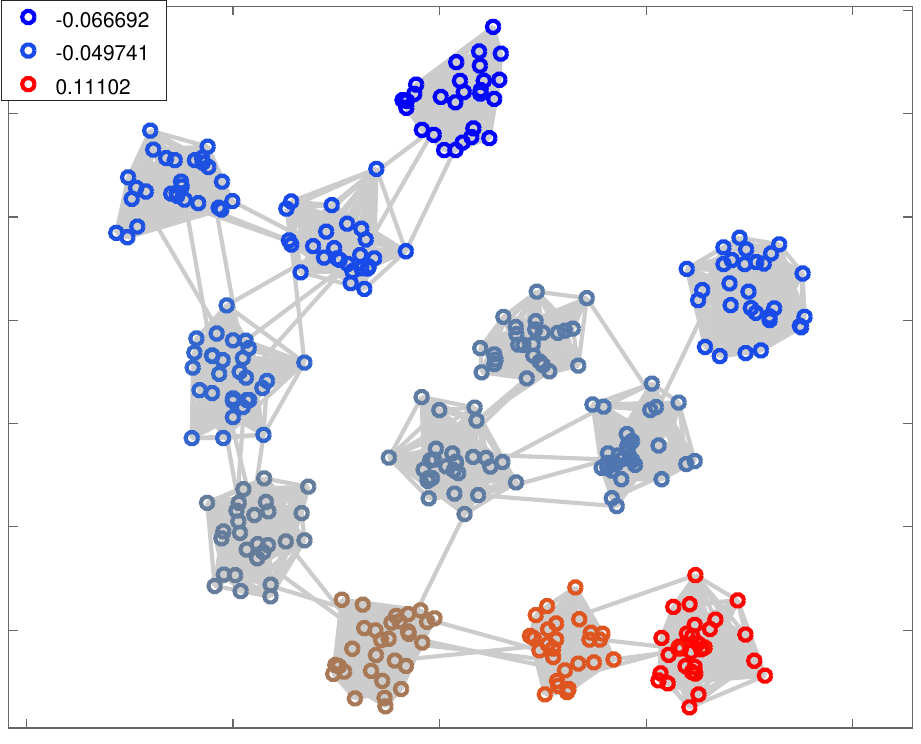} \label{}} \hspace{-1em} &
	\subfloat[$\widetilde{\mathbf{u}}_4$.]{\includegraphics[width=0.165\textwidth, clip=true, trim=95 210 75 225]
		{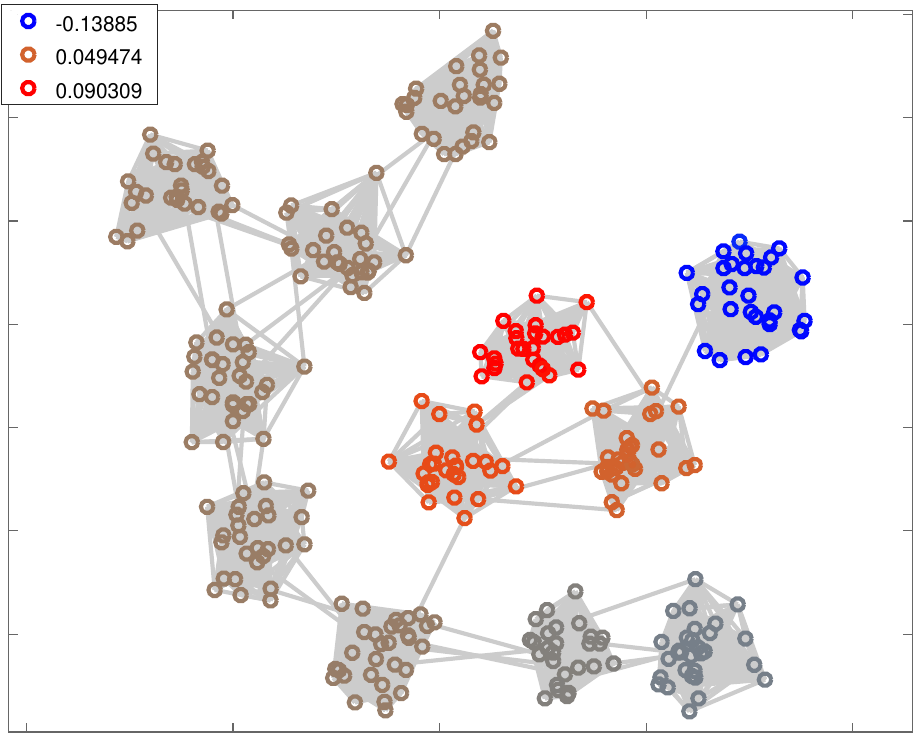} \label{}} \\
	\subfloat[$\widetilde{\mathbf{u}}_5$.]{\includegraphics[width=0.165\textwidth, clip=true, trim=95 220 75 215]
		{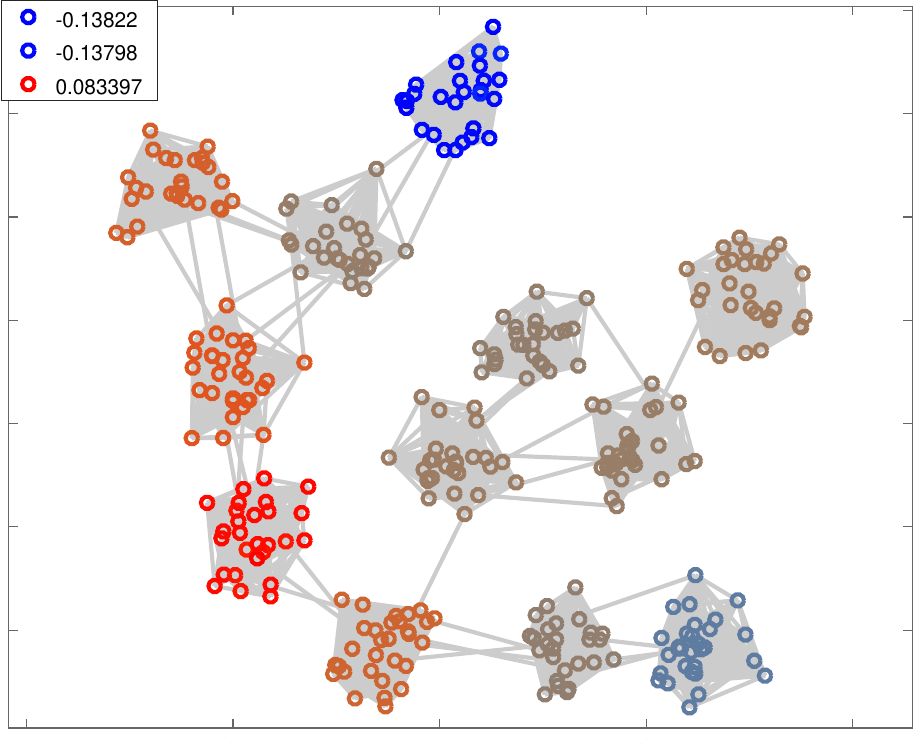} \label{}} \hspace{-1em} &
	\subfloat[$\widetilde{\mathbf{u}}_6$.]{\includegraphics[width=0.165\textwidth, clip=true, trim=95 220 75 215]
		{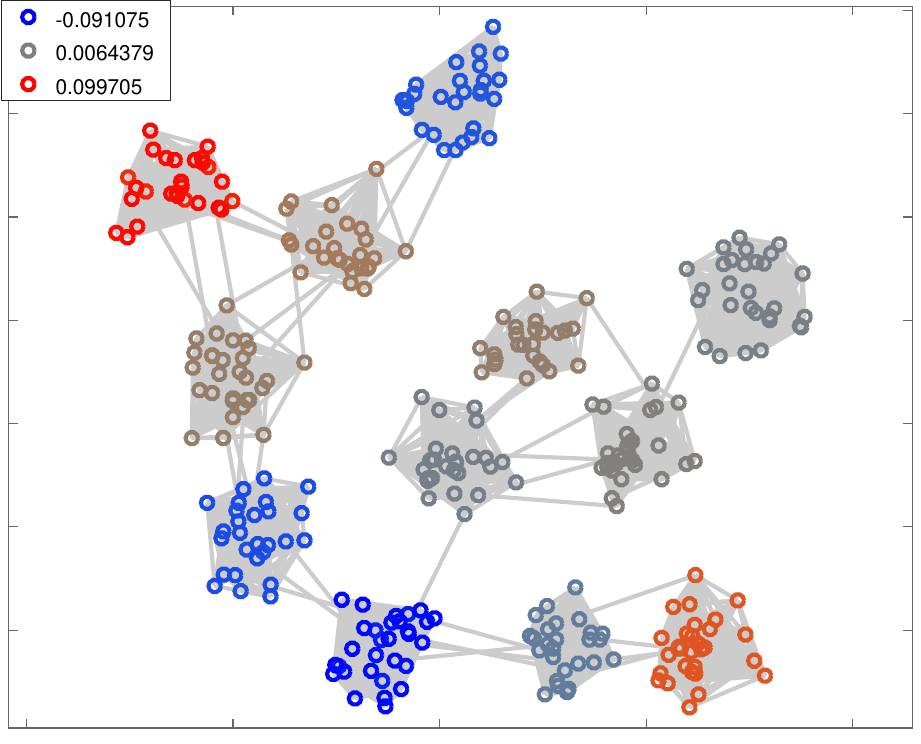} \label{}} \hspace{-1em} &
	\subfloat[$\widetilde{\mathbf{u}}_7$.]{\includegraphics[width=0.165\textwidth, clip=true, trim=95 220 75 215]
		{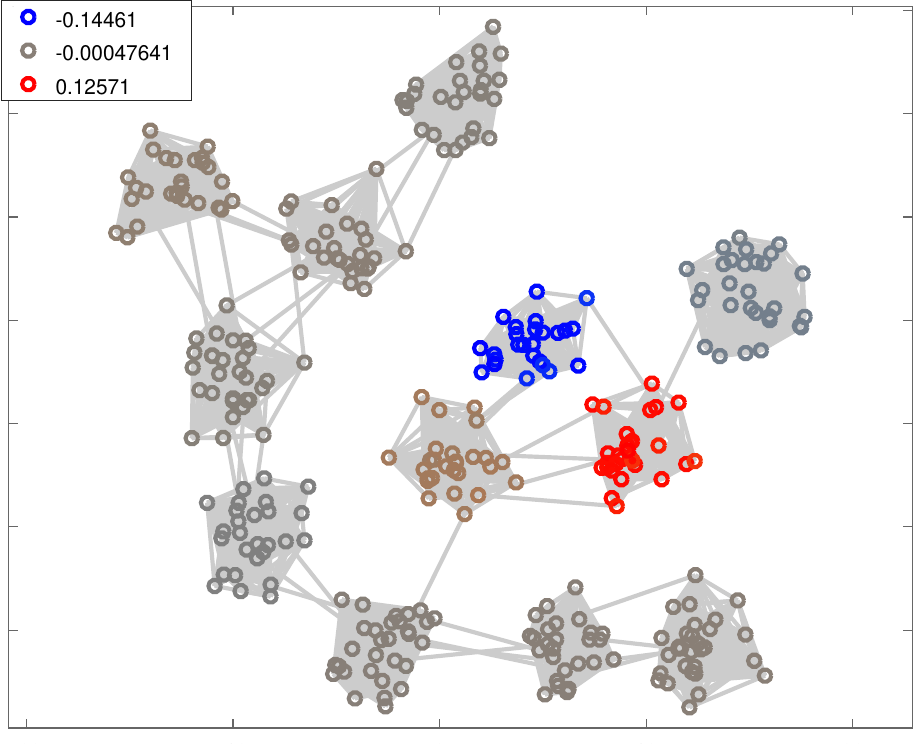} \label{}} \hspace{-1em} &
	\subfloat[$\widetilde{\mathbf{u}}_8$.]{\includegraphics[width=0.165\textwidth, clip=true, trim=95 210 75 225]
		{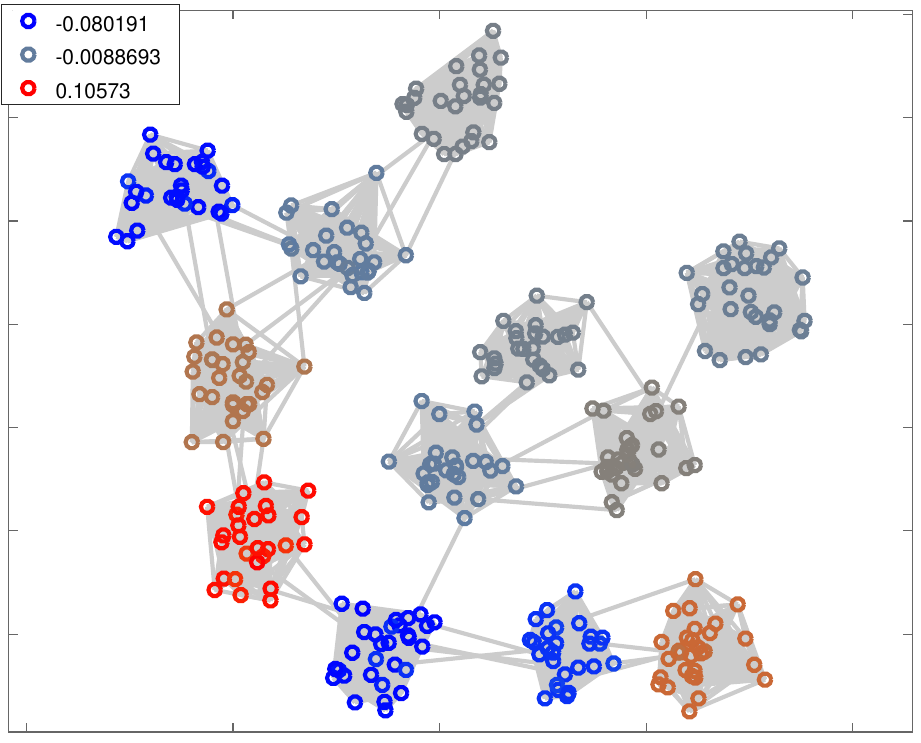} \label{}} \\
	\subfloat[$\widetilde{\mathbf{u}}_9$.]{\includegraphics[width=0.165\textwidth, clip=true, trim=95 220 75 215]
		{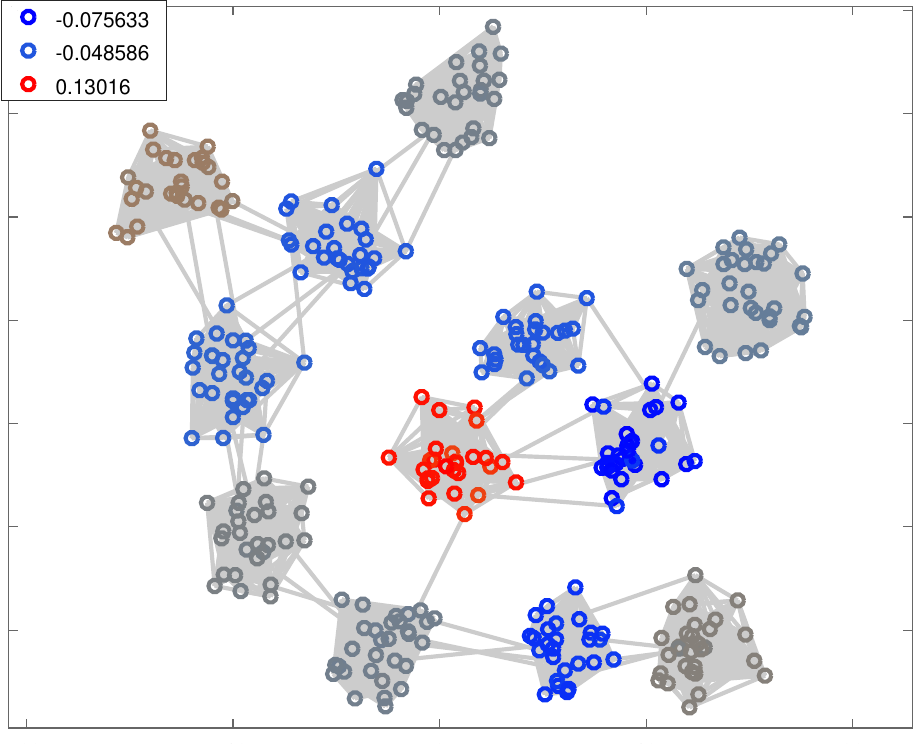} \label{}} \hspace{-1em} &
	\subfloat[$\widetilde{\mathbf{u}}_{10}$.]{\includegraphics[width=0.165\textwidth, clip=true, trim=95 220 75 215]
		{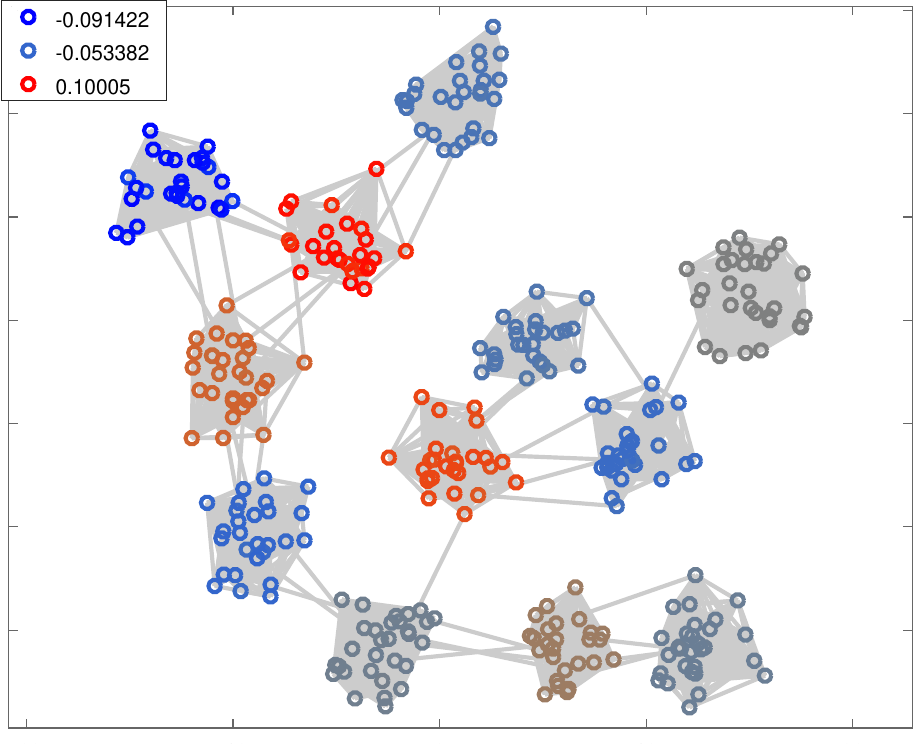} \label{}} \hspace{-1em} &
	\subfloat[$\widetilde{\mathbf{u}}_{11}$.]{\includegraphics[width=0.165\textwidth, clip=true, trim=95 220 75 215]
		{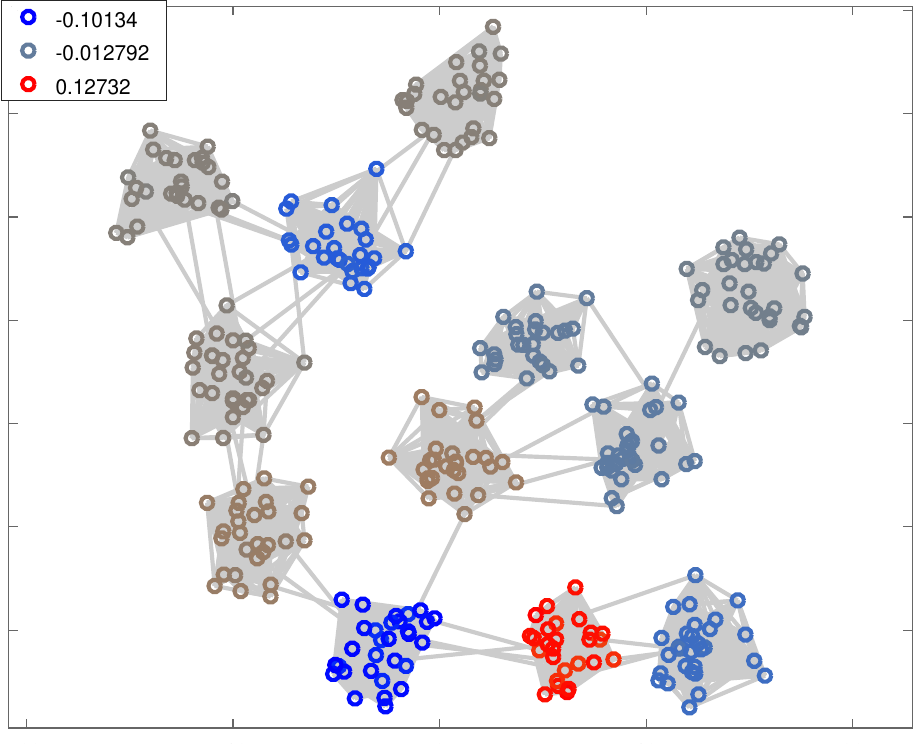} \label{}} \hspace{-1em} &
	\subfloat[$\widetilde{\mathbf{u}}_{12}$.]{\includegraphics[width=0.165\textwidth, clip=true, trim=95 210 75 225]
		{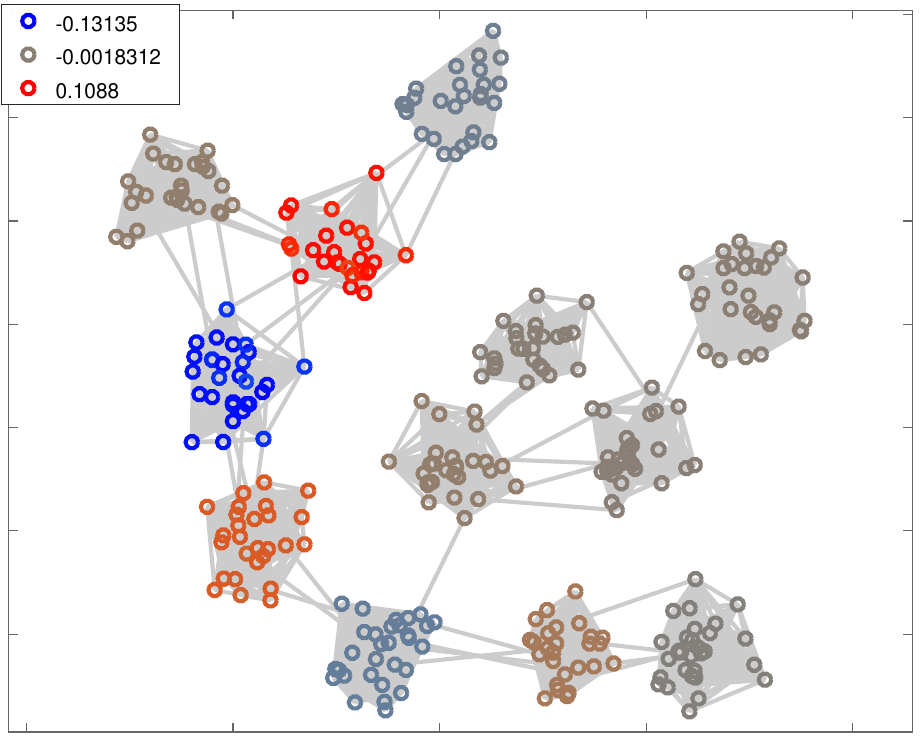} \label{}}
	\end{tabular}
\end{tabular}
\caption{The graph $\widetilde{G}$, the spectrum of its Laplacian, and its first $12$ eigenvectors (c-n) visualized as distribution over the vertices (red is positive, blue is negative).}
\label{fig:connected-graph}
\end{figure}

We then construct $\widetilde{G}$ by establishing randomized edges between the components of $G$. The weight on every inter-cluster edge is also chosen to be $1$. Figure~\ref{fig:connected-graph-graph} shows the immersion of the resultant graph.


Direct computation reveals that for these graphs, 
$\widetilde{\lambda}_{q+1} = 18.436$ and 
$\frac{1}{q} \sum_{j=1}^q \couplingtext_{\widetilde{G}}(G_j) = 0.5417$.
%
The L.H.S. of \eqref{eq:graph-subspace-inequality-kappa-zero} is 
$d_\mathrm{\subspacetext} \left( \spn(\mathbf{u}_{\{1,2,\cdots,12\}}) , \spn(\widetilde{\mathbf{u}}_{\{1,2,\cdots,12\}}) \right) = 2.516\times 10^{-2}$, 
while the R.H.S. is
$\frac{ 
	\sqrt{\frac{1}{q} \sum_{j=1}^q \couplingtext_{\widetilde{G}}(G_j)} 
}{ 
	\widetilde{\lambda}_{q+1}
} = 3.992\times 10^{-2}$, thus validating the result of Proposition~\ref{prop:general-coupling-result}.

Again, $\max_{j\in\{1,\cdots,q\}} \maxextdeg_{\widetilde{G}}(G_j) = 3$ and $ \frac{\lambda_{q+1}}{4} = 4.6091$,
thus satisfying the condition for Proposition~\ref{prop:null-space-purturbation-known-L-spectrum}.
The 
R.H.S. in \eqref{eq:null-space-purturbation-known-L-spectrum-1} is 
$\frac{ 
 \sqrt{ \frac{1}{q}~ \sum_{j=1}^q 
	\couplingtext_{\widetilde{G}}(G_j) }
}{ 
	\lambda_{q+1} ~-~ 
	2 \, \max_{k\in\{1,\cdots,q\}} \maxextdeg_{\widetilde{G}}(G_k) 
} = 6.036 \times 10^{-2}$, thus validating the result of the proposition.

\begin{figure}[h] \vspace{-1em}
	\centering
\begin{tabular}{cccc}
	\subfloat[${\mathbf{u}}'_1$.]{\includegraphics[width=0.165\textwidth, clip=true, trim=95 210 75 225]
		{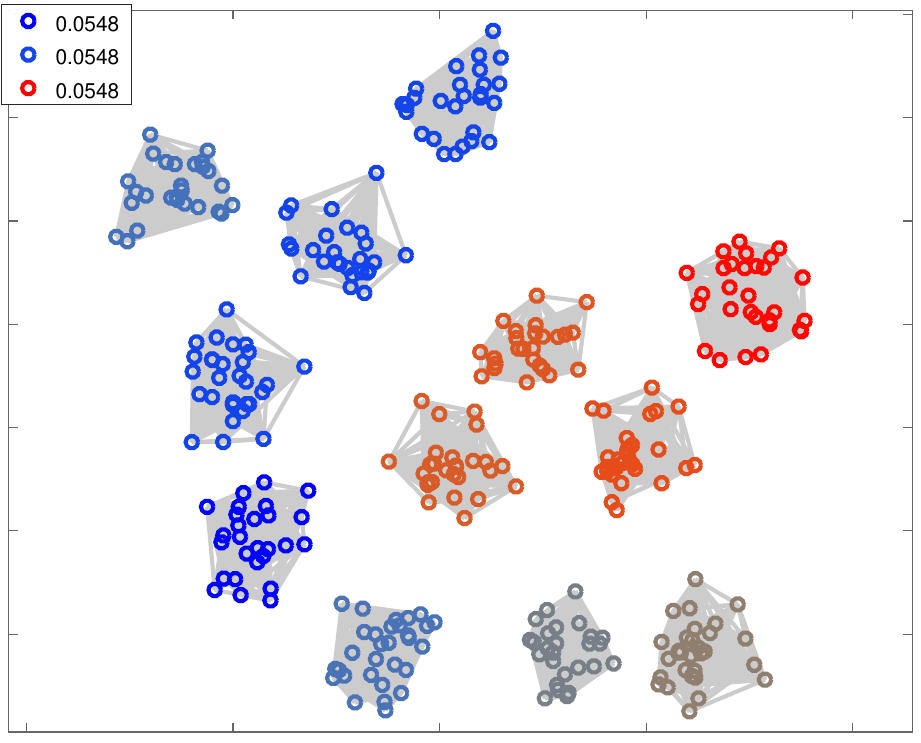} \label{}} \hspace{-0em} &
	\subfloat[${\mathbf{u}}'_2$.]{\includegraphics[width=0.165\textwidth, clip=true, trim=95 210 75 225]
		{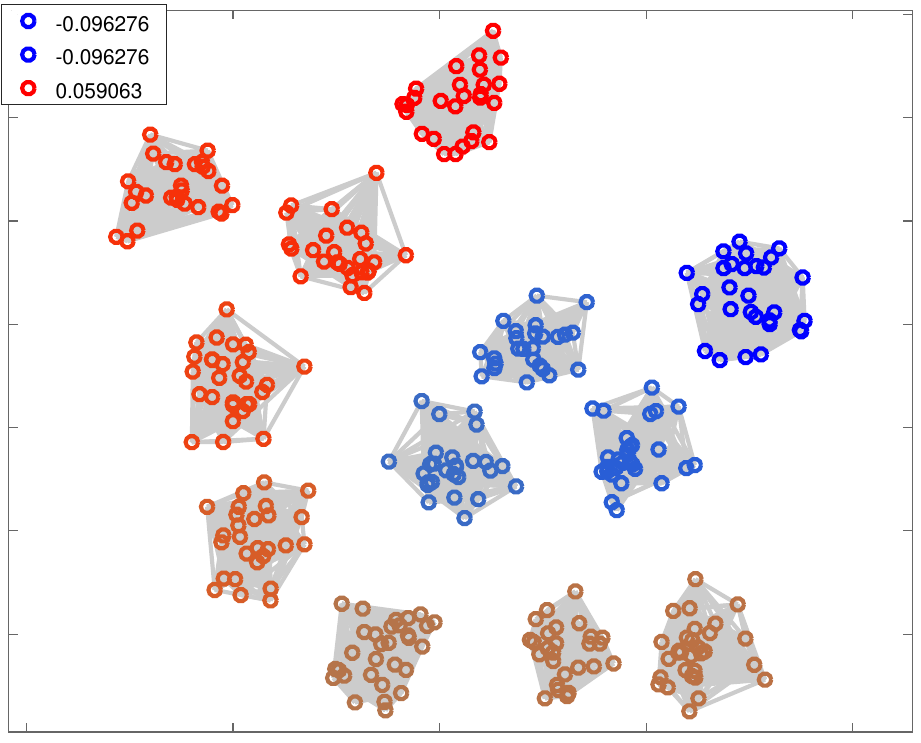} \label{}} \hspace{-0em} &
	\subfloat[${\mathbf{u}}'_3$.]{\includegraphics[width=0.165\textwidth, clip=true, trim=95 220 75 215]
		{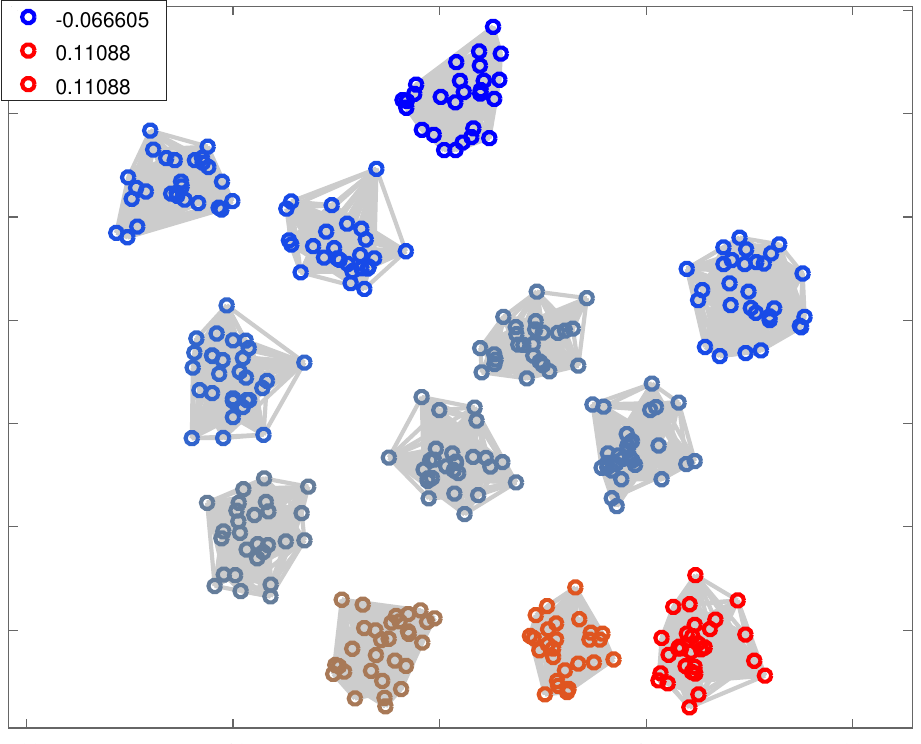} \label{}} \hspace{-0em} &
	\subfloat[${\mathbf{u}}'_4$.]{\includegraphics[width=0.165\textwidth, clip=true, trim=95 210 75 225]
		{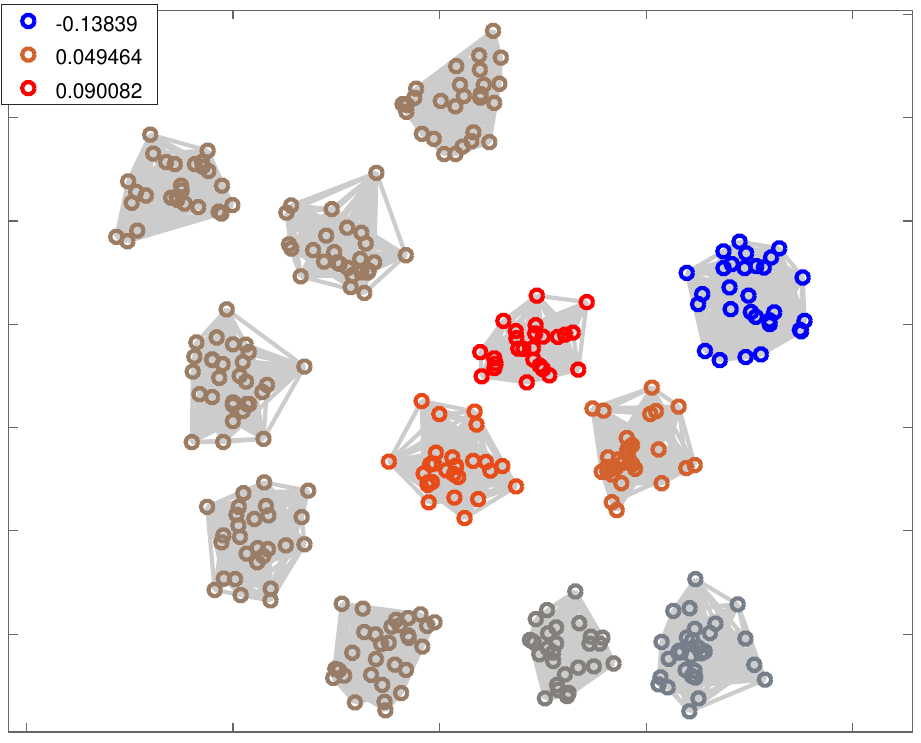} \label{}} \\
	\subfloat[${\mathbf{u}}'_5$.]{\includegraphics[width=0.165\textwidth, clip=true, trim=95 220 75 215]
		{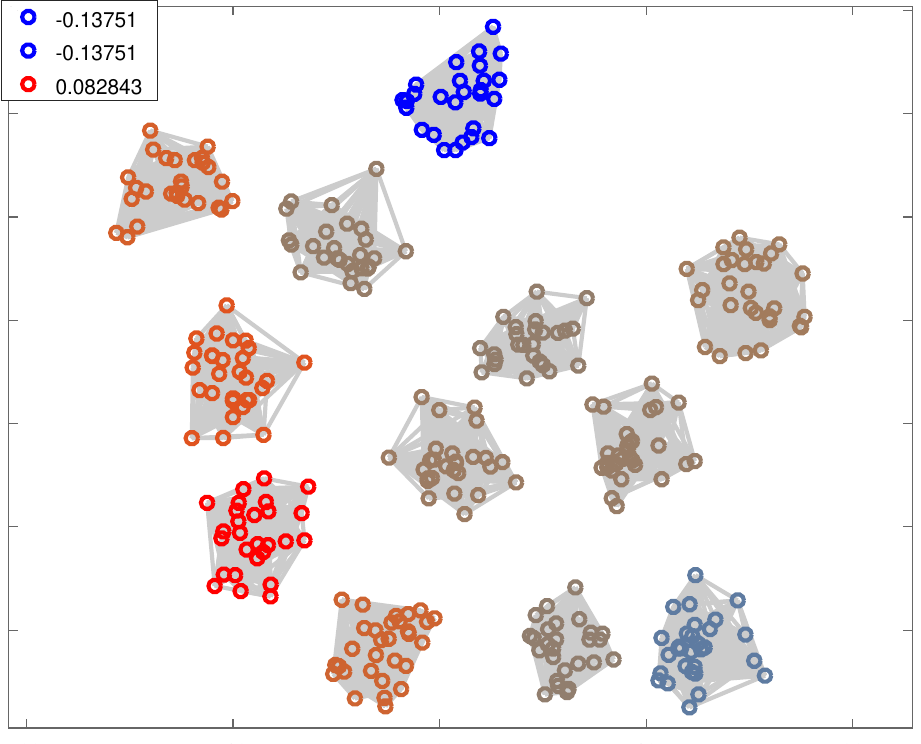} \label{}} \hspace{-0em} &
	\subfloat[${\mathbf{u}}'_6$.]{\includegraphics[width=0.165\textwidth, clip=true, trim=95 210 75 225]
		{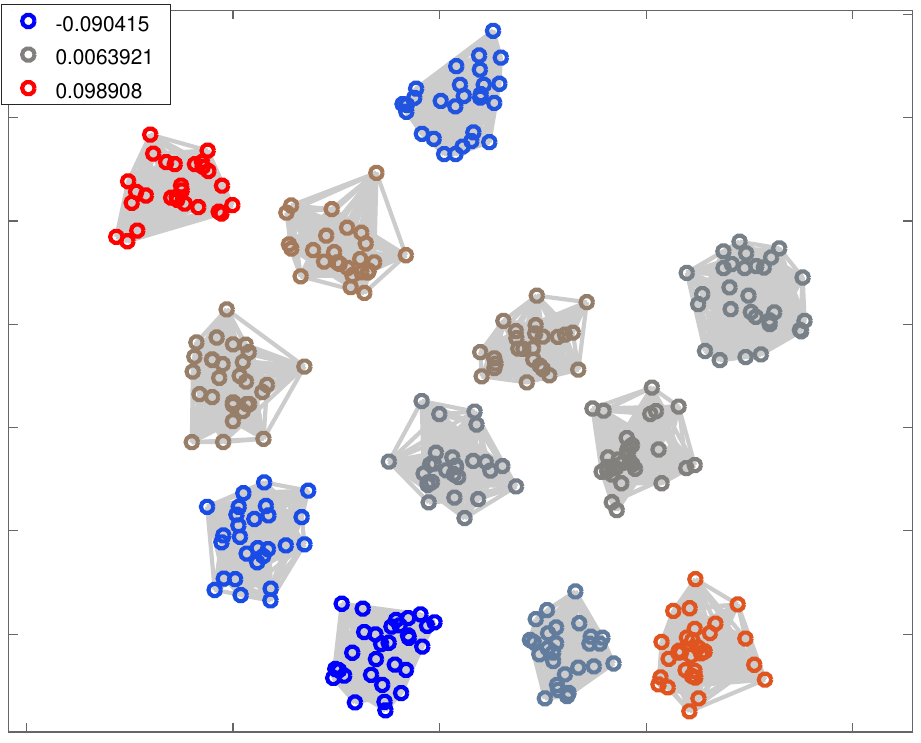} \label{}} \hspace{-0em} &
	\subfloat[${\mathbf{u}}'_7$.]{\includegraphics[width=0.165\textwidth, clip=true, trim=95 210 75 225]
		{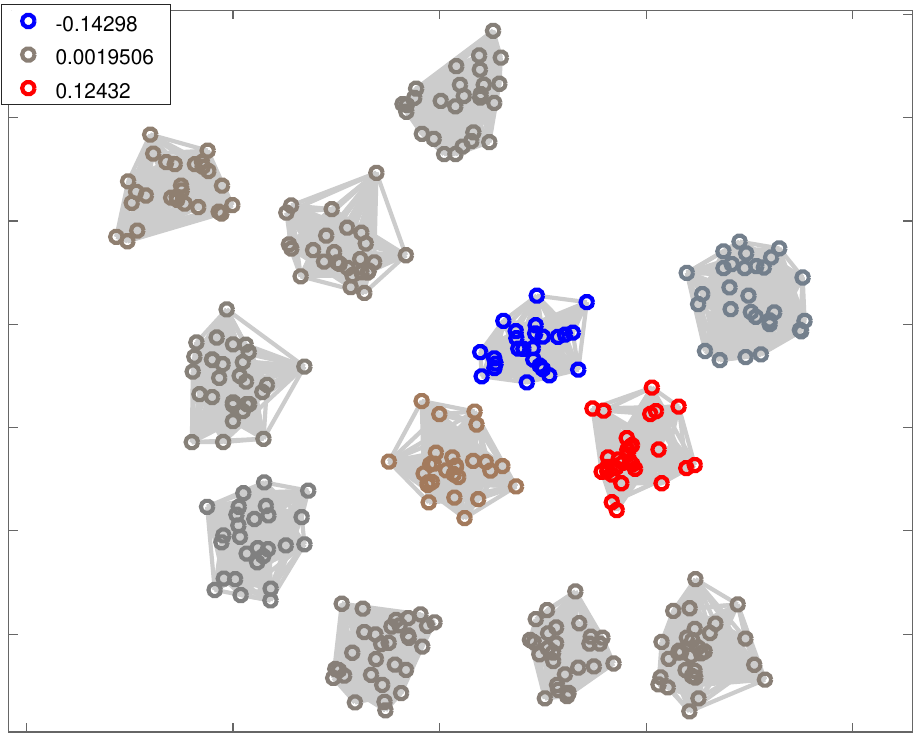} \label{}} \hspace{-0em} &
	\subfloat[${\mathbf{u}}'_8$.]{\includegraphics[width=0.165\textwidth, clip=true, trim=95 210 75 225]
		{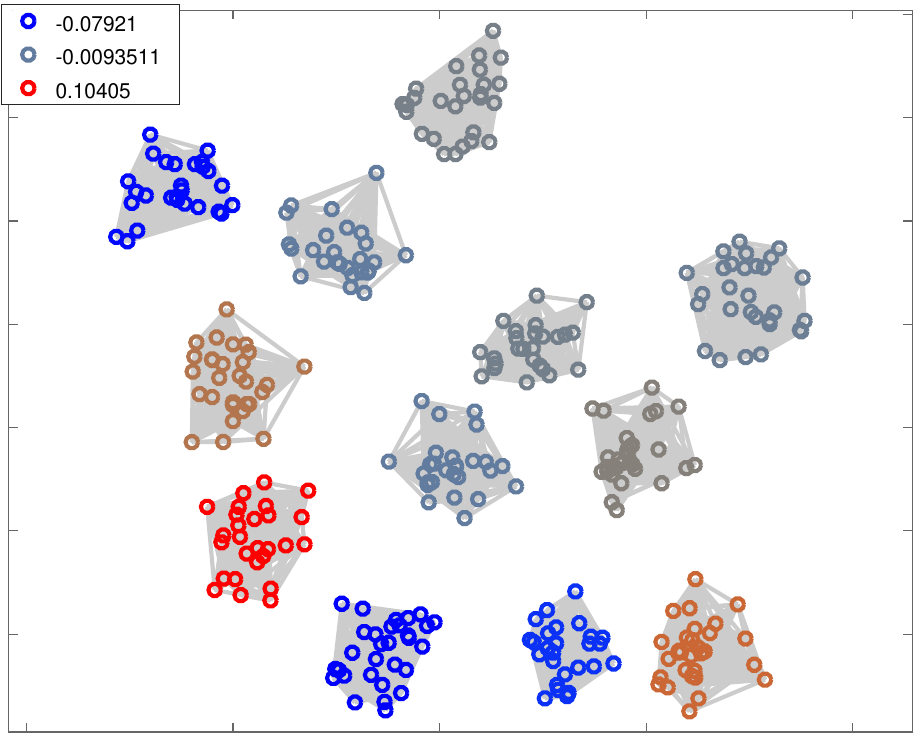} \label{}} \\
	\subfloat[${\mathbf{u}}'_9$.]{\includegraphics[width=0.165\textwidth, clip=true, trim=95 210 75 225]
		{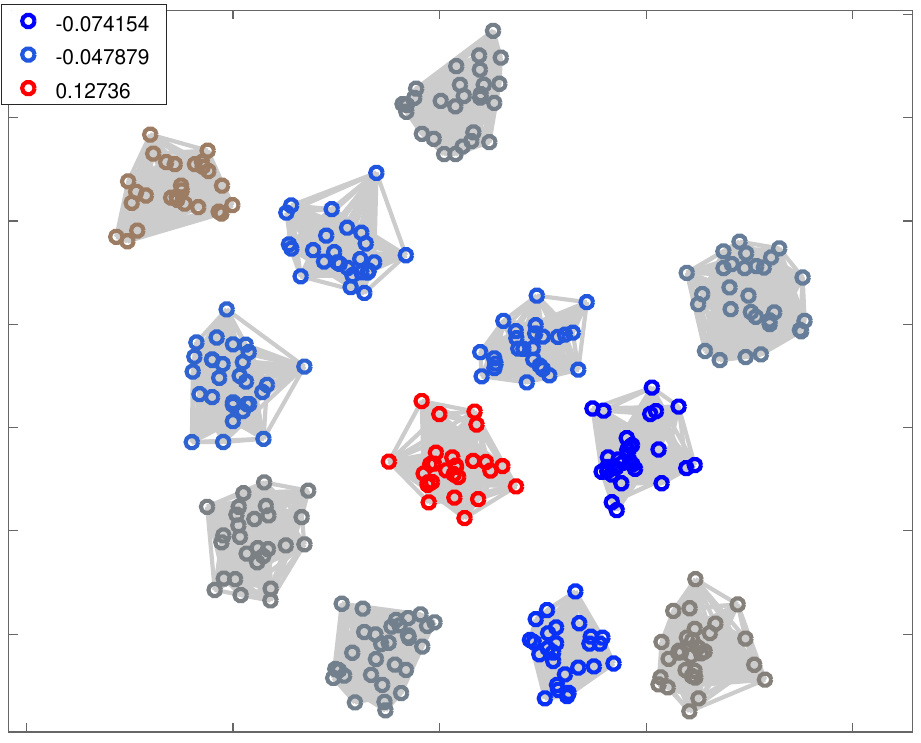} \label{}} \hspace{-0em} &
	\subfloat[${\mathbf{u}}'_{10}$.]{\includegraphics[width=0.165\textwidth, clip=true, trim=95 210 75 225]
		{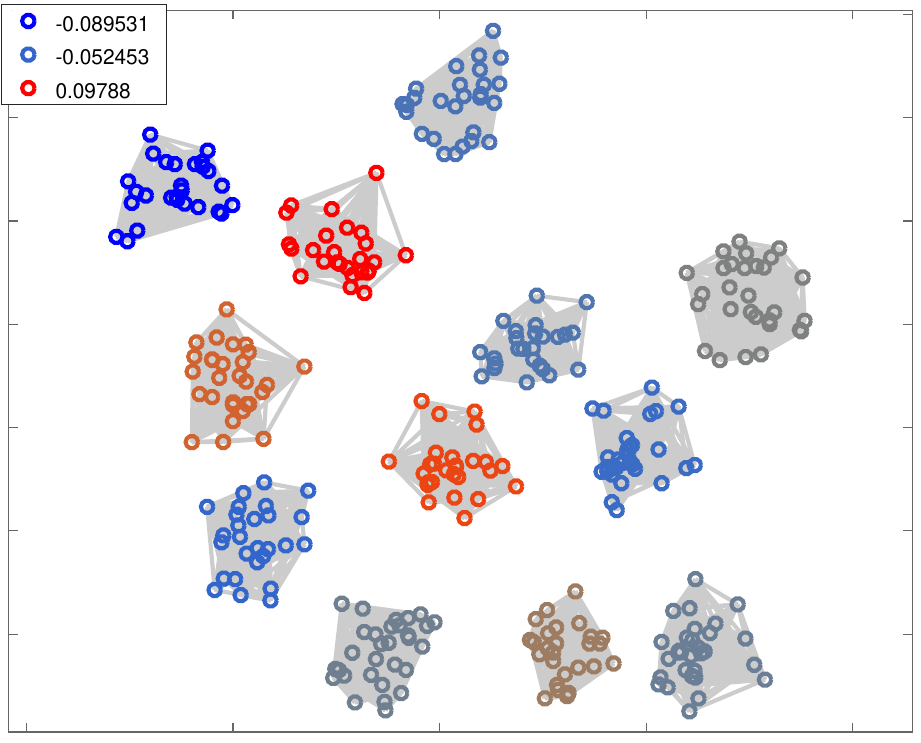} \label{}} \hspace{-0em} &
	\subfloat[${\mathbf{u}}'_{11}$.]{\includegraphics[width=0.165\textwidth, clip=true, trim=95 220 75 215]
		{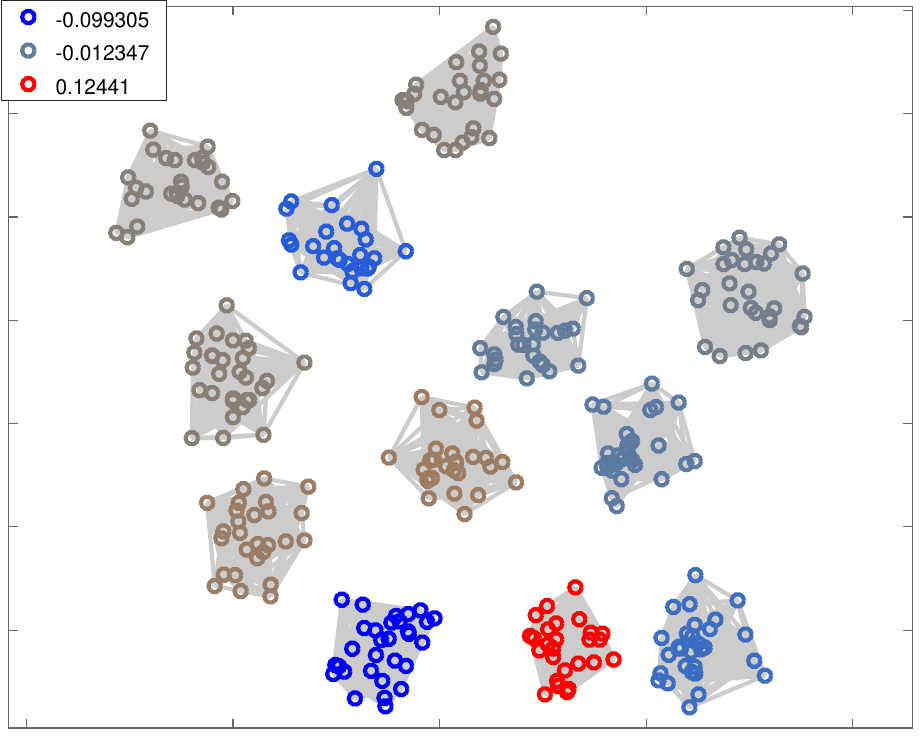} \label{}} \hspace{-0em} &
	\subfloat[${\mathbf{u}}'_{12}$.]{\includegraphics[width=0.165\textwidth, clip=true, trim=95 210 75 225]
		{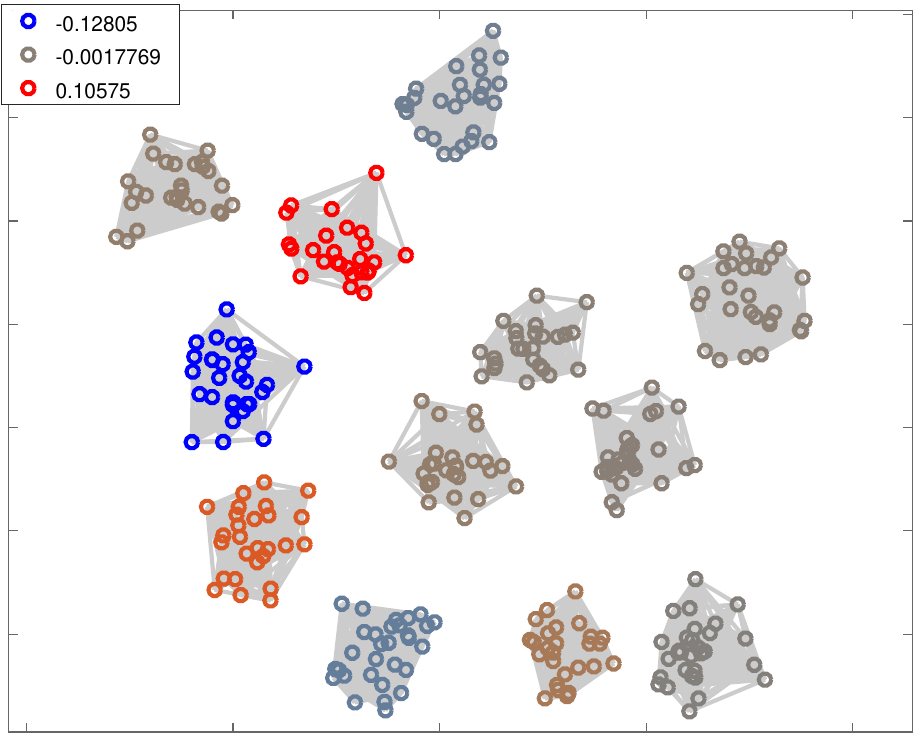} \label{}}
\end{tabular}
\caption{The basis $\{\mathbf{u}'_j\}_{j=1,2,\cdots,q}$ of the null-space of $L$ visualized as distributions over the vertices (red is positive, blue is negative). Compare this with Figures~\ref{fig:connected-graph}(c-n).}
\label{fig:disjoint-graph-modes}
\end{figure}

Since the chosen basis, $\{\mathbf{u}_j\}_{j=1,2,\cdots,q}$, for the null-space of $L$ constitutes of distributions such that $\mathbf{u}_j$ is uniform and positive over vertices of $G_j$ and zero everywhere else, this basis is not ideal for a visual comparison with $\{\widetilde{\mathbf{u}}_j\}_{j=1,2,\cdots,q}$.
For a visual comparison between $\spn(\mathbf{u}_{\{1,2,\cdots,q\}})$ and $\spn(\widetilde{\mathbf{u}}_{\{1,2,\cdots,q\}})$, we choose a basis for the null-space of $L$ that is closest to $\{\widetilde{\mathbf{u}}_j\}_{j=1,2,\cdots,q}$:
Define the $q\times q$ matrix $R = \left([\mathbf{u}_1, \mathbf{u}_2, \cdots, \mathbf{u}_q]\right)^{+} [\widetilde{\mathbf{u}}_1, \widetilde{\mathbf{u}}_2, \cdots, \widetilde{\mathbf{u}}_q]$, where $(\cdot)^{+}$ indicates the Moore-Pesrose pseudoinverse.
We need to chose an unitary matrix that is close to $R$. This is given by taking the SVD of $R = V \Sigma W^\conjsym$ and defining $R' = V W^\conjsym$. Then a basis for $\spn(\mathbf{u}_{\{1,2,\cdots,q\}})$ is defined 
by the columns of 
$[\mathbf{u}_1, \mathbf{u}_2, \cdots, \mathbf{u}_q] R' =: [\mathbf{u}'_1, \mathbf{u}'_2, \cdots, \mathbf{u}'_q]$.
Figure~\ref{fig:disjoint-graph-modes} shows these vectors as distributions over the vertices of $G$.

\section*{References} 
\vspace{0.5em} 
\bibliographystyle{plain}
\bibliography{references}

\end{document}